\input amssym.def
\input amssym
\input epsf.tex
\input epsfig.sty
\documentstyle[11pt]{article}
\begin{document}
\bibliographystyle{alpha}

\newtheorem{thm}{Theorem}
\newtheorem{defin}[thm]{Definition}
\newtheorem{lemma}[thm]{Lemma}
\newtheorem{propo}[thm]{Proposition}
\newtheorem{cor}[thm]{Corollary}
\newtheorem{conj}[thm]{Conjecture}
\newtheorem{exa}[thm]{Example}

\centerline{\LARGE \bf Abelian Categories of}

\centerline{\LARGE \bf Modules over a (Lax) Monoidal Functor}
\vspace*{1cm}

\centerline{\parbox{2.5in}{D.\ N.\ Yetter \\ Department of Mathematics \\
Kansas State University \\ Manhattan, KS 66506}}
\vspace*{1cm}

{\small
\noindent{\bf Abstract:} In \cite{CY.def} Crane and Yetter introduced a
deformation theory for monoidal categories.  The related deformation theory
for monoidal functors introduced by Yetter in \cite{Y.def} is a proper
generalization of Gerstenhaber's deformation theory for associative algebras
\cite{G1,G2,GS}.  In the present paper we solidify the analogy between lax
monoidal functors and associative algebras by showing that under suitable
conditions, categories of functors with an action of a lax monoidal
functor are abelian categories.  The deformation complex of a 
monoidal functor is
generalized to an analogue of the Hochschild complex with coefficients
in a bimodule, and the deformation complex of a monoidal natural transformation
is shown to be a special case.  It is shown further that the cohomology of a 
monoidal functor $F$ with coefficients in an $F,F$-bimodule is given by 
right derived functors.}

\clearpage 
\vspace*{2cm} \section{Introduction \label{intro}}
\vspace*{1cm}

In \cite{CY.def} Crane and 
Yetter introduced an infinitesimal deformation theory
for monoidal categories.  Yetter \cite{Y.def} introduced related deformation
theories for monoidal functors and braided monoidal categories.  
In \cite{Y.def,Y.book}
Yetter showed that the rich structure present in Gerstenhaber's deformation
theory of associative algebras is again present in these new theories and
gave the hitherto missing proof that all obstructions are closed.  
For monoidal functors, as for associative algebras, first order deformations
are classified by $H^2$, while obstructions to higher order deformations
lie in $H^3$.  For monoidal categories, the first order deformations lie
in $H^3$, while obstructions lie in $H^4$.

One problem with categorical deformation theories not found in algebraic
deformation theories is the fact that the natural definition of the
cochain groups as module of natural transformations does not live
comfortably in a natural setting for a calculus of derived functors:  the
source category for the functors between which natural transformations
are taken changes from degree to degree.  The structures and results
of this paper arose from the need to remedy this apparent defect.

In sections 2 through 6 we follow closely \cite{Y.book}, except that we
treat strong monoidal functors as special cases of lax rather than oplax
monoidal functors.

\clearpage
  \section{Review of Categorical Deformations \label{def.rev}}
\vspace*{1cm}

Throughout, we work in a setting in which all categories are $R$-linear 
for some
fixed commutative ring 
$R$ and in which all functors are linear (in each variable).  As was
observed in \cite{Y.FTTD,Y.book}, given an $R$-linear category $\cal C$
and an $R$-algebra $A$, we can
form a category ${\cal C}\otimes A$ by {\em extension of scalars}:

\[ Ob({\cal C}\otimes A) = Ob({\cal C}) \]

\noindent and

\[ Hom_{{\cal C}\otimes A}(X,Y) = Hom_{\cal C}(X,Y)\otimes_R A \]
 
\noindent extending both composition and any structural functors
by (multi)linearity.  If the algebra is an $m$-adically complete local
ring, we can also form an $m$-adic completion $\widehat{{\cal C}\otimes A}$,
by $m$-adically completing all of the homsets and extending composition
and any structural functors by continuity.

As in \cite{Y.FTTD,Y.book}, we 
denote ${\cal C}\otimes R[\epsilon]/\!<\!\epsilon^{n+1}\!>\!$ by
${\cal C}^{(n)}$, and $\widehat{{\cal C}\otimes R[[x]]}$ by 
${\cal C}^{(\infty)}$.

Similarly, if we have a functor $F:{\cal C}\rightarrow {\cal D}$, we may
obtain a functor from ${\cal C}\otimes A$ to ${\cal D}\otimes A$ by
extending the definition of $F$ on maps by linearity, likewise in the
$m$-adically complete setting.  In the case of $A =  
R[\epsilon]/\!<\!\epsilon^{n+1}\!>\!$ or $A = R[[x]]$, we will denote the 
resulting functor by $F^{(n)}$ or $F^{(\infty)}$, respectively.

We will always work in a context where we have a specifed category 
$\cal C$, a specified functor $F:{\cal C}\rightarrow {\cal D}$, or a 
specified natural transformation 
$\theta:F\Rightarrow G$ for $F,G:{\cal C}\rightarrow {\cal D}$. 
The category(ies) will be equipped with some structure functors
and structural natural transformations between them and satisfying
specified coherence conditions. The functor(s) will be equipped with
some
structural natural transformations 
relating them and the structure
functors of $\cal C$ and $\cal D$ and specified 
coherence conditions relating these
with the structural natural transformations of $\cal C$ and $\cal D$.
The natural transformation $\theta$ will be equipped with 
specified coherence conditions relating
$\theta$ and the structural natural transformations of $F,G,{\cal C}$ and
$\cal D$.  We call such a category (resp. functor, natural transformation)
a {\em category with structure} (resp. {\em functor with structure},
{\em natural transformation with structure}). 

In the case of a natural transformation with structure, we consider the
natural transformation itself to be a ``structural natural transformation''.
We then make:

\begin{defin}
An {\em $n^{th}$ order deformation} of a category with structure $\cal C$ 
(resp. a functor with structure $F:{\cal C}\rightarrow {\cal D}$, 
a natural transformation with structure $\theta:F\Rightarrow G$)
is an assignment to each structural natural transformation 
$\phi:\Phi \Rightarrow \Psi$ of $\cal C$
(resp. of $\cal C$, $\cal D$, and $F$) of a structural natural
transformation 

\[ \phi^{(n)} = \phi + \phi_1 \epsilon + \ldots + \phi_n \epsilon^n \]

\noindent such that $\phi^{(n)}:\Phi^{(n)}\Rightarrow \Psi^{(n)}$ and 
${\cal C}^{(n)}$ (resp. $F^{(n)}$, $\theta^{(n)}$) is a category with
structure (resp. functor with structure, natural transformation with
structure) satisfying the same coherence conditions as $\cal C$ (resp.
$F$, $\theta$).

A {\em formal deformation} is defined similarly by a formal power series.
\end{defin}

Among these we can distinguish one which always exists:

\begin{defin}
The {\em trivial deformation} of $\cal C$ (resp. 
$F:{\cal C}\rightarrow {\cal D}$, $\theta:F\Rightarrow G$) is given by 
letting $\phi^{(n)} = \phi$
for all structural natural transformations.
\end{defin}

In terms of this strict notion of triviality, we can distinguish important 
special cases of deformations of a functor or natural transformation.  
Observe that a deformation of a functor with structure (resp. natural
transformation with structure) induces ``forgetful'' deformations of
its source and target (resp. its source and target and their source and
target).  We can thus consider cases in which certain of the induced
``forgetful'' deformations are trivial:

\begin{defin}
A deformation of a functor with structure 
$F:{\cal C}\rightarrow {\cal D}$ between two
categories with structure is {\em purely functorial} (resp. {\em fibred},
{\em cofibred}) if the induced deformations on the source and target
(resp. target, source) are trivial.

A deformation of a natural transformation with structure
$\theta:F\Rightarrow G$ between two functors with structure is
{\em purely transformational} (resp. {\em 1-fibred}, {\em 1-cofibred},
{\em purely functorial}, {\em 2-fibred}, {\em 2-cofibred})
if the source and target (resp. target, source, source and target of
the parallel functors, target of the parallel functors, source of the
parallel functors) are trivial.
\end{defin}

In the case to which we will apply this very general notion --- semigroupal
and monoidal categories, and the various types of monoidal and semigroupal
functors --- there is always a good notion of equivalence between deformations.
Let us review the relevant definitions and some theorems and lemmas which
will be needed:

\begin{defin}
A {\em mon\-oid\-al category} \index{category, monoidal} \index{monoidal category}
$\cal C$ is a category $\cal C$ equipped
with a functor $\otimes :{\cal C} \times {\cal C} \rightarrow {\cal C}$
and an object $I$, together with natural isomorphisms 
$\alpha : \otimes (\otimes \times 1_{\cal C}) \Rightarrow \otimes (1_{\cal C}
\times \otimes)$, $\rho : \otimes I \Rightarrow 1_{\cal C}$ and 
$\lambda : I\otimes \Rightarrow 1_{\cal C}$, satisfying the
pentagon and triangle coherence conditions of Figure \ref{mon.cat} and the
bigon ($\rho_I = \lambda_I$) coherence condition \index{coherence condition} 
\index{pentagon}
(cf. \cite{CWM}).\footnote{It can be shown that the bigon condition is 
redundant.}  Similarly, a {\em semigroupal
category}\index{semigroupal category}\index{category, semigroupal} 
is a category equipped with only $\otimes$ and $\alpha$, satisfying
the pentagon of Figure \ref{mon.cat}.  
\end{defin}
 
\begin{figure}[htb] \centering

\setlength{\unitlength}{3947sp}%
\begingroup\makeatletter\ifx\SetFigFont\undefined
\def\x#1#2#3#4#5#6#7\relax{\def\x{#1#2#3#4#5#6}}%
\expandafter\x\fmtname xxxxxx\relax \def\y{splain}%
\ifx\x\y   
\gdef\SetFigFont#1#2#3{%
  \ifnum #1<17\tiny\else \ifnum #1<20\small\else
  \ifnum #1<24\normalsize\else \ifnum #1<29\large\else
  \ifnum #1<34\Large\else \ifnum #1<41\LARGE\else
     \huge\fi\fi\fi\fi\fi\fi
  \csname #3\endcsname}%
\else
\gdef\SetFigFont#1#2#3{\begingroup
  \count@#1\relax \ifnum 25<\count@\count@25\fi
  \def\x{\endgroup\@setsize\SetFigFont{#2pt}}%
  \expandafter\x
    \csname \romannumeral\the\count@ pt\expandafter\endcsname
    \csname @\romannumeral\the\count@ pt\endcsname
  \csname #3\endcsname}%
\fi
\fi\endgroup
\begin{picture}(4287,6180)(76,-5386)
\thinlines
\put(1876,-4036){\vector( 1, 0){1125}}
\put(1426,-4186){\vector( 3,-4){675}}
\put(3226,-4261){\vector(-3,-4){675}}
\put(901,-4111){\makebox(0,0)[lb]{\smash{\SetFigFont{10}{14.4}{rm}$(A \otimes I) \otimes B$}}}
\put(3151,-4111){\makebox(0,0)[lb]{\smash{\SetFigFont{10}{14.4}{rm}$A \otimes (I \otimes B)$}}}
\put(2251,-3961){\makebox(0,0)[lb]{\smash{\SetFigFont{10}{14.4}{rm}$\alpha$}}}
\put(1126,-4861){\makebox(0,0)[lb]{\smash{\SetFigFont{10}{14.4}{rm}$\rho \otimes B$}}}
\put(3151,-4861){\makebox(0,0)[lb]{\smash{\SetFigFont{10}{14.4}{rm}$A \otimes \lambda$}}}
\put(2101,-5386){\makebox(0,0)[lb]{\smash{\SetFigFont{10}{14.4}{rm}$A \otimes B$}}}
\put(901,-886){\vector( 1,-2){450}}
\put(2251,-2086){\vector( 1, 0){900}}
\put(3751,-1861){\vector( 2, 3){600}}
\put(3076,464){\vector( 1,-1){900}}
\put(1051,-436){\vector( 1, 1){900}}
\put(601,-2161){\makebox(0,0)[lb]{\smash{\SetFigFont{10}{14.4}{rm}$(A \otimes (B \otimes C)) \otimes D$}}}
\put(1126,-61){\makebox(0,0)[lb]{\smash{\SetFigFont{10}{14.4}{rm}$\alpha$}}}
\put(526,-1411){\makebox(0,0)[lb]{\smash{\SetFigFont{10}{14.4}{rm}$\alpha \otimes D$}}}
\put(2626,-2311){\makebox(0,0)[lb]{\smash{\SetFigFont{10}{14.4}{rm}$\alpha$}}}
\put( 76,-736){\makebox(0,0)[lb]{\smash{\SetFigFont{10}{14.4}{rm}$((A \otimes B) \otimes C) \otimes D$}}}
\put(3526,-736){\makebox(0,0)[lb]{\smash{\SetFigFont{10}{14.4}{rm}$A \otimes (B \otimes (C \otimes D))$}}}
\put(3301,-2161){\makebox(0,0)[lb]{\smash{\SetFigFont{10}{14.4}{rm}$A \otimes ((B \otimes C) \otimes D)$}}}
\put(4276,-1486){\makebox(0,0)[lb]{\smash{\SetFigFont{10}{14.4}{rm}$A \otimes \alpha$}}}
\put(1876,614){\makebox(0,0)[lb]{\smash{\SetFigFont{10}{14.4}{rm}$(A \otimes B) \otimes (C \otimes D)$}}}
\put(3901, 14){\makebox(0,0)[lb]{\smash{\SetFigFont{10}{14.4}{rm}$\alpha$}}}
\end{picture}

\caption{Coherence Conditions for Monoidal Categories \label{mon.cat}}
\end{figure}

\begin{defin} \index{monoidal functor, lax} \index{lax monoidal functor}
A {\em lax mon\-oid\-al functor} $F:{\cal C}\rightarrow {\cal D}$ between
two mon\-oid\-al categories $\cal C$ and $\cal D$ is a functor $F$ between
the underlying categories, equipped with a natural transformation
\[ \tilde{F}:  F(-)\otimes F(-)\rightarrow F(-\otimes -)\]

\noindent and a map
$F_I: I\rightarrow F(I)$, satisfying the hexagon and two squares of Figure 
\ref{lax.mon.f}.

If $\tilde{F}$ and $F_I$ are isomorphisms, $F$ is a {\em strong monoidal 
functor}.  If $F_I$ and all components of $\tilde{F}$ are identity maps,
$F$ is a {\em strict monoidal functor}. 

Lax, strong and strict semigroupal functors are defined similarly.
\index{semigroupal functor}

\end{defin}

\begin{figure}[htb] \centering

\setlength{\unitlength}{3947sp}%
\begingroup\makeatletter\ifx\SetFigFont\undefined
\def\x#1#2#3#4#5#6#7\relax{\def\x{#1#2#3#4#5#6}}%
\expandafter\x\fmtname xxxxxx\relax \def\y{splain}%
\ifx\x\y   
\gdef\SetFigFont#1#2#3{%
  \ifnum #1<17\tiny\else \ifnum #1<20\small\else
  \ifnum #1<24\normalsize\else \ifnum #1<29\large\else
  \ifnum #1<34\Large\else \ifnum #1<41\LARGE\else
     \huge\fi\fi\fi\fi\fi\fi
  \csname #3\endcsname}%
\else
\gdef\SetFigFont#1#2#3{\begingroup
  \count@#1\relax \ifnum 25<\count@\count@25\fi
  \def\x{\endgroup\@setsize\SetFigFont{#2pt}}%
  \expandafter\x
    \csname \romannumeral\the\count@ pt\expandafter\endcsname
    \csname @\romannumeral\the\count@ pt\endcsname
  \csname #3\endcsname}%
\fi
\fi\endgroup
\begin{picture}(4125,6930)(151,-6136)
\thinlines
\put(2101,-4936){\makebox(1.6667,11.6667){\SetFigFont{5}{6}{rm}.}}
\put(2101,-4936){\vector( 1, 0){1125}}
\put(3526,-5611){\vector( 0, 1){450}}
\put(1726,-5551){\vector( 0, 1){450}}
\put(3001,-5851){\vector(-1, 0){750}}
\put(1426,-5011){\makebox(0,0)[lb]{\smash{\SetFigFont{10}{14.4}{rm}$F(A\otimes I)$}}}
\put(3376,-5011){\makebox(0,0)[lb]{\smash{\SetFigFont{10}{14.4}{rm}$F(A)$}}}
\put(3301,-5911){\makebox(0,0)[lb]{\smash{\SetFigFont{10}{14.4}{rm}$F(A)\otimes I$}}}
\put(1351,-5911){\makebox(0,0)[lb]{\smash{\SetFigFont{10}{14.4}{rm}$F(A)\otimes F(I)$}}}
\put(2476,-4786){\makebox(0,0)[lb]{\smash{\SetFigFont{10}{14.4}{rm}$F(\rho)$}}}
\put(3676,-5536){\makebox(0,0)[lb]{\smash{\SetFigFont{10}{14.4}{rm}$\rho$}}}
\put(2251,-6136){\makebox(0,0)[lb]{\smash{\SetFigFont{10}{14.4}{rm}$F(A)\otimes F_0$}}}
\put(1051,-5461){\makebox(0,0)[lb]{\smash{\SetFigFont{10}{14.4}{rm}$F_{A,I}$}}}
\put(2101,-3061){\makebox(1.6667,11.6667){\SetFigFont{5}{6}{rm}.}}
\put(2101,-3061){\vector( 1, 0){1125}}
\put(3526,-3736){\vector( 0, 1){450}}
\put(1726,-3736){\vector( 0, 1){450}}
\put(3031,-3976){\vector(-1, 0){750}}
\put(1426,-3136){\makebox(0,0)[lb]{\smash{\SetFigFont{10}{14.4}{rm}$F(I\otimes A)$}}}
\put(3376,-3136){\makebox(0,0)[lb]{\smash{\SetFigFont{10}{14.4}{rm}$F(A)$}}}
\put(3301,-4036){\makebox(0,0)[lb]{\smash{\SetFigFont{10}{14.4}{rm}$I\otimes F(A)$}}}
\put(1351,-4036){\makebox(0,0)[lb]{\smash{\SetFigFont{10}{14.4}{rm}$F(I)\otimes F(A)$}}}
\put(2476,-2911){\makebox(0,0)[lb]{\smash{\SetFigFont{10}{14.4}{rm}$F(\lambda)$}}}
\put(3676,-3661){\makebox(0,0)[lb]{\smash{\SetFigFont{10}{14.4}{rm}$\lambda$}}}
\put(2251,-4261){\makebox(0,0)[lb]{\smash{\SetFigFont{10}{14.4}{rm}$F_0\otimes F(A)$}}}
\put(1051,-3586){\makebox(0,0)[lb]{\smash{\SetFigFont{10}{14.4}{rm}$F_{I,A}$}}}
\put(1726,464){\vector( 1, 0){1725}}
\put(2026,-2161){\vector( 1, 0){1275}}
\put(1126,-631){\vector( 0, 1){900}}
\put(4051,-571){\vector( 0, 1){825}}
\put(1126,-1966){\vector( 0, 1){750}}
\put(4051,-1951){\vector( 0, 1){825}}
\put(601,389){\makebox(0,0)[lb]{\smash{\SetFigFont{10}{14.4}{rm}$F([A\otimes B]\otimes C)$}}}
\put(601,-1036){\makebox(0,0)[lb]{\smash{\SetFigFont{10}{14.4}{rm}$F(A\otimes B)\otimes F(C)$}}}
\put(2326,614){\makebox(0,0)[lb]{\smash{\SetFigFont{10}{14.4}{rm}$F(\alpha)$}}}
\put(376,-286){\makebox(0,0)[lb]{\smash{\SetFigFont{10}{14.4}{rm}$F_{A\otimes B,C}$}}}
\put(4276,-361){\makebox(0,0)[lb]{\smash{\SetFigFont{10}{14.4}{rm}$F_{A,B\otimes C}$}}}
\put(151,-1636){\makebox(0,0)[lb]{\smash{\SetFigFont{10}{14.4}{rm}$F_{A,B}\otimes F(C)$}}}
 \put(3676,389){\makebox(0,0)[lb]{\smash{\SetFigFont{10}{14.4}{rm}$F(A\otimes [B\otimes C])$}}}
\put(3601,-961){\makebox(0,0)[lb]{\smash{\SetFigFont{10}{14.4}{rm}$F(A)\otimes F(B\otimes C)$}}}
\put(3526,-2236){\makebox(0,0)[lb]{\smash{\SetFigFont{10}{14.4}{rm}$F(A)\otimes [F(B)\otimes F(C)]$}}}
\put(376,-2236){\makebox(0,0)[lb]{\smash{\SetFigFont{10}{14.4}{rm}$[F(A)\otimes F(B)]\otimes F(C)$}}}
\put(2551,-2461){\makebox(0,0)[lb]{\smash{\SetFigFont{10}{14.4}{rm}$\alpha$}}}
\put(4276,-1561){\makebox(0,0)[lb]{\smash{\SetFigFont{10}{14.4}{rm}$F(A)\otimes F_{B,C}$}}}
\end{picture}

\caption{\label{lax.mon.f} Coherence Conditions for a Lax Monoidal Functor}
\end{figure}
 
We will have no cause to consider oplax monoidal functors in this work.

We refer to the components of the natural transformations and maps specified
in these definitions, and to their inverses (if any), as {\em structure maps}.\index{natural
transformation}\index{transformation, natural}\index{structure maps}
Likewise, a map which is obtained from some other map $f$ by forming an 
iterated mon\-oid\-al product of $f$ with identity maps for various objects is 
called a {\em prolongation} of $f$.\index{prolongation}
  Sometimes by abuse of terminology 
prolongations of structure maps are themselves referred to as structure maps.
  
We will also refer to a diagram obtained by applying the same iterated 
mon\-oid\-al product with identity maps to every map of a given diagram as
a prolongation\index{prolongation} of the given diagram.

It is, of course, a matter of taste whether one defines strong mon\-oid\-al functors
as lax mon\-oid\-al functors with invertible structure maps, as here, or as
oplax mon\-oid\-al functors with invertible structure maps as in
\cite{Y.def,Y.book}.

Crucial to the construction of our deformation theories are the coherence
theorem of Mac Lane \cite{Mac.coh} and a non-symmetric variant of the 
coherence theorem of Epstein \cite{Ep.coh}, which we
will soon state in the most convenient form for our purposes.

\begin{defin}
For any set $S$, $S\downarrow MonCat$ (resp. $S\downarrow SGCat$) 
is the category whose objects are (small)
mon\-oid\-al (resp. semigroupal) categories equipped with a map from 
$S$ to their set of objects, and 
whose arrows are strict mon\-oid\-al functors whose object maps commute with the
map from $S$.

$S\downarrow LaxSGFun$ (resp.  
$S\downarrow StrongSGFun$) is the category whose objects are lax (resp.  
strong) semigroupal functors 
between a pair of semigroupal categories, the source of which is equipped with
a map from $S$ to its set of arrows, and whose arrows are pairs of strict
mon\-oid\-al functors forming commuting squares and commuting with the map from
$S$.
\end{defin}

Observe that $S\downarrow MonCat$, (resp. $S\downarrow SGCat$,
$S\downarrow LaxSGFun$, 
$S\downarrow StrongSGFun$) is a
category of models of an essentially algebraic theory, and thus by general
principles has an initial object. We refer to this initial object as the
{\em free mon\-oid\-al category} (resp. {\em semigroupal category, lax semigroupal
functor, strong semigroupal functor}) {\em on} $S$.

\begin{defin} \index{diagram, formal} \index{formal diagram}
A {\em formal diagram in the theory of mon\-oid\-al categories (resp. semigroupal
categories)} is a diagram in the free mon\-oid\-al (resp. semigroupal) category on
$S$ for some set $S$.

A {\em formal diagram in the theory of lax (resp. strong) semigroupal 
functors} is a diagram in the target category 
of the free lax (resp. strong) semigroupal functor on $S$ 
for some set $S$.
\end{defin}

The coherence theorem of Mac Lane \cite{Mac.coh} may then be stated as

\begin{thm} \label{MacLane.coh} \index{coherence theorem, Mac Lane's} 
\index{Mac Lane's coherence theorem}
Every formal diagram in the theory of mon\-oid\-al categories commutes.  
Consequently, any diagram which is the image of a formal diagram under a
(strict mon\-oid\-al) functor commutes.
\end{thm}

The same proof carries the weaker result:

\begin{thm}
Every formal diagram in the theory of semigroupal categories commutes.
Consequently, any diagram which is the image of a formal diagram under
a (strict semigroupal) functor commutes.\end{thm}

Epstein \cite{Ep.coh} proved a coherence theorem only for lax semigroupal 
functors
between symmetric semigroupal categories, but the same proof will carry the
result:

\begin{thm} Every formal diagram in the theory of lax (resp. strong)
semigroupal functors commutes.  Every formal diagram in the theory of
strong monoidal functors commutes.
Consequently, any diagram which is a functorial image of such a formal
diagram under a (strict mon\-oid\-al) functor commutes.
\end{thm}

These coherence theorems are the basis for a very useful notion and notational
convention:  throughout our discussion of categorical deformation theory we 
will use {\em padded composition operators} $\lceil \;\; \rceil$.
These operators are an embodiment of the coherence theorems of Mac Lane
\index{Mac Lane's coherence theorem} \index{coherene theorem, Mac Lane's}
\cite{Mac.coh}
and Epstein \cite{Ep.coh}\index{Epstein's coherence theorem} \index{coherence
theorem, Epstein's}.  

\begin{defin} \index{padded composition} \index{$\lceil \rceil$}
Given a mon\-oid\-al category $\cal C$ (resp. a semigroupal or strong
monoidal functor 
$F:{\cal X}\rightarrow {\cal C}$), and a sequence of maps $f_1,\ldots ,f_n$
in $\cal C$ such that the source of $f_{i+1}$ is isomorphic (resp. maps) 
to the target of $f_i$
by a composition of prolongations of structure maps (i.e. by a 
formal diagram with
underlying diagram a chain of composable maps), we let  

\[ \lceil f_1,\ldots ,f_n \rceil \]

\noindent denote the composite $a_0 f_1 a_1 f_2 \ldots a_{n-1} f_n a_n$, where
the $a_i$'s are composites of prolongations of structure maps and the
following hold:
\begin{enumerate}
\item The source of
$a_0$ is reduced (no tensorands of $I$) and completely left-parenthesized  
(resp. 
reduced and completely left-parenthesized and free from images of mon\-oid\-al
products under $F$).
\item The target of $a_n$ is reduced and completely right-parenthesized
(resp. reduced and completely right-parenthesized and free from products 
both of whose factors are images under $F$).
\item The composite is well-defined.
\end{enumerate}

\end{defin}

The fact that this defines a well-defined map will be a consequence of 
the coherence 
theorems.
However, if the $f_i$ are simply maps, $\lceil \;\; \rceil$ may not 
be well-defined in the event that there
are ``accidental coincidences''.  
In our circumstance, the maps in the sequences to which the padded composition
operator is applied will always be components of natural transformations with a
particular structure:

\begin{defin} \index{paracoherent} \index{natural transformation, paracoherent}
Given a mon\-oid\-al category $\cal C$ (resp. a mon\-oid\-al functor 
$F:{\cal C}\rightarrow
{\cal D}$), a natural transformation is {\em $\cal C$-paracoherent} 
(resp. {\em $F$-paracoherent}) if its source and target functors are iterated 
prolongations of the structure functors $\otimes$, $I$, and $1_{\cal C}$ 
(resp. $\otimes$, $F$,
$I$, $1_{\cal C}$, and $ 1_{\cal D}$), where $I$ is regarded as a 
functor from the trivial
one object category.
\end{defin}

In the case where the maps in the sequence are specified not merely as 
maps, but
as components of particular paracoherent natural transformations, their 
sources and targets are given an explicit structure as images of 
iterated prolongations of 
structure functors. Thus we may require that the ``padding'' maps given
in terms of the structural natural transformations be (components of) natural
transformations between the appropriate functors.

A number of elementary lemmas hold for these operators.  Proofs are left
to the reader.

\begin{lemma} \label{first.ceil.lemma}

\[ \lceil f_1\ldots f_n \rceil = 
\lceil \lceil f_1\ldots f_k\rceil \lceil f_{k+1}\ldots
  f_n\rceil \rceil .\]

\end{lemma}

\begin{lemma}

If $\lceil \; \rceil$ is applied in the case of a monoidal category or
strong monoidal functor

\[ \lceil f_1\ldots g\otimes I \ldots f_n\rceil =  
\lceil f_1\ldots g \ldots f_n\rceil =
\lceil f_1\ldots I\otimes g \ldots f_n\rceil .\]

\end{lemma}

\begin{lemma}

\[ \lceil f_1\ldots f_n \rceil = 
\lceil f_1\ldots \lceil f_k\ldots f_l \rceil \ldots f_n \rceil
 .\]

\end{lemma}

\begin{lemma}

\[ \lceil f_1\ldots g\otimes h \ldots f_n\rceil =  
  \lceil f_1\ldots \lceil g \rceil \otimes h \ldots f_n\rceil =
		\lceil f_1\ldots g\otimes \lceil h \rceil \ldots f_n\rceil .\]

\end{lemma}

\begin{lemma} \label{last.ceil.lemma}
If $\phi_{X_1,\ldots ,X_n}$ is a $\cal C$-paracoherent natural transformation
(resp. $F$-paracoherent natural transformation, for $F$ a strong mon\-oid\-al functor), 
then so
is $\phi_{X_1,\ldots, I, \ldots, X_n}$, where I is inserted in the $i^{th}$ position,
and similarly if $I$ is inserted in the $i^{th}$ position for all $i\in T\subset \{1,\ldots n\}$. Moreover, 
in this latter case $\lceil \phi_{\ldots} \rceil$ is a paracoherent
natural transformation from the fully left-parenthesized product 
of $X_{i_1}\ldots X_{i_k}$  (resp. the
fully left-parenthesized product of  $F(X_{i_1}), \ldots , F(X_{i_k})$) 
to the fully 
right-parenthesized product of (resp. $F$ of the fully 
right-parenthesized product) of $X_{i_1}\ldots X_{i_k}$, 
where 

\[ \{i_1,\ldots ,i_k\} =
\{1,\ldots ,n \} \setminus T \]

\noindent and $i_1 < i_2 < \ldots i_n$.
\end{lemma}

From these lemmas we deduce a final lemma:

\begin{lemma} \label{unit.ceil.commuting.lemma}
If $\psi_{A,B,C}, \phi_{A,B,C}:[A\otimes B]\otimes C\rightarrow A\otimes[B\otimes C]$
are natural transformations, then

\[ \lceil [\phi_{A,I,I}\otimes B] \psi_{A,I,B} \rceil = 
\lceil \psi_{A,I,B} [\phi_{A,I,I}\otimes B] \rceil \]

\noindent and

\[ \lceil [A\otimes \phi_{I,I,B}] \psi_{A,I,B} \rceil = 
\lceil \psi_{A,I,B} [A\otimes \phi_{I,I,B}] \rceil .\]

\end{lemma}

\noindent {\bf proof:} First, apply
Lemma \ref{first.ceil.lemma}. Then use the naturality of $\lceil \psi_{A,I,B} \rceil:A\otimes B \rightarrow A\otimes B$ and the source and target data for $\phi_{A,I,I}$ and $\phi_{I,I,B}$, as given by Lemma \ref{last.ceil.lemma}. $\Box$
\medskip

Finally, we make

\begin{defin} \index{monoidal natural transformation} 
\index{natural transformation, monoidal}
A {\em mon\-oid\-al natural transformation} is a natural transformation
$\phi:F\Rightarrow G$ between mon\-oid\-al functors which satisfies

\[ \tilde{G}_{A,B}(\phi_{A\otimes B}) = \phi_A \otimes \phi_B (\tilde{F}_{A,B})
\]

\noindent and $F_0 = G_0(\phi_I)$.  A {\em semigroupal natural transformation}
between semigroupal functors is defined similarly.
\end{defin}

\noindent and

\begin{defin} \index{monoidal equivalence} \index{equivalence, monoidal}
 A {\em mon\-oid\-al equivalence} between mon\-oid\-al categories $\cal C$
and $\cal D$ is an equivalence of categories in which the functors 
$F:{\cal C}\rightarrow {\cal D}$ and $G:{\cal D}\rightarrow {\cal C}$ are equipped with
the structure of mon\-oid\-al functors, and the natural isomorphisms
$\phi:FG\Rightarrow Id_{\cal C}$ and $\psi:GF\Rightarrow Id_{\cal D}$ are both
mon\-oid\-al natural transformations.  If there exists a mon\-oid\-al equivalence between
$\cal C$ and $\cal D$, we say that $\cal C$ and $\cal D$ are {\em mon\-oid\-ally
equivalent}. \end{defin}

Using this, we can now define equivalences between deformations of monoidal and
semigroupal
categories and functors.

\begin{defin}
Two $n^{th}$ order deformations of a monoidal (semigroupal) category
$\cal C$ are equivalent if there exists a monoidal (semigroupal) functor
between them, whose underlying functor is $Id_{{\cal C}^{(n)}}$ and whose
structural natural transformations reduces modulo
$\epsilon$ to the identity natural transformations.
\end{defin}

\begin{defin}
Two $n^{th}$ order purely functorial deformations $\tilde{F}$ and $\hat{F}$
of a monoidal (semigroupal) functor
$F:{\cal C}\rightarrow {\cal D}$ are equivalent if there exists a
monoidal (semigroupal) 
natural isomorphism $\vartheta:\tilde{F}\Rightarrow \hat{F}$
which reduces modulo $\epsilon$ to the identity natural transformation.
\end{defin}

Observe that in this last case, it is a matter of indifference whether we
are considering strong or lax monoidal functors.

For purely transformational monoidal (semigroupal) natural transformations
equivalence of deformations is simply equality.

\clearpage 
  \section{Deformation Complexes of 
Semigroupal Categories and Functors}
\vspace*{1cm}

We define a cochain complex associated to any semigroupal
category or semigroupal functor:
 
\begin{defin} The {\em deformation complex} 
\index{deformation complex of a lax semigroupal functor} 
\index{complex, deformation}
of a strong or lax semigroupal functor
$(F:{\cal C}\rightarrow {\cal D}, \tilde{F})$ is the cochain complex 

\[(X^\bullet(F), \delta),\]

\noindent where

\[ X^n(F) = {\rm Nat}(^n\otimes(F^n), F(\otimes^n)) \]   

\noindent and

\begin{eqnarray*}
\delta(\phi)_{A_0,\ldots ,A_n} & = & 
	\lceil F(A_0)\otimes \phi_{A_1,\ldots ,A_n} \rceil \\
 & &\;\;\; + \sum_{i=1}^n (-1)^i \lceil
\phi_{A_0,\ldots,A_{i-1}\otimes A_i,\ldots ,A_n} \rceil \\
 & &\;\;\; + (-1)^{n+1}
\lceil \phi_{A_0,\ldots ,A_{n-1}}\otimes F(A_n) \rceil .
\end{eqnarray*}

\end{defin}

A similar definition applies for oplax semigroupal functors.
Observe that here we regard strong semigroupal and monoidal functors
as special cases of lax semigroupal functors, while in \cite{Y.book} they
are treated as special cases of oplax functors.

\begin{defin}
The {\em deformation complex} of a semigroupal category $({\cal C},\otimes,
\alpha)$ is $(X^\bullet(Id_{\cal C}), \delta)$.  We also denote this by
$(X^\bullet({\cal C}), \delta)$.
\end{defin}

The motivation for these definitions is given in \cite{CY.def} and
\cite{Y.def}, or can be readily discovered by the reader by computing
by hand the conditions on the term $\alpha^{(1)}$ or $\tilde{F}^{(1)}$  
in a first order
deformation.

We also have

\begin{thm} \label{cat.first}
The first-order deformations of a semigroupal category 
$\cal C$ are classified 
up to equivalence by $H^3({\cal C})$. \index{deformations, first-order}
\index{deformations, classification theorems for}
\end{thm}

\noindent{\bf sketch of proof:} Consider two first-order deformations
$\tilde{\alpha} = \alpha + \alpha^{(1)}\epsilon$ and
$\hat{\alpha} = \alpha + a^{(1)}\epsilon$
of $\cal C$. Consider also a semigroupal functor whose underlying 
functor is the
identity functor and whose structural transformation is of the form

\[ 1_{A\otimes B} + \phi_{A,B}\epsilon: A\otimes B \rightarrow A \otimes B. \]

Now, write out the coherence condition for semigroupal functors in this
case, and look at the degree 1 terms.  The resulting
equation is nothing more than

\[  \alpha^{(1)} - a^{(1)} = \delta(\phi) .\]

\noindent$\Box$\medskip

In \cite{Y.def} it is shown that

\begin{thm} \label{fun.first} The purely functorial
first-order deformations of a semigroupal functor
$F:{\cal C} \rightarrow {\cal D}$ are classified up to equivalence by
$H^2(F)$. \index{deformations, first-order}
\index{deformations, classification theorems for}
\end{thm}

The proof, similar to that of the previous theorem, is given in \cite{Y.def}.

Obstructions to extending $n^{th}$ order deformations to $(n+1)^{st}$ order
deformations will be discussed below.

First, however, we will show that the deformation complexes 
already defined together with standard machinery from homological algebra
suffice to handle the cases of fibred and total deformations and
deformations of braided mon\-oid\-al categories.

Compositions of natural transformations with functors induce
two cochain maps whenever we have a semigroupal functor $F:{\cal C} \rightarrow
{\cal D}$ (whether lax, oplax, or strong):\index{cochain map, associated
to a mon\-oid\-al functor}

\[ \lceil F(-)\rceil:X^\bullet ({\cal C})\rightarrow X^\bullet (F) \]

\noindent and

\[ \lceil (-)_{F^\bullet} \rceil: X^\bullet ({\cal D})\rightarrow 
X^\bullet (F).
 \]

Recalling that $X^\bullet({\cal C}) = X^\bullet(Id_{\cal C})$, we see that the 
cochain maps just defined are, in fact, special cases of more general
families defined for composable pairs of functors

\[ {\cal C}\stackrel{F}{\longrightarrow}{\cal D}\stackrel{G}{\longrightarrow}
{\cal E} \]

\noindent namely,

\[ \lceil G(-)\rceil :X^\bullet(F)\rightarrow X^\bullet(G(F)) \]

\noindent and

\[ \lceil (-)_{F^\bullet}\rceil :X^\bullet(G)\rightarrow X^\bullet(G(F)) . \]

To consider deformations of braided mon\-oid\-al categories 
the following is also useful:  
given any $K$-linear semigroupal category $\cal C$, 
there is a ``diagonal'' cochain map \index{cochain map, diagonal}

\[ \Delta:X^\bullet ({\cal C})\rightarrow X^\bullet({\cal C}\boxtimes {\cal C})
  \]

\noindent given by:  

\[ \Delta(\phi) = \phi \boxtimes \lceil Id \rceil + \lceil Id \rceil \boxtimes
\phi \;\;.\]

These cochain maps allow us to assemble the simpler deformation complexes for
semigroupal categories and semigroupal functors into complexes whose
cohomology is related to more general types of deformations.

Recall the construction of a cone over a cochain map:

\begin{defin} \index{cone} 
Given a map of cochain complexes $u^\bullet :A^\bullet \rightarrow B^\bullet$,
the cone on $u^\bullet$ is the cochain complex

\[  (C_u^\bullet, d_u) = \left( B^\bullet \oplus A^{\bullet+1}, 
			\left[ \begin{array}{cc} 
				d_B & 0 \\
				u & -d_A  \end{array} \right] \right). \]

\end{defin}


Here we adopt the convention that elements of direct sums are written as
row vectors with entries in the summands, and that
arrays of maps act on the right
by matrix multiplication (with the action of maps in lieu of scalar 
multiplication). Note that this is consistent with our notational convention: 
maps {\em act} on the right on elements (improperly) 
thought of as maps, unless parentheses denoting application intervene.

In the next section use cones on the cochain maps discussed above to classify
first order deformations more general than those classified by Theorems
\ref{cat.first} and \ref{fun.first}

\clearpage 
  \section{First Order Deformations}
\vspace*{1cm}

Let us now consider the problem of classifying first order fibred 
deformations of semigroupal 
functors.\index{deformation, first-order}\index{deformation, fibred} 
If we have a lax (resp. oplax, strong)
semigroupal functor $[F,\tilde{F}]:({\cal C},\otimes ,\alpha )\rightarrow
({\cal D},\otimes , a)$, and we replace $\tilde{F}^{(0)} = \tilde{F}$ with
$\tilde{F}^{(0)} + \tilde{F}^{(1)}\epsilon$ and $\alpha^{(0)} = \alpha$
with $\alpha^{(0)} + \alpha^{(1)}\epsilon$ for $\epsilon^2 = 0$, the conditions
for the new coherence diagrams to commute become

\[ \delta(\alpha^{(1)}) = 0 \]

\noindent and

\[ \delta(\tilde{F}^{(1)}) + \lceil F(\alpha) \rceil = 0, \]

\noindent as can be verified readily by computing the $\epsilon$-degree
1 terms going around the pentagon and hexagon coherence diagrams.

It then follows directly that the pair $[\tilde{F}^{(1)},\alpha^{(1)}]$ is
a 2-cocycle in

\[ (C_{\lceil F(-) \rceil}^\bullet, d_{-\lceil F(-) \rceil}). \]

Now, consider the condition that two such 2-cocycles 
$[\tilde{F}^{(1)}_1,\alpha^{(1)}_1]$ and $[\tilde{F}^{(1)}_2,\alpha^{(1)}_2]$
are equivalent. Let $F_1:{\cal C}_1\rightarrow {\cal D}$ and 
$F_2:{\cal C}_2\rightarrow {\cal D}$ denote the semigroupal functors from the
corresponding deformations (suppressing here the naming of structural
maps). In particular, there is a structure map which makes 
$Id_{{\cal C}^{(2)}}$ into a (necessarily strong) semigroupal functor
and which reduces modulo $\epsilon$ to the identity natural transformation. 
Second, there is a semigroupal natural isomorphism $\psi$ from $F_1$ to
$F_2({\frak I})$ which reduces modulo $\epsilon$ to the identity,
where $\frak I$ is the identity functor on ${\cal C}^{(2)}$ made into a 
semigroupal functor by given structure
map. 

Denoting the structural map for  $Id_{{\cal C}^{(2)}}$ by $id + \iota^{(1)}\epsilon$ and 
letting $\psi = id + \psi^{(1)}\epsilon$, the coherence conditions become

\begin{eqnarray*}
\lefteqn{ [\tilde{F} + \tilde{F}^{(1)}_{2 A,B}\epsilon ]
([id + \iota^{(1)}_{A,B}\epsilon]
(id_{F(A\otimes B)} + \psi^{(1)}_{A\otimes B}\epsilon)) = } \\
&  & 
[[id_{F(A)} + \psi^{(1)}_{A}\epsilon]\otimes [id_{F(B)} 
+ \psi^{(1)}_{B}\epsilon]]
(\tilde{F} + \tilde{F}^{(1)}_{2 A,B}\epsilon) 
\end{eqnarray*}

\noindent and
 
\begin{eqnarray*}
\lefteqn{[id_{F(A)}\otimes 
	[id_{F(B\otimes C)} + \iota^{(1)}_{B,C}\epsilon]]} \\
\lefteqn{([id_{F(A\otimes[B\otimes C])} + 
\iota^{(1)}_{A,B\otimes C}\epsilon](F(\alpha + 
\alpha^{(1)}_{1 A,B,C}\epsilon))) =  }\\
& & [F(\alpha + 
\alpha^{(1)}_{2 A,B,C}\epsilon)]([[id_{F(A\otimes B)} + \iota^{(1)}_{A,B}]\otimes
id_{F(C)}] \\ 
& & \;\;([id_{F([A\otimes B]\otimes C)} + 
\iota^{(1)}_{A\otimes B,C}\epsilon])) \;\;. 
\end{eqnarray*}

	Using the bilinearity of composition and $\otimes$, the coherence
conditions on the original maps, and the condition $\epsilon^2 = 0$, these
readily reduce to

\[ \tilde{F}^{(1)}_1  - \tilde{F}^{(1)}_2 = \iota^{(1)} - \delta(\psi^{(1)}) 
\]

\noindent and

\[ \alpha^{(1)}_1 - \alpha^{(1)}_2 = \delta(\iota^{(1)}) \;\;. \]

We have thus demonstrated

\begin{thm} The 
first order fibred deformations of a semigroupal functor
$F:{\cal C} \rightarrow {\cal D}$ are classified up to equivalence by
the third cohomology of the cone $C^\bullet_{\lceil F(-) \rceil} = 
X_{\rm fibred}^\bullet(F)$. \index{deformation, fibred}  \index{deformations,
classification theorems for}
\end{thm}

A similar analysis shows

\begin{thm} The 
first order total deformations of a sem\-i\-group\-al functor
$F:{\cal C} \rightarrow {\cal D}$ are classified up to equivalence by
the third cohomology of the cone $C^\bullet_{\lceil F(p_1) \rceil  - \lceil (p_2)_{F^\bullet}\rceil} = X_{\rm total}^\bullet(F)$.\index{deformation, total}
  \index{deformations, classification theorems for}
\end{thm}

The case of total deformations of a multiplication (or equivalently, 
deformations
of a braided mon\-oid\-al category) presents another subtlety:  the
source and target must be deformed in tandem.  

\begin{propo}
If ${\cal C}^{(n)}, \otimes, \alpha^{(0)} + \alpha^{(1)}\epsilon + \ldots 
+ \alpha^{(n)}\epsilon^n$ 
is an $n^{th}$-order deformation of $({\cal C}, \otimes, \alpha)$ and

\[ \beta^{(k)} = \sum_{i=0}^k \alpha^{(i)}\boxtimes \alpha^{(k-i)} \;\;, \]

\noindent then
$([{\cal C}\boxtimes{\cal C}]^{(n)}, \otimes \boxtimes \otimes, 
\beta^{(0)} + \beta^{(1)}\epsilon + \ldots \beta^{(n)}\epsilon^n) $
is an $n^{th}$-order deformation of $({\cal C}\boxtimes{\cal C}, 
\otimes \boxtimes \otimes, \alpha \boxtimes \alpha)$.
We call this deformation the {\em diagonal deformation} of 
${\cal C}\boxtimes{\cal C}$. \index{deformation, diagonal} \index{diagonal
deformation}
\end{propo}

\noindent{\bf proof:}
Observe first that ${\cal C}\boxtimes {\cal C}$ is defined with respect
to the commutative ring $R$, and that $[{\cal C}\boxtimes_R {\cal C}]^{(n)}$
is canonically isomorphic to 
${\cal C}^{(n)}\boxtimes_{R[\epsilon]/\!<\!\epsilon^n\!>\!}{\cal C}^{(n)}$.

The diagonal deformation is then simply the 
$R[\epsilon]/\!<\!\epsilon^n\!>\!$-linearized version of the 
diagonal semigroupal
structure induced on ${\cal C}^{(n)}\times {\cal C}^{(n)}$ by the (deformed)
semigroupal structure on ${\cal C}^{(n)}$.  The formula for the $\beta^{(k)}$'s
is derived by simply collecting terms according to their degree in $\epsilon$.
$\Box$
\medskip

\begin{defin} 
A {\em coarse deformation} of a multiplication is a total deformation of 
the semigroupal functor such that
the deformation of the source ${\cal C}\boxtimes {\cal C}$ is the diagonal
deformation induced by the deformation of the target. \index{deformation,
coarse} \index{coarse deformation of a multiplication} 
A {\em deformation} of a multiplication (and thus of a braided mon\-oid\-al
category) is a coarse deformation which is equipped with natural isomorphisms
as required to make it into a multiplication.  \index{deformation of a 
multiplication}
\end{defin}

We will consider the behavior
of units in general in Section \ref{unit.chapter}, so we here confine ourselves
to consider the appropriate deformation complex for coarse deformations of
multiplications:

Consider the composite cochain map 

\[ \phi:X^\bullet({\cal C})\stackrel{(\Delta,Id)}{\longrightarrow} 
	X^\bullet({\cal C}\boxtimes {\cal C})\oplus X^\bullet({\cal C})
	\stackrel{\lceil \Phi(p_1) \rceil - \lceil (p_2)_{\Phi^\bullet} \rceil}
	{\longrightarrow} X^\bullet(\Phi) .\]

\noindent An argument similar to that given above for fibred deformations
shows that: 

\begin{thm} The first order coarse deformations of a multiplication 
$\Phi:{\cal C}\boxtimes {\cal C}\rightarrow {\cal C}$ are classified up to
equivalence by the third cohomology of the cone 
$C_\phi^\bullet = X_{\rm coarse}^\bullet(\Phi)$. \index{deformations,
classification theorems for}
\end{thm}

Since $\phi$ is defined as a composite, something more remains to be said:
if we consider our cochain complexes as objects in the homotopy category
$K^+(R)$ or the derived category $D^+(R)$, the octahedral property ensures
the existence of an exact triangle relating 
$X_{\rm coarse}^\bullet(\Phi)$, 
$X_{\rm total}^\bullet(\Phi)$ and $C_{(\Delta,Id)}$, and thus of a 
long-exact sequence in cohomology.

\clearpage 
  \section{Obstructions and the Cup 
Product and  Pre-Lie Structures on 
$X^\bullet(F)$}
\vspace*{1cm}

The cochain complex associated to any of the types of
semigroupal functors shares
many of the properties of the Hochschild 
complex\index{Hochschild complex}\index{complex, Hochschild}
of an associative 
algebra $A$ with coefficients in $A$, which were described by Gerstenhaber
\cite{G1,G2,GS}.
In particular, we have two products defined on cochains.  The first, 
the cup product, \index{cup product} \index{product, cup}

\[ -\cup - : X^n(F) \times X^m(F)\rightarrow X^{n+m}(F) \;, \]

\noindent is given by

\[ G\cup H_{A_1, \ldots A_{n+m}} = \lceil 
G_{A_1, \ldots ,A_n}\otimes H_{A_{n+1}, \ldots ,A_{n+m}} \rceil \;. \]

\noindent The second, the composition product, \index{product, composition}
\index{composition product}

\[\langle -,- \rangle
:X^n(F)\times X^m(F)\rightarrow X^{n+m-1}(F) \]

\noindent is given by

\begin{eqnarray*} \lefteqn{\langle G,H \rangle_{A_1, \ldots A_{n+m-1}} 
= } \\
 & & \sum (-1)^{mi} \lceil (G_{A_1,\ldots ,A_i,A_{i+1}\otimes \ldots \otimes A_{i+n}, A_{i+n+1},\ldots
A_{n+m-1}}) \\
 & & \hspace*{.5cm}
F(A_1)\otimes \ldots F(A_i) \otimes H_{A_{i+1},
\ldots ,A_{i+n}}\otimes F(A_{i+n+1})\otimes \\ 
 & & \hspace*{.5cm} \ldots \otimes F(A_{n+m-1}) \rceil 
\end{eqnarray*}

\noindent in the case of strong and lax semigroupal functors.  A similar
formula is given in \cite{Y.book} for the oplax case.

\begin{propo}
The product $\langle -,- \rangle$
comes from a ``pre-Lie system'',  in the terminology
of Gerstenhaber \cite{G1}, given by \index{pre-Lie system}

\begin{eqnarray*}
\lefteqn{\langle G,H \rangle^{(i)}_{A_1, \ldots A_{n+m-1}} 
= } \\ 
 & & \lceil (G_{A_1,\ldots ,A_i,A_{i+1}\otimes \ldots \otimes A_{i+n}, 
A_{i+n+1},\ldots
A_{n+m-1}}) F(A_1)\otimes \ldots \\
 & & \hspace*{.5cm}
 F(A_i) \otimes H_{A_{i+1},
\ldots ,A_{i+n}}\otimes F(A_{i+n+1})\otimes \ldots \otimes F(A_{n+m-1}) \rceil
\end{eqnarray*}

\noindent in the case of strong and lax semigroupal functors, where 
$X^n(F)$ has degree $n-1$.  
\end{propo}

Again in \cite{Y.book} the oplax case is treated as well.

\noindent{\bf proof:} 
First, note that the ambiguities of parenthesization in the semigroupal 
products
in this definition are rendered irrelevant by the $\lceil \;\; \rceil$ on
each term, by virtue of the coherence theorems
for semigroupal functors.

It is obvious that the product is given by a sum of
these terms with the correct signs for the construction of a Lie bracket from
a pre-Lie system, so actually the content of the proposition is that the
$\langle - , - \rangle^{(i)}$'s satisfy the definition of a pre-Lie system. 
\index{pre-Lie system}
That is, for $G \in X^m(F)$, $H \in X^n(F)$ and $K \in X^p(F)$, we have

\[ \langle \langle G,H \rangle^{(i)}, K \rangle^{(j)} = \left\{
	\begin{array}{ll}
	 \langle \langle G,K \rangle^{(j)}, H \rangle^{(i+p-1)} & \mbox{if 
		$0 \leq j \leq i-1$} \\
	  \langle G, \langle H, K \rangle^{(j-i)} \rangle^{(i)} & \mbox{if
		$i \leq j \leq n$}
	\end{array} \right. \]

\noindent (recall that a $k$-chain has degree $k-1$).

This is a simple computational check.  In verifying
the first case,  naturality
will allow the prolongations of $K$ and $H$ to commute. $\Box$\medskip

Now, suppose we have an ${M-1}^{st}$ order deformation

\[ \tilde{\alpha} = \alpha^{(0)} + \alpha^{(1)}\epsilon + \ldots + 
          \alpha^{(M-1)}\epsilon^{M-1} \;. \]

As was shown in \cite{CY.def}, the obstruction to extending this to an $M^{th}$
order deformation is the 4-cochain \index{obstruction}

\begin{eqnarray*}\omega_{A,B,C,D}^{(M)} & = & 
	\sum_{\mbox{\scriptsize $\begin{array}{c}
				i + j  =  M \\
				0\leq i,j  <  M \end{array}$}}
\lceil \alpha_{A\otimes B,C,D}^{(i)} \alpha_{A,B,C\otimes D}^{(j)} \rceil \\
  &  &\;  - \hspace{-.25cm} \sum_{\mbox{\scriptsize $\begin{array}{c}
                   i + j + k  =  M \\
		            0\leq i,j,k  <  M \end{array}$}}\hspace*{-.25cm}
   \lceil [\alpha_{A,B,C}\otimes D] \alpha_{A,B\otimes C, D}^{(j)} [A\otimes 
	\alpha_{B,C,D}^{(k)}] \rceil  .\end{eqnarray*}

The deformation extends precisely when this cochain is a coboundary, in which
case $\alpha^{(M)}$ may be any solution to $\delta (\alpha^{(M)}) = \omega^{(M)}$.

In \cite{Y.book} it is shown that

\begin{thm} \index{obstruction}
For all $M$, the obstruction $\omega^{(M)}$ is a 4-cocycle.  Thus, an 
$(M-1)^{st}$ order deformation extends to an $M^{th}$ order deformation if
and only if the cohomology class $[\omega^{(M)}] \in H^4(\cal C)$ vanishes.
\end{thm}

\noindent {\bf proof:}  The proof is essentially computational.  Following
\cite{Y.book} we introduce notation for the summands of the coboundary
of a 3-cochain $\phi_{A,B,C}$ ---

\begin{eqnarray*}
\partial_0 \phi_{A,B,C,D} & = & A\otimes \phi_{B,C,D} \\
\partial_1 \phi_{A,B,C,D} & = & \phi_{A\otimes B,C,D} \\
\partial_2 \phi_{A,B,C,D} & = & \phi_{A,B\otimes C,D} \\
\partial_3 \phi_{A,B,C,D} & = & \phi_{A,B,C\otimes D} \\
\partial_4 \phi_{A,B,C,D} & = & \phi_{A,B,C}\otimes D 
\end{eqnarray*}

--- and of a 4-cochain $\psi_{A,B,C,D}$ ---

\begin{eqnarray*}
\underline{\partial}_0 \psi_{A,B,C,D,E} & = & A\otimes \psi_{B,C,D,E} \\
\underline{\partial}_1 \psi_{A,B,C,D,E} & = & \psi_{A\otimes B,C,D,E} \\
\underline{\partial}_2 \psi_{A,B,C,D,E} & = & \psi_{A,B\otimes C,D,E} \\
\underline{\partial}_3 \psi_{A,B,C,D,E} & = & \psi_{A,B,C\otimes,D,E} \\
\underline{\partial}_4 \psi_{A,B,C,D,E} & = & \psi_{A,B,C,D\otimes E} \\
\underline{\partial}_5 \psi_{A,B,C,D,E} & = & \psi_{A,B,C,D}\otimes E 
\end{eqnarray*}

(We include the underline stroke only for ease of reading, not out of any
logical necessity.)

We then have

\[ \delta (\phi )_{A,B,C,D} = \sum_{i=0}^4 (-1)^{i+1} \partial_i \phi_{A,B,C,D} \]

\noindent for 3-cochains $\phi$, and

\[ \delta (\psi )_{A,B,C,D,E} = \sum_{i=0}^5 (-1)^{i+1} 
           \underline{\partial}_i \psi_{A,B,C,D,E} \]

\noindent for 4-cochains $\psi$.

In this notation the obstruction cochain $\omega^{M}$ becomes

\[ \omega^{(M)} = \hspace*{-1mm}\sum_{\mbox{\scriptsize $\begin{array}{c}
					i + j  =  M \\
					0\leq i,j  <  M \end{array}$}} 
\hspace*{-1.5mm} \lceil \partial_1 \alpha^{(i)} \partial_3 \alpha^{(j)} \rceil
-  \hspace*{-1mm}\sum_{\mbox{\scriptsize $\begin{array}{c}
                   i + j + k  =  M \\
		            0\leq i,j,k  <  M \end{array}$}}
   \hspace*{-1.5mm} \lceil \partial_4 \alpha^{(i)} \partial_2 \alpha^{(j)} 
	\partial_0 \alpha^{(k)} \rceil \; ,
 \]

\noindent while the vanishing of the obstruction $\omega^{(N)}$ (for
$N < M$) becomes

\begin{eqnarray*} 0 & = & \delta \alpha^{(N)} +\omega^{(N)} \\
 & = &  \hspace*{-1mm}\sum_{\mbox{\scriptsize $\begin{array}{c}
		i + j  =  N \\
		0\leq i,j  \leq  N \end{array}$}}
 \hspace*{-1.5mm} \lceil \partial_1 \alpha^{(i)} \partial_3 \alpha^{(j)} \rceil
-  \hspace*{-1mm}\sum_{\mbox{\scriptsize $\begin{array}{c}
                   i + j + k  =  N \\
	            0\leq i,j,k  \leq  N \end{array}$}}
  \hspace*{-1.5mm} \lceil \partial_4 \alpha^{(i)} \partial_2 \alpha^{(j)} 
	\partial_0 \alpha^{(k)} \rceil \; .
 \end{eqnarray*}

\noindent We wish to show that

\[  \sum_{i=0}^5 (-1)^{i+1} 
           \underline{\partial}_i \omega_{A,B,C,D,E}^{(M)} = 0 \; .\]

Observe that $\omega^{(1)} = 0$ and $\delta(\alpha^{(1)}) = 0$, so we
may proceed by induction under the assumption that $\omega^{(N)}$ and
$\alpha^{(N)}$
satisfy 

\[ \delta(\alpha^{(N)}) + \omega^{(N)} = 0 \] 

\noindent for  $N < M$.

\begin{figure}[htb] 
\epsfig{file=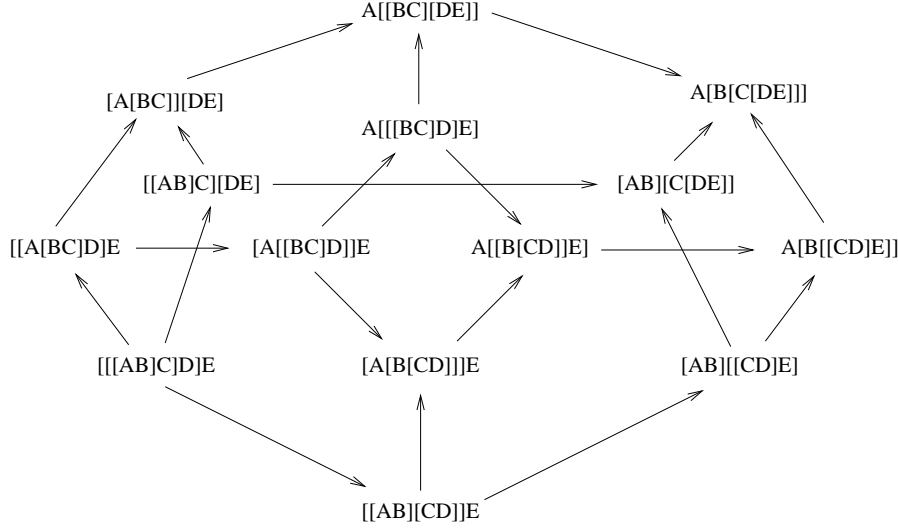,width=4.7in}
\caption{The Associahedron \label{associahedron}}
\end{figure}

It is convenient to picture the summands of the left-hand side in terms of compositions of maps along the boundaries of faces of the ``associahedron'' (or
3-dimensional Stasheff polytope) \cite{Stasheff} given in Figure 
\ref{associahedron}. \index{associahedron} \index{Stasheff polytope}

Suppose we have an $(M-1)^{st}$ order deformation of a semigroupal category with
structure map. Observe that each summand 

\[ \underline{\partial}_i 
\omega_{A,B,C,D,E}^{(M)}\]

\noindent essentially
represents the sum
of all composites with total degree $M$ along the three-edge directed path 
minus
the sum of all composites with total degree $M$ along the two-edge directed
path on the boundary of one of the pentagonal faces.  (Here, degree
refers to the power of $\epsilon$ whose coefficient is given by the
composite.)  

This is ``essentially'' the content of each summand, but one must remember
that the context $\lceil \; \rceil$ is not contentless---the summands are
actually composites of the differences just described with various 
structure maps (prolongations of $\alpha^{(0)}$) with the
property that all sources are

\[ [[[A\otimes B]\otimes C]\otimes D]\otimes E \]

\noindent  and all targets are

\[ A\otimes [B\otimes [C\otimes [D\otimes E]]] \, .\]

The odd-index summands correspond to the pentagonal faces on the bottom
of the associahedron as shown in Figure \ref{associahedron}, while the
even-index summands correspond to those on the top.  The square faces of
the associahedron correspond to families of naturality squares, one for
each possible pair of degrees.

In fact, it will suffice to compute 
$[\underline{\partial}_1 + \underline{\partial}_3 + 
\underline{\partial}_5](\omega^{(M)})$:

\begin{lemma} \label{closed.lemma} Suppose for all $N < M$ we have
$\delta(\alpha^{(N)}) + \omega^{(N)} = 0$. Then

\begin{eqnarray*} 
\lefteqn{[\underline{\partial}_1 + \underline{\partial}_3 + 
\underline{\partial}_5](\omega^{(M)})  = }\\
&  & \sum_{\mbox{\scriptsize $\begin{array}{c}
		i + j +k =  M \\
		0\leq i,j,k <  M \end{array}$}}
   \lceil \underline{\partial}_1 \partial_1 \alpha^{(i)} \underline{\partial}_1 
\partial_3 \alpha^{(j) } \underline{\partial}_4 
\partial_3 \alpha^{(k) }\rceil \\
&  & \; - \hspace*{-1.4cm} \sum_{\mbox{\scriptsize $\begin{array}{c}
                   i + j + k + l + m + n  =  M \\
	            0\leq i,j,k,l,m,n  <  M \end{array}$}}  \hspace*{-1.4cm}
   \lceil \underline{\partial}_5 \partial_4 \alpha^{(i)}  \underline{\partial}_5 \partial_2 \alpha^{(j)} \underline{\partial}_5 \partial_0 \alpha^{(k)} \underline{\partial}_3 \partial_2 \alpha^{(l)} \underline{\partial}_3 \partial_0 \alpha^{(m)} \underline{\partial}_0 \partial_0 \alpha^{(n)}\rceil \; .
 \end{eqnarray*}

\end{lemma}

This lemma suffices to complete
the proof of the theorem, since the lemma and calculation by which it is 
derived are precisely
dual to a corresponding statement and derivation concerning 
$[\underline{\partial}_0 + \underline{\partial}_2 + 
\underline{\partial}_4](\omega^{(M)})$. The value derived for this last
expression is

\[ \sum_{\mbox{\scriptsize $\begin{array}{c}
			i + j +k =  M \\
			0\leq i,j,k <  M \end{array}$}}
   \lceil \underline{\partial}_1 \partial_1 \alpha^{(i)} \underline{\partial}_4 
\partial_1 \alpha^{(j) } \underline{\partial}_4 
\partial_3 \alpha^{(k) }\rceil \]
\[ -  \hspace*{-.5cm} \sum_{\mbox{\scriptsize $\begin{array}{c}
                   i + j + k + l + m + n  =  M \\
                   0\leq i,j,k,l,m,n  <  M \end{array}$}} \hspace*{-1.4cm}
   \lceil \underline{\partial}_5 \partial_4 \alpha^{(i)}  \underline{\partial}_2 \partial_4 \alpha^{(j)} \underline{\partial}_2 \partial_2 \alpha^{(k)} \underline{\partial}_0 \partial_4 \alpha^{(l)} \underline{\partial}_0 \partial_2 \alpha^{(m)} \underline{\partial}_0 \partial_0 \alpha^{(n)}\rceil \; .
 \]

Once coincidences of different names for the same map (all of which may
be read off from the associahedron) are taken into account, this expression
differs from that computed in the lemma only in the third and fourth 
factors of the composites in the second summation.  The terms, however,
may be matched one-to-one by swapping the indices $k$ and $l$ into pairs
that are equal by virtue of naturality, thus completing the proof.

Thus, it suffices to prove Lemma \ref{closed.lemma}.  The ambitious reader may 
reconstruct the proof by realizing that the vanishing of earlier obstructions 
is just what is needed to ``fuse'' the summands corresponding to paths round 
two adjacent faces of the associahedron into a similar expression 
corresponding to paths round the union of the two faces. The less ambitious
reader is referred to \cite{Y.book}, where complete details are given.$\Box$
\medskip

We now turn to the question of obstructions for fibred and total deformations
of mon\-oid\-al functors, and for deformations of 
multiplications on mon\-oid\-al categories
(or equivalently, of braided mon\-oid\-al categories). \index{obstructions}
\index{deformation, total} \index{deformation, fibred} \index{deformation
of multiplications on mon\-oid\-al categories}

Since fibred deformations and deformations of multiplications are special cases
of total deformations, defined by restricting the deformation of the target to
be trivial or the deformation of the source to be the diagonal 
deformation induced
by the deformation of the target, respectively, it suffices to consider obstructions
in the case of total deformations.  We begin by giving an explicit formula for
these obstructions, and then show that they are closed.

Recall that the appropriate deformation complex for total deformations of a strong
(or lax) semigroupal
functor
\[ F:{\cal C}\rightarrow {\cal D}, \phi:F(-\otimes -)\Rightarrow F(-)\otimes F(-)\]  

\noindent is 

\[ X_{\bf total}^\bullet (F) = C^\bullet_{\lceil F(p_1) \rceil - \lceil (p_2)_{F^\bullet}\rceil} = X^\bullet(F)\oplus X^{\bullet+1}({\cal C})\oplus X^{\bullet+1}({\cal D}) \]

\noindent with coboundary given by

\[ \left[ \begin{array}{ccc}
\delta_F & 0 & 0 \\
\lceil F(-) \rceil & -\delta_{\cal C} & 0 \\
-\lceil (-)_{F^\bullet}\rceil & 0 & -\delta_{\cal D} 
\end{array} \right] .\]

Thus, a cochain will have coboundary which vanishes in each of the second and
third co\"{o}rdinates if and only if its second and third co\"{o}rdinates are 
cocycles in
$X^{\bullet+1}({\cal C})$ and $ X^{\bullet+1}({\cal D})$, respectively.  
Similarly,
it is easy to see that the obstruction cochain for a total deformation must 
have
as second and third co\"{o}rdinates the obstructions for the 
deformations of the
source and target category, respectively.

Thus, we are left to consider the value of the first co\"{o}rdinate of the
obstruction, and the value of the first co\"{o}rdinate of the coboundary.  
Consider
the hexagonal coherence diagram for oplax mon\-oid\-al functors given in 
Figure \ref{lax.mon.f}, with the maps replaced by their deformed versions.

Calculating the difference of the degree $n$ terms of the two directions around
the diagram gives

\begin{eqnarray*}
\lefteqn{ \sum_{i+j+k=n} 
\lceil a^{(i)}_{F(A),F(B),F(C)}  
[F(A)\otimes \Phi_{B,C}^{(j)} ]\Phi_{A,B\otimes C}^{(k)}\rceil -} \\
& & 
\sum_{i+j+k=n} \lceil  [\Phi_{A,B}^{(i)}\otimes F(C)]\Phi_{A\otimes B,C}^{(j)}
F(\alpha^{(k)}_{A,B,C}) \rceil ,
 \end{eqnarray*}

\noindent where $\alpha$ and $a$ are the associators for $\cal C$ and $\cal D$,
respectively.  This must vanish for $n = 1$ for first order 
total deformations:  the vanishing 
is simply the cocycle condition in $X_{\bf total}^\bullet (F)$.  
For a deformation
to extend to an $N^{th}$ order deformation this quantity must vanish for
all $n \leq N$, and indeed in addition to the vanishing of the corresponding
second and third co\"{o}rdinates, this condition is sufficient.  
Separating out the
terms in which the index $^{(n)}$ occurs, we find that the vanishing 
conditions are precisely the condition that 
$[\phi^{(n)}, \alpha^{(n)}, a^{(n)}]$ 
cobounds $[\Omega^{(n)}, \omega^{(n)}, o^{(n)}]$, where

\begin{eqnarray*}
\Omega^{(n)}  &= &  
\sum_{\mbox{\scriptsize $\begin{array}{c} i+j+k = n \\ i,j,k < n 
\end{array}$} } \lceil a^{(i)}_{F(A),F(B),F(C)}  
[F(A)\otimes \Phi_{B,C}^{(j)} ]\Phi_{A,B\otimes C}^{(k)}\rceil \\
& & \; - 
\sum_{\mbox{\scriptsize $\begin{array}{c} i+j+k = n \\ i,j,k < n 
\end{array}$} } \lceil  [\Phi_{A,B}^{(i)}\otimes F(C)]\Phi_{A\otimes B,C}^{(j)}
F(\alpha^{(k)}_{A,B,C}) \rceil \end{eqnarray*}

\noindent and $\omega^{(n)}$ and $o^{(n)}$ are the obstructions to the
extension of the deformations of the source and target categories, 
respectively.

\begin{figure}[htb]
\epsfig{file=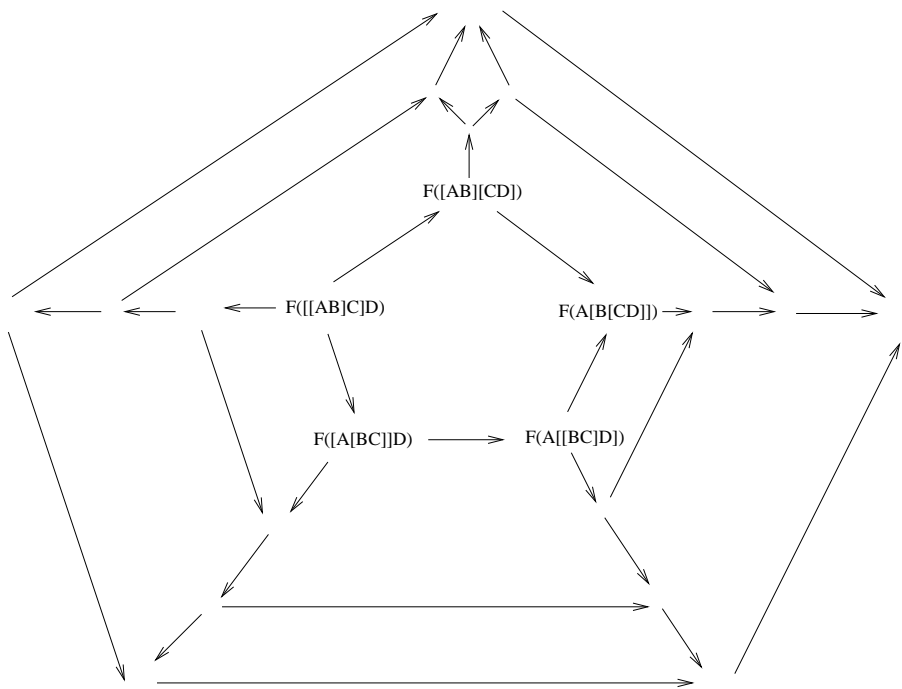,width=4.7in}
\caption{\label{Chinese.lantern} The ``Chinese Lantern''}
\end{figure}

All that remains to show is that the first co\"{o}rdinate of the coboundary of
$[\Omega^{(n)}, \omega^{(n)}, o^{(n)}]$ vanishes.  We leave the details of the
proof to the reader.  The method is identical to that applied in the case of 
the obstructions for deformations of a semigroupal category, except that the
associahedron must be replaced with the diagram given in 
Figure \ref{Chinese.lantern}.  In Figure \ref{Chinese.lantern} 
we have suppressed all object
and arrow labels except for the objects on the inner pentagon which are 
written with the null infix in place of $\otimes$ to save space.\footnote{The 
diagram of Figure \ref{Chinese.lantern} was given the
name ``the Chinese lantern'' due to its resemblance to a paper lantern when
it is drawn in perspective as the edges of a 3-dimensional polytope with
the innermost and outermost pentagons a parallel horizontal faces. The diagram
had occurred previously in work of Stasheff, who did not name it, and has
been dubbed the ``multiplicahedron'' by other recent researchers.}  
The labels
can be recovered by labeling all radial maps with prolongations of $\Phi$
and all maps parallel to those between the labeled objects with prolongations
of (functorial images of) $\alpha$.  All hexagons are prolongations of
the coherence hexagon for semigroupal functors, and all squares, except the
diamond-shaped one in the top center, are naturality squares.  The
diamond is a functoriality square. \index{Chinese lantern (commutative
diagram)}

\clearpage 
 \section{Units} \label{unit.chapter}
\vspace*{1cm}

Thus far we have discussed deformations of semigroupal 
categories and
functors. In \cite{Y.book} it is shown that there is very little 
additional content to the deformation theory of monoidal categories
and functors.  We recall the relevant results form \cite{Y.book}, and
refer the reader there for proofs:

\begin{thm} \index{deformation of a mon\-oid\-al category}
Every semigroupal deformation of a mon\-oid\-al category 
\[({\cal C},\otimes, I, \alpha,
\rho, \lambda )\]

\noindent becomes a mon\-oid\-al category when equip\-ped with 
unit transformations
$\tilde{\rho}$ and $\tilde{\lambda}$ given by 

\[ \tilde{\nu} = \sum_i \nu^{(i)} \epsilon^i \]

\[ \tilde{\lambda}_A = \sum_i \lambda_A^{(i)} \epsilon^i \]

\[ \tilde{\rho}_A = \sum_i \rho_A^{(i)} \epsilon^i,\]

\noindent where

\[ \lambda_A^{(n)} = \sum_{i+j = n} \lceil \alpha_{A,I,I}^{(i)} 
[A\otimes \nu^{(j)}] \rceil 
\]

\[ \rho_B^{(n)} = 
\sum_{i+j = n} \lceil \beta_{I,I,B}^{(i)} [\nu^{(j)}\otimes B] \rceil,\]

\noindent $\nu^{(i)}:I\rightarrow I$ is any family of maps satisfying
$\nu^{(0)} = \lambda_I = \rho_I$,
and $\tilde{\alpha}^{-1} = \sum_i \beta^{(i)} \epsilon^i$.
\end{thm}

For strong mon\-oid\-al functors, we have  

\begin{thm} \index{deformation} \index{monoidal functor, strong}
\index{functor, strong monoidal} \index{strong monoidal functor}
If $(F:{\cal C}\rightarrow {\cal D}, \Phi, F_o)$ is a 
strong mon\-oid\-al functor, then every 
semigroupal deformation of $F$ extends uniquely to a deformation as a mon\-oid\-al 
functor.
\end{thm}

The final conditions involving units 
are the conditions
in the definition of a multiplication on a mon\-oid\-al category.  
This condition, however,
is trivially satisfied by any deformation, 
since the isomorphism giving the structure
of the given multiplication still provides 
the necessary structure after deformation.

\clearpage 
  \section{Modules over Monoidal Categories and Functors}
\vspace*{1cm}

The key to the construction of the deformation complexes for monoidal
categories and monoidal functors was the fact that both can serve in
different ways as analogues of rings.  In order to move the deformation
complexes of Section \ref{def.rev} into the realm of classical homological 
algebra, we need to consider the analogues of modules in each case.  
In the case of categories, the author introduced the required notion
in \cite{Y.TiPTiC}.

\begin{defin}
If $\cal C$ and $\cal D$ are monoidal categories, {\em a (strong) left
$\cal C$-module $\cal M$} is a category equipped with a functor
$\triangleright:{\cal C} \times {\cal M} \rightarrow {\cal M}$ (written as an
infix) and
natural isomorphisms 

\[ a_l:(A\otimes B)\triangleright X \rightarrow A\triangleright 
(B\triangleright X) \] 

\noindent and 

\[ \ell:I\triangleright X \rightarrow X \] 

\noindent satisfying

\begin{center}
\setlength{\unitlength}{3947sp}%
\begingroup\makeatletter\ifx\SetFigFont\undefined
\def\x#1#2#3#4#5#6#7\relax{\def\x{#1#2#3#4#5#6}}%
\expandafter\x\fmtname xxxxxx\relax \def\y{splain}%
\ifx\x\y   
\gdef\SetFigFont#1#2#3{%
  \ifnum #1<17\tiny\else \ifnum #1<20\small\else
  \ifnum #1<24\normalsize\else \ifnum #1<29\large\else
  \ifnum #1<34\Large\else \ifnum #1<41\LARGE\else
     \huge\fi\fi\fi\fi\fi\fi
  \csname #3\endcsname}%
\else
\gdef\SetFigFont#1#2#3{\begingroup
  \count@#1\relax \ifnum 25<\count@\count@25\fi
  \def\x{\endgroup\@setsize\SetFigFont{#2pt}}%
  \expandafter\x
    \csname \romannumeral\the\count@ pt\expandafter\endcsname
    \csname @\romannumeral\the\count@ pt\endcsname
  \csname #3\endcsname}%
\fi
\fi\endgroup
\begin{picture}(4287,6180)(76,-5386)
\thinlines
\put(1876,-4036){\vector( 1, 0){1125}}
\put(1426,-4186){\vector( 3,-4){675}}
\put(3226,-4261){\vector(-3,-4){675}}
\put(901,-4111){\makebox(0,0)[lb]{\smash{\SetFigFont{10}{14.4}{rm}$(A \otimes I) \triangleright X$}}}
\put(3151,-4111){\makebox(0,0)[lb]{\smash{\SetFigFont{10}{14.4}{rm}$A \triangleright (I \triangleright X)$}}}
\put(2251,-3961){\makebox(0,0)[lb]{\smash{\SetFigFont{10}{14.4}{rm}$a_l$}}}
\put(1126,-4861){\makebox(0,0)[lb]{\smash{\SetFigFont{10}{14.4}{rm}$\rho \triangleright X$}}}
\put(3151,-4861){\makebox(0,0)[lb]{\smash{\SetFigFont{10}{14.4}{rm}$A \triangleright \ell$}}}
\put(2101,-5386){\makebox(0,0)[lb]{\smash{\SetFigFont{10}{14.4}{rm}$A \triangleright X$}}}
\put(901,-886){\vector( 1,-2){450}}
\put(2251,-2086){\vector( 1, 0){900}}
\put(3751,-1861){\vector( 2, 3){600}}
\put(3076,464){\vector( 1,-1){900}}
\put(1051,-436){\vector( 1, 1){900}}
\put(601,-2161){\makebox(0,0)[lb]{\smash{\SetFigFont{10}{14.4}{rm}$(A \otimes (B \otimes C)) \triangleright X$}}}
\put(1126,-61){\makebox(0,0)[lb]{\smash{\SetFigFont{10}{14.4}{rm}$a_l$}}}
\put(526,-1411){\makebox(0,0)[lb]{\smash{\SetFigFont{10}{14.4}{rm}$\alpha \otimes X$}}}
\put(2626,-2311){\makebox(0,0)[lb]{\smash{\SetFigFont{10}{14.4}{rm}$a_l$}}}
\put( 76,-736){\makebox(0,0)[lb]{\smash{\SetFigFont{10}{14.4}{rm}$((A \otimes B) \otimes C) \triangleright X$}}}
\put(3526,-736){\makebox(0,0)[lb]{\smash{\SetFigFont{10}{14.4}{rm}$A \triangleright (B \triangleright (C \triangleright X))$}}}
\put(3301,-2161){\makebox(0,0)[lb]{\smash{\SetFigFont{10}{14.4}{rm}$A \triangleright ((B \otimes C) \triangleright X)$}}}
\put(4276,-1486){\makebox(0,0)[lb]{\smash{\SetFigFont{10}{14.4}{rm}$A \triangleright a_l$}}}
\put(1876,614){\makebox(0,0)[lb]{\smash{\SetFigFont{10}{14.4}{rm}$(A \otimes B) \triangleright (C \triangleright X)$}}}
\put(3901, 14){\makebox(0,0)[lb]{\smash{\SetFigFont{10}{14.4}{rm}$a_l$}}}
\end{picture}
\end{center}

{\em A (strong) right $\cal D$-module $\cal M$} is a category equipped
wtih a functor $\triangleleft:{\cal M}\times {\cal C}\rightarrow {\cal M}$ (also
as an infix) and 
natural isomorphisms $a_r:(X\triangleleft C)\triangleleft D \rightarrow X\triangleleft (C\otimes D)$
and $r:X\triangleleft I \rightarrow X$ satisfying

\begin{center}
\setlength{\unitlength}{3947sp}%
\begingroup\makeatletter\ifx\SetFigFont\undefined
\def\x#1#2#3#4#5#6#7\relax{\def\x{#1#2#3#4#5#6}}%
\expandafter\x\fmtname xxxxxx\relax \def\y{splain}%
\ifx\x\y   
\gdef\SetFigFont#1#2#3{%
  \ifnum #1<17\tiny\else \ifnum #1<20\small\else
  \ifnum #1<24\normalsize\else \ifnum #1<29\large\else
  \ifnum #1<34\Large\else \ifnum #1<41\LARGE\else
     \huge\fi\fi\fi\fi\fi\fi
  \csname #3\endcsname}%
\else
\gdef\SetFigFont#1#2#3{\begingroup
  \count@#1\relax \ifnum 25<\count@\count@25\fi
  \def\x{\endgroup\@setsize\SetFigFont{#2pt}}%
  \expandafter\x
    \csname \romannumeral\the\count@ pt\expandafter\endcsname
    \csname @\romannumeral\the\count@ pt\endcsname
  \csname #3\endcsname}%
\fi
\fi\endgroup
\begin{picture}(4287,6180)(76,-5386)
\thinlines
\put(1876,-4036){\vector( 1, 0){1125}}
\put(1426,-4186){\vector( 3,-4){675}}
\put(3226,-4261){\vector(-3,-4){675}}
\put(901,-4111){\makebox(0,0)[lb]{\smash{\SetFigFont{10}{14.4}{rm}$(Y \triangleleft I) \triangleleft B$}}}
\put(3151,-4111){\makebox(0,0)[lb]{\smash{\SetFigFont{10}{14.4}{rm}$Y \triangleleft (I \otimes B)$}}}
\put(2251,-3961){\makebox(0,0)[lb]{\smash{\SetFigFont{10}{14.4}{rm}$a_r$}}}
\put(1126,-4861){\makebox(0,0)[lb]{\smash{\SetFigFont{10}{14.4}{rm}$r \triangleleft B$}}}
\put(3151,-4861){\makebox(0,0)[lb]{\smash{\SetFigFont{10}{14.4}{rm}$Y \triangleleft \lambda$}}}
\put(2101,-5386){\makebox(0,0)[lb]{\smash{\SetFigFont{10}{14.4}{rm}$Y \triangleleft B$}}}
\put(901,-886){\vector( 1,-2){450}}
\put(2251,-2086){\vector( 1, 0){900}}
\put(3751,-1861){\vector( 2, 3){600}}
\put(3076,464){\vector( 1,-1){900}}
\put(1051,-436){\vector( 1, 1){900}}
\put(601,-2161){\makebox(0,0)[lb]{\smash{\SetFigFont{10}{14.4}{rm}$(Y \triangleleft (B \otimes C)) \otimes D$}}}
\put(1126,-61){\makebox(0,0)[lb]{\smash{\SetFigFont{10}{14.4}{rm}$a_r$}}}
\put(526,-1411){\makebox(0,0)[lb]{\smash{\SetFigFont{10}{14.4}{rm}$a_r \triangleleft D$}}}
\put(2626,-2311){\makebox(0,0)[lb]{\smash{\SetFigFont{10}{14.4}{rm}$a_r$}}}
\put( 76,-736){\makebox(0,0)[lb]{\smash{\SetFigFont{10}{14.4}{rm}$((Y \triangleleft B) \triangleleft C) \triangleleft D$}}}
\put(3526,-736){\makebox(0,0)[lb]{\smash{\SetFigFont{10}{14.4}{rm}$Y \triangleleft (B \otimes (C \otimes D))$}}}
\put(3301,-2161){\makebox(0,0)[lb]{\smash{\SetFigFont{10}{14.4}{rm}$Y \triangleleft ((B \otimes C) \otimes D)$}}}
\put(4276,-1486){\makebox(0,0)[lb]{\smash{\SetFigFont{10}{14.4}{rm}$Y \triangleleft \alpha$}}}
\put(1876,614){\makebox(0,0)[lb]{\smash{\SetFigFont{10}{14.4}{rm}$(Y \triangleleft B) \triangleleft (C \otimes D)$}}}
\put(3901, 14){\makebox(0,0)[lb]{\smash{\SetFigFont{10}{14.4}{rm}$a_r$}}}
\end{picture}
\end{center}

{\em A (strong) $\cal C$, $\cal D$-bimodule} is a category $\cal M$
equipped with both a left $\cal C$-module structure and a right
$\cal D$-module structure and a natural isomorphism
$a_m:(A \triangleright X)\triangleleft C \rightarrow A\triangleright (X \triangleleft C)$ and satisfying

\begin{center}

\setlength{\unitlength}{3947sp}%
\begingroup\makeatletter\ifx\SetFigFont\undefined
\def\x#1#2#3#4#5#6#7\relax{\def\x{#1#2#3#4#5#6}}%
\expandafter\x\fmtname xxxxxx\relax \def\y{splain}%
\ifx\x\y   
\gdef\SetFigFont#1#2#3{%
  \ifnum #1<17\tiny\else \ifnum #1<20\small\else
  \ifnum #1<24\normalsize\else \ifnum #1<29\large\else
  \ifnum #1<34\Large\else \ifnum #1<41\LARGE\else
     \huge\fi\fi\fi\fi\fi\fi
  \csname #3\endcsname}%
\else
\gdef\SetFigFont#1#2#3{\begingroup
  \count@#1\relax \ifnum 25<\count@\count@25\fi
  \def\x{\endgroup\@setsize\SetFigFont{#2pt}}%
  \expandafter\x
    \csname \romannumeral\the\count@ pt\expandafter\endcsname
    \csname @\romannumeral\the\count@ pt\endcsname
  \csname #3\endcsname}%
\fi
\fi\endgroup
\begin{picture}(4362,6705)(1,-5911)
\thinlines
\put(901,-886){\vector( 1,-2){450}}
\put(2251,-2086){\vector( 1, 0){900}}
\put(3751,-1861){\vector( 2, 3){600}}
\put(3076,464){\vector( 1,-1){900}}
\put(1051,-436){\vector( 1, 1){900}}
\put(601,-2161){\makebox(0,0)[lb]{\smash{\SetFigFont{10}{14.4}{rm}
$(A \triangleright (B \triangleright X)) \triangleleft D$}}}
\put( 76,-736){\makebox(0,0)[lb]{\smash{\SetFigFont{10}{14.4}{rm}$((A \otimes B) \triangleright X) \triangleleft D$}}}
\put(3526,-736){\makebox(0,0)[lb]{\smash{\SetFigFont{10}{14.4}{rm}$A \triangleright (B \triangleright (X \triangleleft D))$}}}
\put(3301,-2161){\makebox(0,0)[lb]{\smash{\SetFigFont{10}{14.4}{rm}$A \triangleright ((B \triangleright X) \triangleleft D)$}}}
\put(1876,614){\makebox(0,0)[lb]{\smash{\SetFigFont{10}{14.4}{rm}$(A \otimes B) \triangleright (X \triangleleft D)$}}}
\put(3901, 14){\makebox(0,0)[lb]{\smash{\SetFigFont{10}{14.4}{rm}$a_l$}}}
\put(4276,-1486){\makebox(0,0)[lb]{\smash{\SetFigFont{10}{14.4}{rm}$A \triangleright a_m$}}}
\put(2626,-2311){\makebox(0,0)[lb]{\smash{\SetFigFont{10}{14.4}{rm}$a_m$}}}
\put(451,-1411){\makebox(0,0)[lb]{\smash{\SetFigFont{10}{14.4}{rm}$a_l \triangleleft D$}}}
\put(1051,-61){\makebox(0,0)[lb]{\smash{\SetFigFont{10}{14.4}{rm}$a_m$}}}
\put(826,-4486){\vector( 1,-2){450}}
\put(2176,-5686){\vector( 1, 0){900}}
\put(3676,-5461){\vector( 2, 3){600}}
\put(3001,-3136){\vector( 1,-1){900}}
\put(976,-4036){\vector( 1, 1){900}}
\put(526,-5761){\makebox(0,0)[lb]{\smash{\SetFigFont{10}{14.4}{rm}$(A \triangleright (X \triangleleft C)) \triangleleft D$}}}
\put(  1,-4336){\makebox(0,0)[lb]{\smash{\SetFigFont{10}{14.4}{rm}$((A \triangleright X) \triangleleft C) \triangleleft D$}}}
\put(3451,-4336){\makebox(0,0)[lb]{\smash{\SetFigFont{10}{14.4}{rm}$A \triangleright (X \triangleleft (C \otimes D))$}}}
\put(3226,-5761){\makebox(0,0)[lb]{\smash{\SetFigFont{10}{14.4}{rm}$A \triangleright ((X \triangleleft C) \triangleleft D)$}}}
\put(1801,-2986){\makebox(0,0)[lb]{\smash{\SetFigFont{10}{14.4}{rm}$(A \triangleright X) \triangleleft (C \otimes D)$}}}
\put(1051,-3661){\makebox(0,0)[lb]{\smash{\SetFigFont{10}{14.4}{rm}$a_r$}}}
\put(3826,-3586){\makebox(0,0)[lb]{\smash{\SetFigFont{10}{14.4}{rm}$a_m$}}}
\put(4201,-5086){\makebox(0,0)[lb]{\smash{\SetFigFont{10}{14.4}{rm}$A \triangleright a_r$}}}
\put(2551,-5911){\makebox(0,0)[lb]{\smash{\SetFigFont{10}{14.4}{rm}$a_m$}}}
\put(301,-5011){\makebox(0,0)[lb]{\smash{\SetFigFont{10}{14.4}{rm}$a_m \triangleleft D$}}}
\end{picture}
\end{center}

\end{defin}

We will be primarily concerned with the special case where 
${\cal C} = {\cal D} = {\cal M}$ with $\triangleright = \triangleleft = \otimes$
and the obvious identification of structure maps with the structure maps
of the monoidal structure, but most of our results hold more generally.

Other examples may be found in \cite{Y.TiPTiC}.

More important for our purposes are the corresponding notions for (lax)
monoidal functors:

\begin{defin}
Suppose $\cal C$, $\cal D$, and $\cal E$ are monoidal categories and
that 

\[ (F,\tilde{F},F_0):{\cal C}\rightarrow {\cal E}\] 

\noindent and
\[ (G,\tilde{G},G_0):{\cal D}\rightarrow {\cal E}\] 
\noindent are monoidal functors.
Suppose, moreover, that 
$\cal M$ is a left $\cal C$-module (resp. right $\cal D$-module,
$\cal C$, $\cal D$-bimodule). {\em  A left $F$-module over $\cal M$
(resp. a right $G$-module over $\cal M$, an $F, G$-bimodule over $\cal M$)}
is then a functor $M:{\cal M}\rightarrow {\cal E}$ equipped with
natural transformations $\mu_l:F(-)\otimes M(-) \Rightarrow M(-\triangleright - )$
(resp. $\mu_r:M(-) \otimes G(-) \Rightarrow M(- \triangleleft -)$, both 
$\mu_l$ and $\mu_r$) and satisfying

\begin{center}
\setlength{\unitlength}{3947sp}%
\begingroup\makeatletter\ifx\SetFigFont\undefined
\def\x#1#2#3#4#5#6#7\relax{\def\x{#1#2#3#4#5#6}}%
\expandafter\x\fmtname xxxxxx\relax \def\y{splain}%
\ifx\x\y   
\gdef\SetFigFont#1#2#3{%
  \ifnum #1<17\tiny\else \ifnum #1<20\small\else
  \ifnum #1<24\normalsize\else \ifnum #1<29\large\else
  \ifnum #1<34\Large\else \ifnum #1<41\LARGE\else
     \huge\fi\fi\fi\fi\fi\fi
  \csname #3\endcsname}%
\else
\gdef\SetFigFont#1#2#3{\begingroup
  \count@#1\relax \ifnum 25<\count@\count@25\fi
  \def\x{\endgroup\@setsize\SetFigFont{#2pt}}%
  \expandafter\x
    \csname \romannumeral\the\count@ pt\expandafter\endcsname
    \csname @\romannumeral\the\count@ pt\endcsname
  \csname #3\endcsname}%
\fi
\fi\endgroup
\begin{picture}(4275,5640)(526,-5086)
\thinlines
\put(2491,329){\vector( 1, 0){975}}
\put(4441,119){\vector( 0,-1){630}}
\put(4426,-976){\vector( 0,-1){720}}
\put(2176,-1951){\vector( 1, 0){1410}}
\put(1246,-961){\vector( 0,-1){720}}
\put(1231,134){\vector( 0,-1){630}}
\put(601,239){\makebox(0,0)[lb]{\smash{\SetFigFont{10}{14.4}{rm}$[F(A) \otimes F(B)] \otimes M(X)$}}}
\put(3601,239){\makebox(0,0)[lb]{\smash{\SetFigFont{10}{14.4}{rm}$F(A) \otimes [F(B) \otimes M(X)]$}}}
\put(751,-811){\makebox(0,0)[lb]{\smash{\SetFigFont{10}{14.4}{rm}$F(A \otimes B) \otimes M(X)$}}}
\put(3901,-811){\makebox(0,0)[lb]{\smash{\SetFigFont{10}{14.4}{rm}$F(A) \otimes M(B \triangleright X) $}}}
\put(3901,-2011){\makebox(0,0)[lb]{\smash{\SetFigFont{10}{14.4}{rm}$M(A \triangleright [B \triangleright X])$}}}
\put(751,-2011){\makebox(0,0)[lb]{\smash{\SetFigFont{10}{14.4}{rm}$M([A \otimes B] \triangleright X)$}}}
\put(2776,464){\makebox(0,0)[lb]{\smash{\SetFigFont{10}{14.4}{rm}$\alpha$}}}
\put(2476,-2236){\makebox(0,0)[lb]{\smash{\SetFigFont{10}{14.4}{rm}$M(a_l)$}}}
\put(4726,-211){\makebox(0,0)[lb]{\smash{\SetFigFont{10}{14.4}{rm}$Id \otimes \mu_l$}}}
\put(4801,-1411){\makebox(0,0)[lb]{\smash{\SetFigFont{10}{14.4}{rm}$\mu_l$}}}
\put(676,-1336){\makebox(0,0)[lb]{\smash{\SetFigFont{10}{14.4}{rm}$\mu_l$}}}
\put(526,-211){\makebox(0,0)[lb]{\smash{\SetFigFont{10}{14.4}{rm}$\tilde{F} \otimes Id$}}}
\put(2401,-3586){\vector( 1, 0){825}}
\put(3751,-3811){\vector( 0,-1){675}}
\put(3376,-4711){\vector(-1, 0){1200}}
\put(1876,-3811){\vector( 0,-1){675}}
\put(1576,-3661){\makebox(0,0)[lb]{\smash{\SetFigFont{10}{14.4}{rm}$I \otimes M(X)$}}}
\put(3376,-3661){\makebox(0,0)[lb]{\smash{\SetFigFont{10}{14.4}{rm}$F(I) \otimes M(X)$}}}
\put(3451,-4786){\makebox(0,0)[lb]{\smash{\SetFigFont{10}{14.4}{rm}$M(I \otimes X)$}}}
\put(1651,-4786){\makebox(0,0)[lb]{\smash{\SetFigFont{10}{14.4}{rm}$M(X)$}}}
\put(1426,-4186){\makebox(0,0)[lb]{\smash{\SetFigFont{10}{14.4}{rm}$\lambda$}}}
\put(3976,-4186){\makebox(0,0)[lb]{\smash{\SetFigFont{10}{14.4}{rm}$\mu_l$}}}
\put(2476,-5086){\makebox(0,0)[lb]{\smash{\SetFigFont{10}{14.4}{rm}$M(\lambda)$}}}
\put(2476,-3361){\makebox(0,0)[lb]{\smash{\SetFigFont{10}{14.4}{rm}$F_0 \otimes Id$}}}
\end{picture}
\end{center}

(resp.  a similar hexagon and square for the right action, 
all four of these and

\begin{center}
\setlength{\unitlength}{3947sp}%
\begingroup\makeatletter\ifx\SetFigFont\undefined
\def\x#1#2#3#4#5#6#7\relax{\def\x{#1#2#3#4#5#6}}%
\expandafter\x\fmtname xxxxxx\relax \def\y{splain}%
\ifx\x\y   
\gdef\SetFigFont#1#2#3{%
  \ifnum #1<17\tiny\else \ifnum #1<20\small\else
  \ifnum #1<24\normalsize\else \ifnum #1<29\large\else
  \ifnum #1<34\Large\else \ifnum #1<41\LARGE\else
     \huge\fi\fi\fi\fi\fi\fi
  \csname #3\endcsname}%
\else
\gdef\SetFigFont#1#2#3{\begingroup
  \count@#1\relax \ifnum 25<\count@\count@25\fi
  \def\x{\endgroup\@setsize\SetFigFont{#2pt}}%
  \expandafter\x
    \csname \romannumeral\the\count@ pt\expandafter\endcsname
    \csname @\romannumeral\the\count@ pt\endcsname
  \csname #3\endcsname}%
\fi
\fi\endgroup
\begin{picture}(4575,2790)(376,-2311)
\thinlines
\put(2476,164){\vector( 1, 0){1125}}
\put(4651,-61){\vector( 0,-1){750}}
\put(4651,-1186){\vector( 0,-1){675}}
\put(1276,-61){\vector( 0,-1){675}}
\put(1276,-1186){\vector( 0,-1){600}}
\put(2101,-2011){\vector( 1, 0){1725}}
\put(751,-1036){\makebox(0,0)[lb]{\smash{\SetFigFont{10}{14.4}{rm}$M(A \triangleright X) \otimes G(B)$}}}
\put(4051,-1036){\makebox(0,0)[lb]{\smash{\SetFigFont{10}{14.4}{rm}$F(A) \otimes M(X \triangleleft B)$}}}
\put(751,-2086){\makebox(0,0)[lb]{\smash{\SetFigFont{10}{14.4}{rm}$M([A \triangleright X] \triangleleft B)$}}}
\put(601, 89){\makebox(0,0)[lb]{\smash{\SetFigFont{10}{14.4}{rm}$[F(A) \otimes M(X)] \otimes G(B)$}}}
\put(3751, 89){\makebox(0,0)[lb]{\smash{\SetFigFont{10}{14.4}{rm}$F(A) \otimes [M(X) \otimes G(B)]$}}}
\put(4126,-2086){\makebox(0,0)[lb]{\smash{\SetFigFont{10}{14.4}{rm}$M(A \triangleright [X \triangleleft B])$\,\,\,\, {\em .)} }}}
\put(4876,-436){\makebox(0,0)[lb]{\smash{\SetFigFont{10}{14.4}{rm}$Id \otimes \mu_r$}}}
\put(4951,-1561){\makebox(0,0)[lb]{\smash{\SetFigFont{10}{14.4}{rm}$\mu_l$}}}
\put(376,-361){\makebox(0,0)[lb]{\smash{\SetFigFont{10}{14.4}{rm}$\mu_l \otimes Id$}}}
\put(676,-1561){\makebox(0,0)[lb]{\smash{\SetFigFont{10}{14.4}{rm}$\mu_r$}}}
\put(2551,-2311){\makebox(0,0)[lb]{\smash{\SetFigFont{10}{14.4}{rm}$M(a_m)$}}}
\put(2926,389){\makebox(0,0)[lb]{\smash{\SetFigFont{10}{14.4}{rm}$\alpha$}}}
\end{picture}

\end{center}

\end{defin}

\begin{defin}
{\em A left (resp. right, bi-) module map} from $M$ to $N$
is a natural transformation $f:M\Rightarrow N$ satisfying

\begin{center}
\setlength{\unitlength}{3947sp}%
\begingroup\makeatletter\ifx\SetFigFont\undefined%
\gdef\SetFigFont#1#2#3#4#5{%
  \reset@font\fontsize{#1}{#2pt}%
  \fontfamily{#3}\fontseries{#4}\fontshape{#5}%
  \selectfont}%
\fi\endgroup%
\begin{picture}(3300,2010)(76,-1711)
\thinlines
\put(1726,-1486){\vector( 1, 0){1050}}
\put(3226,-211){\vector( 0,-1){1050}}
\put(1051,-211){\vector( 0,-1){1050}}
\put(1726, 14){\vector( 1, 0){1050}}
\put(601,-61){\makebox(0,0)[lb]{\smash{\SetFigFont{10}{14.4}{rm}$F(A) \otimes M(X)$}}}
\put(601,-1561){\makebox(0,0)[lb]{\smash{\SetFigFont{10}{14.4}{rm}$F(A) \otimes N(X)$}}}
\put(2926,-1561){\makebox(0,0)[lb]{\smash{\SetFigFont{10}{14.4}{rm}$N(A \triangleright X)$}}}
\put(2926,-61){\makebox(0,0)[lb]{\smash{\SetFigFont{10}{14.4}{rm}$M(A \triangleright X)$}}}
\put(2101,164){\makebox(0,0)[lb]{\smash{\SetFigFont{10}{14.4}{rm}$\mu_l$}}}
\put(2026,-1711){\makebox(0,0)[lb]{\smash{\SetFigFont{10}{14.4}{rm}$\nu_l$}}}
\put(3376,-736){\makebox(0,0)[lb]{\smash{\SetFigFont{10}{14.4}{rm}$f_{A \triangleright X}$}}}
\put( 76,-736){\makebox(0,0)[lb]{\smash{\SetFigFont{10}{14.4}{rm}$F(A) \otimes f_X$}}}
\end{picture}
\end{center}

(resp.  as similar square for the right action, both)
\end{defin}

Henceforth by abuse of notation, we will denote monoidal functors
by the name of the functor only.

There are a number of obvious examples, and some not-so-obvious:

First, it is clear that any monoidal functor $F:{\cal C}\rightarrow {\cal D}$ 
equipped with $\mu_l = \mu_r = \tilde{F}$ is an $F,F$-bimodule over $\cal C$.

Second, suppose we have a second monoidal functor 
$G:{\cal C}\rightarrow {\cal D}$ and a natural transformation 
$\phi:F\Rightarrow G$.  $G$ then becomes an $F,F$-bimodule over $\cal C$, when
equipped with the structure $\mu_l = \tilde{G}(\phi \otimes Id_G)$ and
$\mu_r = \tilde{G}(Id_G \otimes \phi)$.

We now turn to the main theorem of this section:

\begin{thm} \label{abelian}
If $\cal E$ is an abelian category equipped with a monoidal structure,
and $F:{\cal C}\rightarrow {\cal E}$ and $G:{\cal D}\rightarrow {\cal E}$
are monoidal functors such that for all $A \in Ob({\cal C})$
$F(A)\otimes -$ is exact (resp. for all $C \in Ob({\cal D})$
$-\otimes G(C)$ is exact, both), then for any left $\cal C$-module
(resp. right $\cal D$-module, $\cal C$,$\cal D$-bimodule) $\cal M$,
the category of left $F$-modules (resp. right $G$-modules, $F,G$-bimodules)
over $\cal M$ is an abelian category.
\end{thm}

\noindent{\bf proof:}  We proceed by showing that the forgetful functor
to the functor category ${\cal E}^{\cal M}$ induces an additive structure
and all necessary limits and colimits, then verify the additional
conditions for abelianness in the form in terms of the ``parallel'' of
a map as given in Popescu \cite{Pop}.
 
Now, for the null object, we can take the constant functor
$0$, since by the exactness hypotheses  $F(A)\otimes -$ and 
$-\otimes g(C)$ preserve $0$, and thus the zero map is the
unique action.  Initial and terminal conditions all follow from the
uniqueness of the zero map and consideration of the diagrams which
assert that it is a left $F$-module map (resp. right $G$-module map, both).

For kernels, consider first the case of a left $F$-module map 
$f:M\Rightarrow N$, we claim
that the kernel in ${\cal E}^{\cal M}$ has a unique left $F$-module structure
such that the inclusion is an $F$-module map, as are all canonical
maps induced by $F$-module maps annihilated by post-composition with $f$.

Now, for each object $X \in Ob({\cal M})$ we have the exact
sequence

\[ 0\rightarrow ker(f_X)\rightarrow M(C) \rightarrow N(C) \]

\noindent but $F(A)\otimes -$ is exact, and $M$ and $N$ are left
$F$-modules.  Thus we have a commutative diagram with exact rows
\bigskip

\begin{center}
\setlength{\unitlength}{3947sp}%
\begingroup\makeatletter\ifx\SetFigFont\undefined%
\gdef\SetFigFont#1#2#3#4#5{%
  \reset@font\fontsize{#1}{#2pt}%
  \fontfamily{#3}\fontseries{#4}\fontshape{#5}%
  \selectfont}%
\fi\endgroup%
\begin{picture}(4875,1380)(1201,-1561)
\thinlines
\put(1426,-286){\vector( 1, 0){225}}
\put(3076,-286){\vector( 1, 0){375}}
\put(4801,-286){\vector( 1, 0){450}}
\put(4426,-1486){\vector( 1, 0){1050}}
\put(3001,-1486){\vector( 1, 0){450}}
\put(3976,-511){\vector( 0,-1){750}}
\put(5926,-511){\vector( 0,-1){750}}
\put(1426,-1486){\vector( 1, 0){525}}
\put(2401,-1261){\vector( 0,-1){0}}
\multiput(2401,-511)(0.00000,-9.03614){83}{\makebox(1.6667,11.6667){\SetFigFont{5}{6}{rm}.}}
\put(1201,-361){\makebox(0,0)[lb]{\smash{\SetFigFont{10}{14.4}{rm}$0$}}}
\put(1801,-361){\makebox(0,0)[lb]{\smash{\SetFigFont{10}{14.4}{rm}$F(A ) \otimes ker(f_X)$}}}
\put(3601,-361){\makebox(0,0)[lb]{\smash{\SetFigFont{10}{14.4}{rm}$F(A) \otimes M(X)$}}}
\put(5401,-361){\makebox(0,0)[lb]{\smash{\SetFigFont{10}{14.4}{rm}$F(A) \otimes N(X)$}}}
\put(1201,-1561){\makebox(0,0)[lb]{\smash{\SetFigFont{10}{14.4}{rm}$0$}}}
\put(2101,-1561){\makebox(0,0)[lb]{\smash{\SetFigFont{10}{14.4}{rm}$ker(f_{A \triangleright X})$}}}
\put(3601,-1561){\makebox(0,0)[lb]{\smash{\SetFigFont{10}{14.4}{rm}$M(A \triangleright X)$}}}
\put(5701,-1561){\makebox(0,0)[lb]{\smash{\SetFigFont{10}{14.4}{rm}$N(A \triangleright X)$}}}
\put(4201,-886){\makebox(0,0)[lb]{\smash{\SetFigFont{10}{14.4}{rm}$\mu_l$}}}
\put(6076,-886){\makebox(0,0)[lb]{\smash{\SetFigFont{10}{14.4}{rm}$\nu_l$}}}
\put(2551,-886){\makebox(0,0)[lb]{\smash{\SetFigFont{10}{14.4}{rm}$\kappa_l$}}}
\end{picture}
\end{center}

\bigskip
\noindent in which the left vertical exists uniquely  by the 5-lemma
\cite{CWM}.

Coherence conditions follow from the uniqueness condition in the 
5-lemma by applying in the hexagon of exact sequences of Figure 
\ref{cohere.by.5}, and a similar triangle of exact sequences for
the unit condition.

\begin{figure}[htb] \centering
\setlength{\unitlength}{3290sp}%
\begingroup\makeatletter\ifx\SetFigFont\undefined%
\gdef\SetFigFont#1#2#3#4#5{%
  \reset@font\fontsize{#1}{#2pt}%
  \fontfamily{#3}\fontseries{#4}\fontshape{#5}%
  \selectfont}%
\fi\endgroup%
\begin{picture}(8700,4635)(301,-4036)
\put(301,-661){\makebox(0,0)[lb]{\smash{\SetFigFont{10}{12.0}{rm}$0$}}}
\put(1201,-661){\makebox(0,0)[lb]{\smash{\SetFigFont{10}{12.0}{rm}$[F(A) \otimes F(B)] \otimes ker(f_X)$}}}
\put(3901,-661){\makebox(0,0)[lb]{\smash{\SetFigFont{10}{12.0}{rm}$[F(A) \otimes F(B)] \otimes M(X)$}}}
\put(6601,-661){\makebox(0,0)[lb]{\smash{\SetFigFont{10}{12.0}{rm}$[F(A) \otimes F(B)] \otimes N(X)$}}}
\put(1201,314){\makebox(0,0)[lb]{\smash{\SetFigFont{10}{12.0}{rm}$0$}}}
\put(7501,314){\makebox(0,0)[lb]{\smash{\SetFigFont{10}{12.0}{rm}$F(A) \otimes [F(B) \otimes N(X)]$}}}
\put(2101,314){\makebox(0,0)[lb]{\smash{\SetFigFont{10}{12.0}{rm}$F(A) \otimes [F(B) \otimes ker(f_X)]$}}}
\put(4801,314){\makebox(0,0)[lb]{\smash{\SetFigFont{10}{12.0}{rm}$F(A) \otimes [F(B) \otimes M(X)]$}}}
\put(1201,-1561){\makebox(0,0)[lb]{\smash{\SetFigFont{10}{12.0}{rm}$0$}}}
\put(2401,-1561){\makebox(0,0)[lb]{\smash{\SetFigFont{10}{12.0}{rm}$F(A) \otimes ker(f_{B \triangleright X})$}}}
\put(5101,-1561){\makebox(0,0)[lb]{\smash{\SetFigFont{10}{12.0}{rm}$F(A) \otimes M(B \triangleright X)$}}}
\put(7801,-1561){\makebox(0,0)[lb]{\smash{\SetFigFont{10}{12.0}{rm}$F(A) \otimes N(B \triangleright X)$}}}
\put(301,-2161){\makebox(0,0)[lb]{\smash{\SetFigFont{10}{12.0}{rm}$0$}}}
\put(1501,-2161){\makebox(0,0)[lb]{\smash{\SetFigFont{10}{12.0}{rm}$F(A \otimes B) \otimes ker(f_X)$}}}
\put(4201,-2161){\makebox(0,0)[lb]{\smash{\SetFigFont{10}{12.0}{rm}$F(A \otimes B) \otimes M(X)$}}}
\put(6901,-2161){\makebox(0,0)[lb]{\smash{\SetFigFont{10}{12.0}{rm}$F(A \otimes B) \otimes N(X)$}}}
\put(1201,-2986){\makebox(0,0)[lb]{\smash{\SetFigFont{10}{12.0}{rm}$0$}}}
\put(2401,-2986){\makebox(0,0)[lb]{\smash{\SetFigFont{10}{12.0}{rm}$ker(f_{A \triangleright [B \triangleright X]})$}}}
\put(5101,-2986){\makebox(0,0)[lb]{\smash{\SetFigFont{10}{12.0}{rm}$M(A \triangleright [B \triangleright X])$}}}
\put(7801,-2986){\makebox(0,0)[lb]{\smash{\SetFigFont{10}{12.0}{rm}$N(A \triangleright [B \triangleright X])$}}}
\put(301,-3811){\makebox(0,0)[lb]{\smash{\SetFigFont{10}{12.0}{rm}$0$}}}
\put(1501,-3811){\makebox(0,0)[lb]{\smash{\SetFigFont{10}{12.0}{rm}$ker(f_{[A \otimes B] \triangleright X})$}}}
\put(4201,-3811){\makebox(0,0)[lb]{\smash{\SetFigFont{10}{12.0}{rm}$M([A \otimes B] \triangleright X)$}}}
\put(6901,-3811){\makebox(0,0)[lb]{\smash{\SetFigFont{10}{12.0}{rm}$N([A \otimes B] \triangleright X)$}}}
\thinlines
\put(3301,164){\vector( 0,-1){1425}}
\put(1951,-811){\vector( 0,-1){1050}}
\put(3301,-1711){\vector( 0,-1){975}}
\put(1951,-2386){\vector( 0,-1){1050}}
\put(2326,-3511){\vector( 3, 2){675}}
\put(2251,-286){\vector( 3, 2){675}}
\put(5926,164){\vector( 0,-1){1425}}
\put(4576,-811){\vector( 0,-1){1050}}
\put(5926,-1711){\vector( 0,-1){975}}
\put(4576,-2386){\vector( 0,-1){1050}}
\put(4951,-3511){\vector( 3, 2){675}}
\put(4876,-286){\vector( 3, 2){675}}
\put(8626,164){\vector( 0,-1){1425}}
\put(7276,-811){\vector( 0,-1){1050}}
\put(8626,-1711){\vector( 0,-1){975}}
\put(7276,-2386){\vector( 0,-1){1050}}
\put(7651,-3511){\vector( 3, 2){675}}
\put(7576,-286){\vector( 3, 2){675}}
\put(1401,389){\vector( 1, 0){550}}
\put(4176,389){\vector( 1, 0){475}}
\put(6801,389){\vector( 1, 0){550}}
\put(501,-586){\vector( 1, 0){550}}
\put(726,-2086){\vector( 1, 0){550}}
\put(1626,-2911){\vector( 1, 0){550}}
\put(651,-3736){\vector( 1, 0){550}}
\put(3276,-586){\vector( 1, 0){475}}
\put(3351,-3736){\vector( 1, 0){475}}
\put(5901,-586){\vector( 1, 0){550}}
\put(5976,-3736){\vector( 1, 0){550}}
\put(1626,-1486){\vector( 1, 0){550}}
\put(4326,-1486){\vector( 1, 0){475}}
\put(4326,-2911){\vector( 1, 0){475}}
\put(6951,-2911){\vector( 1, 0){550}}
\put(6951,-1486){\vector( 1, 0){550}}
\put(3276,-2086){\vector( 1, 0){475}}
\put(5901,-2086){\vector( 1, 0){550}}
\put(8851,-586){\makebox(0,0)[lb]{\smash{\SetFigFont{9}{10.8}{rm}$Id \otimes \nu_l$}}}
\put(8851,-2311){\makebox(0,0)[lb]{\smash{\SetFigFont{9}{10.8}{rm}$\nu_l$}}}
\put(8176,-3436){\makebox(0,0)[lb]{\smash{\SetFigFont{9}{10.8}{rm}$N(a_l)$}}}
\put(7351,-2686){\makebox(0,0)[lb]{\smash{\SetFigFont{9}{10.8}{rm}$\nu_l$}}}
\put(7351,-1186){\makebox(0,0)[lb]{\smash{\SetFigFont{9}{10.8}{rm}$\tilde{F} \otimes Id$}}}
\put(7651,-61){\makebox(0,0)[lb]{\smash{\SetFigFont{9}{10.8}{rm}$\alpha$}}}
\put(6151,-4036){\makebox(0,0)[lb]{\smash{\SetFigFont{9}{10.8}{rm}$f$}}}
\put(7051,-3136){\makebox(0,0)[lb]{\smash{\SetFigFont{9}{10.8}{rm}$f$}}}
\put(6001,-2311){\makebox(0,0)[lb]{\smash{\SetFigFont{9}{10.8}{rm}$Id \otimes f$}}}
\put(6751,-1711){\makebox(0,0)[lb]{\smash{\SetFigFont{9}{10.8}{rm}$Id \otimes f$}}}
\put(6751,464){\makebox(0,0)[lb]{\smash{\SetFigFont{9}{10.8}{rm}$Id \otimes [Id \otimes f]$}}}
\put(6001,-511){\makebox(0,0)[lb]{\smash{\SetFigFont{9}{10.8}{rm}$Id \otimes f$}}}
\put(5326,-3511){\makebox(0,0)[lb]{\smash{\SetFigFont{9}{10.8}{rm}$M(a_l)$}}}
\put(6001,-1936){\makebox(0,0)[lb]{\smash{\SetFigFont{9}{10.8}{rm}$\mu_l$}}}
\put(6001,-136){\makebox(0,0)[lb]{\smash{\SetFigFont{9}{10.8}{rm}$Id \otimes \mu_l$}}}
\put(4876,-61){\makebox(0,0)[lb]{\smash{\SetFigFont{9}{10.8}{rm}$\alpha$}}}
\put(4651,-1186){\makebox(0,0)[lb]{\smash{\SetFigFont{9}{10.8}{rm}$F \otimes Id$}}}
\put(4651,-2686){\makebox(0,0)[lb]{\smash{\SetFigFont{9}{10.8}{rm}$\mu_l$}}}
\put(2101,-61){\makebox(0,0)[lb]{\smash{\SetFigFont{9}{10.8}{rm}$\alpha$}}}
\put(2026,-1186){\makebox(0,0)[lb]{\smash{\SetFigFont{9}{10.8}{rm}$F \otimes Id$}}}
\put(3451,-286){\makebox(0,0)[lb]{\smash{\SetFigFont{9}{10.8}{rm}$Id \otimes \kappa_l$}}}
\put(3376,-1936){\makebox(0,0)[lb]{\smash{\SetFigFont{9}{10.8}{rm}$\kappa_l$}}}
\put(2026,-2686){\makebox(0,0)[lb]{\smash{\SetFigFont{9}{10.8}{rm}$\kappa_l$}}}
\put(2776,-3511){\makebox(0,0)[lb]{\smash{\SetFigFont{9}{10.8}{rm}$ker(f)(a_l)$}}}
\end{picture}

\caption{Coherence for the left action on a kernel \label{cohere.by.5}}
\end{figure}

Thus, the kernel has a unique left $F$-module structure such that
the inclusion is an $F$-module map.

Applying the same type of argument to the right action and its
hexagon and triangles, and to the two additional hexagons for a bimodule
shows that the kernel of a right $G$-module map (resp. $F,G$-bimodule map)
admits a unique $G$-module structure (resp. $F,G$-bimodule structure)
such that the inclusion is a $G$-module map (resp. $F,G,$-bimodule map).

Now, given an $F$-module map $e:L\Rightarrow M$ such that 
$f(e) = 0$, we claim that the induced map in ${\cal E}^{\cal M}$ is 
an $F$-module map. 

Now, for each $X \in Ob({\cal M})$ we have a diagram
\bigskip

\begin{center}
\setlength{\unitlength}{3947sp}%
\begingroup\makeatletter\ifx\SetFigFont\undefined%
\gdef\SetFigFont#1#2#3#4#5{%
  \reset@font\fontsize{#1}{#2pt}%
  \fontfamily{#3}\fontseries{#4}\fontshape{#5}%
  \selectfont}%
\fi\endgroup%
\begin{picture}(3600,1710)(1201,-1561)
\thinlines
\put(1426,-286){\vector( 1, 0){225}}
\put(2551,-286){\vector( 1, 0){600}}
\put(3826,-286){\vector( 1, 0){825}}
\put(2251,-1336){\vector( 4, 3){1200}}
\put(2026,-511){\vector( 0, 1){0}}
\multiput(2026,-1261)(0.00000,9.03614){83}{\makebox(1.6667,11.6667){\SetFigFont{5}{6}{rm}.}}
\put(1201,-361){\makebox(0,0)[lb]{\smash{\SetFigFont{10}{14.4}{rm}$0$}}}
\put(1801,-361){\makebox(0,0)[lb]{\smash{\SetFigFont{10}{14.4}{rm}$ker(f_X)$}}}
\put(3301,-361){\makebox(0,0)[lb]{\smash{\SetFigFont{10}{14.4}{rm}$M(X)$}}}
\put(4801,-361){\makebox(0,0)[lb]{\smash{\SetFigFont{10}{14.4}{rm}$N(X)$}}}
\put(1801,-1561){\makebox(0,0)[lb]{\smash{\SetFigFont{10}{14.4}{rm}$L(X)$}}}
\put(4051,-61){\makebox(0,0)[lb]{\smash{\SetFigFont{10}{14.4}{rm}$f_X$}}}
\put(2626, 14){\makebox(0,0)[lb]{\smash{\SetFigFont{10}{14.4}{rm}$i_X$}}}
\put(2851,-1261){\makebox(0,0)[lb]{\smash{\SetFigFont{10}{14.4}{rm}$e_X$}}}
\put(1501,-961){\makebox(0,0)[lb]{\smash{\SetFigFont{10}{14.4}{rm}$can_X$}}}
\end{picture}
\end{center}

\bigskip
\noindent with an exact top row.

As before, we form a diagram relating the image of this diagram
under $F(A)\otimes -$ and its instance for $A\triangleright X$ as in
Figure \ref{canonical.to.kernel}.  Here the ``back'' rows are both
exact, and the diagram commutes by construction except possibly from 
$F(A)\otimes L(X)$ to $ker(f_{A\triangleright X})$, the commutativity of which
is precisely what is to be shown.

\begin{figure}[htb] \centering
\setlength{\unitlength}{3947sp}%
\begingroup\makeatletter\ifx\SetFigFont\undefined
\def\x#1#2#3#4#5#6#7\relax{\def\x{#1#2#3#4#5#6}}%
\expandafter\x\fmtname xxxxxx\relax \def\y{splain}%
\ifx\x\y   
\gdef\SetFigFont#1#2#3{%
  \ifnum #1<17\tiny\else \ifnum #1<20\small\else
  \ifnum #1<24\normalsize\else \ifnum #1<29\large\else
  \ifnum #1<34\Large\else \ifnum #1<41\LARGE\else
     \huge\fi\fi\fi\fi\fi\fi
  \csname #3\endcsname}%
\else
\gdef\SetFigFont#1#2#3{\begingroup
  \count@#1\relax \ifnum 25<\count@\count@25\fi
  \def\x{\endgroup\@setsize\SetFigFont{#2pt}}%
  \expandafter\x
    \csname \romannumeral\the\count@ pt\expandafter\endcsname
    \csname @\romannumeral\the\count@ pt\endcsname
  \csname #3\endcsname}%
\fi
\fi\endgroup
\begin{picture}(5550,2685)(526,-2086)
\thinlines
\put(751,-961){\vector( 1, 0){675}}
\put(2551,-961){\vector( 1, 0){1200}}
\put(4876,-961){\vector( 1, 0){600}}
\put(5926,164){\vector( 0,-1){900}}
\put(4351, 89){\vector( 0,-1){825}}
\put(3226,-586){\vector( 0,-1){1200}}
\put(2701,314){\vector( 1, 0){1050}}
\put(1951,164){\vector( 0,-1){900}}
\put(2776,-211){\vector(-1, 1){375}}
\put(2851,-1861){\vector(-1, 1){600}}
\put(3451,-211){\vector( 3, 2){450}}
\put(3601,-1786){\vector( 1, 1){600}}
\put(5101,314){\vector( 1, 0){225}}
\put(901,314){\vector( 1, 0){450}}
\put(1501,239){\makebox(0,0)[lb]{\smash{\SetFigFont{10}{14.4}{rm}$F(A) \otimes ker(f_X)$}}}
\put(3901,239){\makebox(0,0)[lb]{\smash{\SetFigFont{10}{14.4}{rm}$F(A) \otimes M(X)$}}}
\put(5401,239){\makebox(0,0)[lb]{\smash{\SetFigFont{10}{14.4}{rm}$F(A) \otimes N(X)$}}}
\put(526,-1036){\makebox(0,0)[lb]{\smash{\SetFigFont{10}{14.4}{rm}$0$}}}
\put(1576,-1036){\makebox(0,0)[lb]{\smash{\SetFigFont{10}{14.4}{rm}$ker(f_{A \triangleright X}$}}}
\put(3976,-1036){\makebox(0,0)[lb]{\smash{\SetFigFont{10}{14.4}{rm}$M(A \triangleright X)$}}}
\put(5626,-1036){\makebox(0,0)[lb]{\smash{\SetFigFont{10}{14.4}{rm}$N(A \triangleright X)$}}}
\put(2776,-511){\makebox(0,0)[lb]{\smash{\SetFigFont{10}{14.4}{rm}$F(A) \otimes L(X)$}}}
\put(2926,-2086){\makebox(0,0)[lb]{\smash{\SetFigFont{10}{14.4}{rm}$L(A \triangleright X)$}}}
\put(526,239){\makebox(0,0)[lb]{\smash{\SetFigFont{10}{14.4}{rm}$0$}}}
\put(2776,464){\makebox(0,0)[lb]{\smash{\SetFigFont{10}{14.4}{rm}$Id \otimes i_X$}}}
\put(4951,464){\makebox(0,0)[lb]{\smash{\SetFigFont{10}{14.4}{rm}$Id \otimes f_X$}}}
\put(4876,-1261){\makebox(0,0)[lb]{\smash{\SetFigFont{10}{14.4}{rm}$f_{A \triangleright X}$}}}
\put(6076,-211){\makebox(0,0)[lb]{\smash{\SetFigFont{10}{14.4}{rm}$\nu_l$}}}
\put(4501,-286){\makebox(0,0)[lb]{\smash{\SetFigFont{10}{14.4}{rm}$\mu_l$}}}
\put(3301,-1336){\makebox(0,0)[lb]{\smash{\SetFigFont{10}{14.4}{rm}$\lambda_l$}}}
\put(1276,-286){\makebox(0,0)[lb]{\smash{\SetFigFont{10}{14.4}{rm}$\kappa_l$}}}
\put(1726,-1786){\makebox(0,0)[lb]{\smash{\SetFigFont{10}{14.4}{rm}$can{A \triangleright X}$}}}
\put(2701, 14){\makebox(0,0)[lb]{\smash{\SetFigFont{10}{14.4}{rm}$Id \otimes can_X$}}}
\put(3676,-286){\makebox(0,0)[lb]{\smash{\SetFigFont{10}{14.4}{rm}$Id \otimes eX$}}}
\put(3901,-1786){\makebox(0,0)[lb]{\smash{\SetFigFont{10}{14.4}{rm}$e_{A \triangleright X}$}}}
\end{picture}

\caption{Canonical maps are module maps \label{canonical.to.kernel}}
\end{figure}

Now, $f_{A\triangleright X}(e_{A\triangleright X}(\lambda_l)) = 0$ (since
$f_{A\triangleright X}(e_{A\triangleright X}) = 0$).  Thus by the universal 
property of $ker(f_{A\triangleright X})$, there exists a unique map 
$\phi:F(A)\otimes L(X)\rightarrow ker(f_{A\triangleright X})$ 
(in $\cal D$) such that $\iota_{A\triangleright X}(\phi) = e_{A\triangleright X}(\lambda_l)$.

On the one hand, it is clear that $\phi = can_{A\triangleright X}(\lambda_l)$
since $\iota_{A\triangleright X}(can_{A\triangleright X}) = e_{A\triangleright X}$.

On the other, we have

\begin{eqnarray*}
\iota_{A\triangleright X}(\kappa_l(Id\otimes can_X)) & = &
		\mu_l(Id\otimes \iota_X(Id\otimes can_X)) \\
	& = & \mu_l(Id\otimes e_X) \\
	& = & e_{A\triangleright X}(\lambda_l)
\end{eqnarray*}

\noindent where the first equality holds by the construction of
the kernel action $\kappa_l$, the second by the definition of $can_X$,
and the third by the fact that $e$ is an $F$-module map.
Thus $\phi = \kappa_l(Id\otimes can_X)$ and the required square
commutes by the uniqueness of $\phi$.

An essentially identical proof gives the corresponding result for
right $G$-modules, and thus for $F,G$-bimodules.

The proof that cokernels of (bi)module maps have (bi)module structures
and that the quotient map and canonical maps are (bi)module maps may
be obtained by dualizing the rows (only) in the proof for kernels.

As regards biproducts, observe that since the universal properties 
follow from the equational conditions on the projections and inclusions,
it suffices to show that 
given two left $F$-modules $M$ and $N$, the biproduct of their 
underlying objects in ${\cal E}^{\cal M}$ admit 
a unique left $F$-module structure such that the inclusions and projections
are left module maps.  

Now, since biproducts are equationally defined $F(A) \otimes -$ preserves
biproducts up to a canonical isomorphism commuting with both projections
and inclusions.  Thus, if we let 

\[ d: F(A) \otimes [M(X) \oplus N(X)] \rightarrow [F(A) \otimes M(X)] \oplus
[F(A) \otimes N(X)] \]

\noindent  denote the canonical arrow, the diagram obtained
from that of Figure \ref{left.mod.biprod} by omitting either projections
or inclusions commutes.  Thus $[\mu_l \oplus \nu_l](d)$ give the unique
left module structure on $M\oplus N$.

Coherence follows from the universal properties of biproducts.

\begin{figure}[htb] \centering
\setlength{\unitlength}{3947sp}%
\begingroup\makeatletter\ifx\SetFigFont\undefined
\def\x#1#2#3#4#5#6#7\relax{\def\x{#1#2#3#4#5#6}}%
\expandafter\x\fmtname xxxxxx\relax \def\y{splain}%
\ifx\x\y   
\gdef\SetFigFont#1#2#3{%
  \ifnum #1<17\tiny\else \ifnum #1<20\small\else
  \ifnum #1<24\normalsize\else \ifnum #1<29\large\else
  \ifnum #1<34\Large\else \ifnum #1<41\LARGE\else
     \huge\fi\fi\fi\fi\fi\fi
  \csname #3\endcsname}%
\else
\gdef\SetFigFont#1#2#3{\begingroup
  \count@#1\relax \ifnum 25<\count@\count@25\fi
  \def\x{\endgroup\@setsize\SetFigFont{#2pt}}%
  \expandafter\x
    \csname \romannumeral\the\count@ pt\expandafter\endcsname
    \csname @\romannumeral\the\count@ pt\endcsname
  \csname #3\endcsname}%
\fi
\fi\endgroup
\begin{picture}(6462,2664)(601,-2140)
\thinlines
\put(1051,-511){\vector( 0,-1){900}}
\put(7051,-511){\vector( 0,-1){900}}
\put(4051,-586){\vector( 0,-1){825}}
\put(1576,-1561){\vector( 1, 0){1425}}
\put(3001,-1711){\vector(-1, 0){1425}}
\put(6526,-1486){\vector(-1, 0){1575}}
\put(4951,-1711){\vector( 1, 0){1500}}
\put(1726,-211){\vector( 1, 0){975}}
\put(2701,-361){\vector(-1, 0){975}}
\put(6451,-211){\vector(-1, 0){900}}
\put(5626,-361){\vector( 1, 0){825}}
\put(4051,164){\vector( 0,-1){300}}
\put(1276,-61){\vector( 4, 1){1800}}
\put(3076,314){\vector(-4,-1){1500}}
\put(6976,-61){\vector(-4, 1){1800}}
\put(5251,314){\vector( 4,-1){1800}}
\put(601,-361){\makebox(0,0)[lb]{\smash{\SetFigFont{10}{14.4}{rm}$F(A) \otimes M(X)$}}}
\put(6601,-361){\makebox(0,0)[lb]{\smash{\SetFigFont{10}{14.4}{rm}$F(A) \otimes N(X)$}}}
\put(676,-1711){\makebox(0,0)[lb]{\smash{\SetFigFont{10}{14.4}{rm}$M(A \triangleright X)$}}}
\put(6676,-1636){\makebox(0,0)[lb]{\smash{\SetFigFont{10}{14.4}{rm}$N(A \triangleright X)$}}}
\put(3226,314){\makebox(0,0)[lb]{\smash{\SetFigFont{10}{14.4}{rm}$F(A) \otimes [M(X) \oplus N(X)]$}}}
\put(2851,-361){\makebox(0,0)[lb]{\smash{\SetFigFont{10}{14.4}{rm}$[F(A) \otimes M(X)] \oplus [F(A) \otimes N(X)]$}}}
\put(3151,-1711){\makebox(0,0)[lb]{\smash{\SetFigFont{10}{14.4}{rm}$M(A \triangleright X) \oplus N(A \triangleright X)$}}}
\put(2101,-1411){\makebox(0,0)[lb]{\smash{\SetFigFont{10}{14.4}{rm}$i_1$}}}
\put(5776,-1411){\makebox(0,0)[lb]{\smash{\SetFigFont{10}{14.4}{rm}$i_2$}}}
\put(2101,-2011){\makebox(0,0)[lb]{\smash{\SetFigFont{10}{14.4}{rm}$p_1$}}}
\put(5476,-2086){\makebox(0,0)[lb]{\smash{\SetFigFont{10}{14.4}{rm}$p_2$}}}
\put(5776,-136){\makebox(0,0)[lb]{\smash{\SetFigFont{10}{14.4}{rm}$i_2$}}}
\put(5851,-586){\makebox(0,0)[lb]{\smash{\SetFigFont{10}{14.4}{rm}$p_2$}}}
\put(2251,-136){\makebox(0,0)[lb]{\smash{\SetFigFont{10}{14.4}{rm}$i_1$}}}
\put(2326,-586){\makebox(0,0)[lb]{\smash{\SetFigFont{10}{14.4}{rm}$p_2$}}}
\put(1726,314){\makebox(0,0)[lb]{\smash{\SetFigFont{10}{14.4}{rm}$Id \otimes i_1$}}}
\put(2626, 14){\makebox(0,0)[lb]{\smash{\SetFigFont{10}{14.4}{rm}$Id \otimes p_1$}}}
\put(5776,389){\makebox(0,0)[lb]{\smash{\SetFigFont{10}{14.4}{rm}$Id \otimes i_2$}}}
\put(4801, 89){\makebox(0,0)[lb]{\smash{\SetFigFont{10}{14.4}{rm}$Id \otimes p_2$}}}
\put(4201,-61){\makebox(0,0)[lb]{\smash{\SetFigFont{10}{14.4}{rm}$d$}}}
\put(1201,-961){\makebox(0,0)[lb]{\smash{\SetFigFont{10}{14.4}{rm}$\mu_l$}}}
\put(7201,-961){\makebox(0,0)[lb]{\smash{\SetFigFont{10}{14.4}{rm}$\nu_l$}}}
\put(4201,-961){\makebox(0,0)[lb]{\smash{\SetFigFont{10}{14.4}{rm}$\mu_l \oplus \nu_l$}}}
\end{picture}

\caption{Left module structure for biproducts \label{left.mod.biprod}}
\end{figure}

The construction of the right module structure for biproducts of right
modules (or bimodules) is entirely similar.

Finally, we must show that the parallel of a map (cf. Popescu \cite{Pop}) 
is invertible in the category of (bi)modules.  It already is invertible
in ${\cal E}^{\cal M}$, so it suffices to show that a left $F$-module map
which is invertible in ${\cal E}^{\cal M}$ is invertible in the category
of left $F$-modules over $\cal M$.  (Right actions will follow by an
essentially identical proof.)

Now, observe that if $f:M \Rightarrow N$ is invertible in ${\cal E}^{\cal M}$,
we must show that the diagram

\begin{center}
\setlength{\unitlength}{3947sp}%
\begingroup\makeatletter\ifx\SetFigFont\undefined%
\gdef\SetFigFont#1#2#3#4#5{%
  \reset@font\fontsize{#1}{#2pt}%
  \fontfamily{#3}\fontseries{#4}\fontshape{#5}%
  \selectfont}%
\fi\endgroup%
\begin{picture}(3300,2010)(76,-1711)
\thinlines
\put(1726,-1486){\vector( 1, 0){1050}}
\put(3226,-211){\vector( 0,-1){1050}}
\put(1726, 14){\vector( 1, 0){1050}}
\put(1201,-211){\vector( 0,-1){1050}}
\put(2101,164){\makebox(0,0)[lb]{\smash{\SetFigFont{10}{14.4}{rm}$\mu_l$}}}
\put(2026,-1711){\makebox(0,0)[lb]{\smash{\SetFigFont{10}{14.4}{rm}$\nu_l$}}}
\put(601,-61){\makebox(0,0)[lb]{\smash{\SetFigFont{10}{14.4}{rm}$F(A) \otimes N(X)$}}}
\put(2926,-61){\makebox(0,0)[lb]{\smash{\SetFigFont{10}{14.4}{rm}$N(A \triangleright X)$}}}
\put(601,-1561){\makebox(0,0)[lb]{\smash{\SetFigFont{10}{14.4}{rm}$F(A) \otimes M(X)$}}}
\put(2926,-1561){\makebox(0,0)[lb]{\smash{\SetFigFont{10}{14.4}{rm}$M(A \triangleright X)$}}}
\put( 76,-736){\makebox(0,0)[lb]{\smash{\SetFigFont{10}{14.4}{rm}$F(A) \otimes f_X^{-1}$}}}
\put(3376,-736){\makebox(0,0)[lb]{\smash{\SetFigFont{10}{14.4}{rm}$f_{A \triangleright X}^{-1}$}}}
\end{picture}
\end{center}

\noindent commutes.  However, since both $F(A) \otimes f_X$ and $f_{A \triangleright X}$
are invertible in $\cal E$, it suffices to see that the diagram

\begin{center}
\setlength{\unitlength}{3947sp}%
\begingroup\makeatletter\ifx\SetFigFont\undefined%
\gdef\SetFigFont#1#2#3#4#5{%
  \reset@font\fontsize{#1}{#2pt}%
  \fontfamily{#3}\fontseries{#4}\fontshape{#5}%
  \selectfont}%
\fi\endgroup%
\begin{picture}(4437,3480)(1201,-3286)
\thinlines
\put(3376,-2236){\vector( 1, 0){1050}}
\put(4876,-961){\vector( 0,-1){1050}}
\put(3376,-736){\vector( 1, 0){1050}}
\put(2851,-961){\vector( 0,-1){1050}}
\put(3751,-586){\makebox(0,0)[lb]{\smash{\SetFigFont{10}{14.4}{rm}$\mu_l$}}}
\put(3676,-2461){\makebox(0,0)[lb]{\smash{\SetFigFont{10}{14.4}{rm}$\nu_l$}}}
\put(2251,-811){\makebox(0,0)[lb]{\smash{\SetFigFont{10}{14.4}{rm}$F(A) \otimes N(X)$}}}
\put(4576,-811){\makebox(0,0)[lb]{\smash{\SetFigFont{10}{14.4}{rm}$N(A \triangleright X)$}}}
\put(2251,-2311){\makebox(0,0)[lb]{\smash{\SetFigFont{10}{14.4}{rm}$F(A) \otimes M(X)$}}}
\put(4576,-2311){\makebox(0,0)[lb]{\smash{\SetFigFont{10}{14.4}{rm}$M(A \triangleright X)$}}}
\put(1726,-1486){\makebox(0,0)[lb]{\smash{\SetFigFont{10}{14.4}{rm}$F(A) \otimes f_X^{-1}$}}}
\put(5026,-1486){\makebox(0,0)[lb]{\smash{\SetFigFont{10}{14.4}{rm}$f{A \triangleright X}^{-1}$}}}
\put(2101,-136){\vector( 1,-1){300}}
\put(5176,-2461){\vector( 1,-1){450}}
\put(1201, 14){\makebox(0,0)[lb]{\smash{\SetFigFont{10}{14.4}{rm}$F(A) \otimes M(X)$}}}
\put(5476,-3286){\makebox(0,0)[lb]{\smash{\SetFigFont{10}{14.4}{rm}$N(A \triangleright X)$}}}
\put(2401,-286){\makebox(0,0)[lb]{\smash{\SetFigFont{10}{14.4}{rm}$F(A) \otimes f_X$}}}
\put(5626,-2686){\makebox(0,0)[lb]{\smash{\SetFigFont{10}{14.4}{rm}$f_{A \triangleright X}$}}}
\end{picture}
\end{center}

\noindent commutes from $F(A) \otimes M(X)$ to $N(A \triangleright X)$.  But calculating
the two paths, we see that this is just the coherence diagram for $f$ as a
left $F$-module map. 

$\Box$

A similar technique will show

\begin{thm} \label{Ab3}
If $\cal E$ is a cocomplete abelian category equipped with a monoidal 
structure,
and \hspace{1cm}
$F:{\cal C}\rightarrow {\cal E}$ and $G:{\cal D}\rightarrow {\cal E}$
are monoidal functors such that for all $A \in Ob({\cal C})$
$F(A)\otimes -$ is exact and cocontinuous (resp. for all $C \in Ob({\cal D})$
$-\otimes G(C)$ is exact and cocontinuous, both), 
then for any left $\cal C$-module
(resp. right $\cal D$-module, $\cal C$,$\cal D$-bimodule) $\cal M$,
the category of left $F$-modules (resp. right $G$-modules, $F,G$-bimodules)
over $\cal M$ is a cocomplete abelian category.
\end{thm}

Although the hypotheses of Theorems \ref{abelian} and \ref{Ab3} may seem rather
restrictive, it should be noted that they hold whenever the target category
$\cal E$ is of the form $k${\bf -v.s.} for $k$ a field.

With the category of $F,F$-bimodules as a setting, we can generalize the 
deformation complex of $F$ to give an analogue of the Hochschild cohomology of
an algebra with coefficients in a bimodule:

Let

\[ X^n(F,M) = Nat( ^n\otimes(F^n), M(\otimes^n)) \]

\noindent with coboundary given by the the obvious generalization
of the formula for the deformation complex.

Deformation complexes for semigroupal categories and functors are then
the special case for $F = M = Id_{\cal C}$ and $M = F$ respectively.  
There is also another special case of interest.  If $M$ is a monoidal functor
$G$ made into an $F,F$-bimodule with the structure induced by a 
monoidal natural transformation $\phi:F\Rightarrow G$ as described above, it
is easy to show that the first order deformations of $\phi$ as a 
monoidal natural transformation are classified by $H^1(F,G)$, and the
obstructions to higher order deformations are classes in $H^2(F,G)$.

\clearpage  
\section{$Lan_{\otimes^n}( ^n\otimes (F^n))$}
\vspace*{1cm}

In this section, we turn to the principal use of Theorems \ref{abelian} and
\ref{Ab3}:  the reduction of the cohomology of a (cocomplete) abelian 
monoidal category (with exact cocontinuous $\otimes$), and of the cohomology 
of a monoidal
functor targetted at such a category, with coefficients in a module
to a calculus of derived functors.  We will need some auxiliary hypotheses,
these, however, will be
satisfied both by algebras viewed as lax monoidal functors, and by strong
monoidal functors targetted at categories of vector-spaces, provided
the source category satisfies the mild size restriction that there exists
a small subcategory such that all objects are colimits of diagrams in the
small subcategory.

It should be observed that even when it is possible, this reduction does 
not in itself solve the problems
of categorical deformations:  the resulting expression for the cohomology
groups does not shed very much light on the behavior of the obstruction
cocyles, which are governed by the pre-Lie structure.  Nonetheless, the
result is important, in that it places categorical deformation theory
comfortably within the realm of classical homological algebra.

Recall that the $n^{th}$ cochain group associated to a lax monoidal functor 
$F:{\cal C}\rightarrow {\cal E}$
with coefficients in the $F$-module $M:{\cal C}\rightarrow {\cal E}$
is defined by

\[ X^n(F,M) = Nat( ^n\otimes(F^n), M(\otimes^n)) .\]

Note that as $n$ varies, the source category for the functors varies.  This 
rather uncomfortable circumstance from the point of view of classical 
homological algebra can be rectified if the functors $ ^n\otimes(F^n)$
admit left Kan extensions along $\otimes^n$ for all $n$.  In this case 
we can redefine $X^n(F,M)$ by

\[ X^n(F,M) = Nat(Lan_{\otimes^n}( ^n\otimes (F^n)), M). \]

Here, regardless of $n$ we have the abelian group of natural transformations
between functors from $\cal C$ to $\cal E$.

Our first goal, then, is to show under suitable hypotheses that for $n \geq 2$
the functors 

\[ Lan_{\otimes^n}( ^n\otimes (F^n)):{\cal C}\rightarrow {\cal E} \]
 
\noindent are in fact $F,F$-bimodules over $\cal C$ (with 
$\triangleright = \triangleleft
= \otimes$)
in a natural way, and are projective as such.  Our second goal is to
show that the
cochain complex $(X^\bullet(F,M),\delta^\bullet)$ arises by applying
$Nat[-,M]$ to a projective resolution made up of these Kan extensions.

Before embarking on the construction of the desired lifts, let us show that the Kan extensions admit an $F,F$-bimodule structure, and examine which maps from iterated tensor products of $F$ correspond to bimodule maps from the corresponding Kan extension.

Assume that $F(X)\otimes -$ and $-\otimes F(X)$ are cocontinuous for all 
$X$.  Now, this being so, we have

\[ F(X)\otimes Lan_{\otimes^n}(^n\otimes(F^n)) = 
         Lan_{\otimes^n}(F(X)\otimes ^n\otimes(F^n)). \]

Thus a left action corresponds to a natural transformation

\[ Lan_{\otimes^n}(F(X)\otimes ^n\otimes(F^n))\Rightarrow
Lan_{\otimes^n}(^n\otimes(F^n))(X\otimes -) \,. \]

By the universal property of the target, this in turn corresponds to
a natural filling of the rectangle

\begin{center} 
\setlength{\unitlength}{3947sp}%
\begingroup\makeatletter\ifx\SetFigFont\undefined
\def\x#1#2#3#4#5#6#7\relax{\def\x{#1#2#3#4#5#6}}%
\expandafter\x\fmtname xxxxxx\relax \def\y{splain}%
\ifx\x\y   
\gdef\SetFigFont#1#2#3{%
  \ifnum #1<17\tiny\else \ifnum #1<20\small\else
  \ifnum #1<24\normalsize\else \ifnum #1<29\large\else
  \ifnum #1<34\Large\else \ifnum #1<41\LARGE\else
     \huge\fi\fi\fi\fi\fi\fi
  \csname #3\endcsname}%
\else
\gdef\SetFigFont#1#2#3{\begingroup
  \count@#1\relax \ifnum 25<\count@\count@25\fi
  \def\x{\endgroup\@setsize\SetFigFont{#2pt}}%
  \expandafter\x
    \csname \romannumeral\the\count@ pt\expandafter\endcsname
    \csname @\romannumeral\the\count@ pt\endcsname
  \csname #3\endcsname}%
\fi
\fi\endgroup
\begin{picture}(4215,2040)(331,-1831)
\thinlines
\put(676,-361){\vector( 0,-1){975}}
\put(976,-136){\vector( 1, 0){1275}}
\put(2926,-136){\vector( 1, 0){1200}}
\put(2776,-1486){\vector( 1, 0){1350}}
\put(4351,-286){\vector( 0,-1){975}}
\put(901,-1486){\vector( 1, 0){1425}}
\put(2326,-781){\line( 0,-1){195}}
\put(2341,-991){\line( 1, 0){195}}
\put(2341,-856){\line( 1, 1){165}}
\put(2506,-691){\line(-1, 0){ 15}}
\put(2416,-991){\line( 1, 1){195}}
\put(2611,-796){\line(-1, 0){ 15}}
\put(601,-211){\makebox(0,0)[lb]{\smash{\SetFigFont{10}{14.4}{rm}${\cal C}^{\boxtimes n}$}}}
\put(2401,-211){\makebox(0,0)[lb]{\smash{\SetFigFont{10}{14.4}{rm}${\cal D}^{\boxtimes n}$}}}
\put(4201,-211){\makebox(0,0)[lb]{\smash{\SetFigFont{10}{14.4}{rm}${\cal D}^{\boxtimes n}$}}}
\put(676,-1561){\makebox(0,0)[lb]{\smash{\SetFigFont{10}{14.4}{rm}$\cal C$}}}
\put(2476,-1561){\makebox(0,0)[lb]{\smash{\SetFigFont{10}{14.4}{rm}$\cal C$}}}
\put(4276,-1561){\makebox(0,0)[lb]{\smash{\SetFigFont{10}{14.4}{rm}$\cal D$}}}
\put(1426, 14){\makebox(0,0)[lb]{\smash{\SetFigFont{10}{14.4}{rm}$F^n$}}}
\put(2806, 29){\makebox(0,0)[lb]{\smash{\SetFigFont{10}{14.4}{rm}$F(X)\times - \boxtimes {\cal D}^{\boxtimes n-1}$}}}
\put(331,-766){\makebox(0,0)[lb]{\smash{\SetFigFont{10}{14.4}{rm}$^n\otimes$}}}
\put(4546,-796){\makebox(0,0)[lb]{\smash{\SetFigFont{10}{14.4}{rm}$\otimes^n$}}}
\put(1216,-1831){\makebox(0,0)[lb]{\smash{\SetFigFont{10}{14.4}{rm}$X\otimes -$}}}
\put(2761,-1831){\makebox(0,0)[lb]{\smash{\SetFigFont{10}{14.4}{rm}$Lan_{^n\otimes}(\otimes^n(F^n))$}}}
\end{picture}

\end{center}

But, we have

\begin{center} 
\setlength{\unitlength}{3947sp}%
\begingroup\makeatletter\ifx\SetFigFont\undefined
\def\x#1#2#3#4#5#6#7\relax{\def\x{#1#2#3#4#5#6}}%
\expandafter\x\fmtname xxxxxx\relax \def\y{splain}%
\ifx\x\y   
\gdef\SetFigFont#1#2#3{%
  \ifnum #1<17\tiny\else \ifnum #1<20\small\else
  \ifnum #1<24\normalsize\else \ifnum #1<29\large\else
  \ifnum #1<34\Large\else \ifnum #1<41\LARGE\else
     \huge\fi\fi\fi\fi\fi\fi
  \csname #3\endcsname}%
\else
\gdef\SetFigFont#1#2#3{\begingroup
  \count@#1\relax \ifnum 25<\count@\count@25\fi
  \def\x{\endgroup\@setsize\SetFigFont{#2pt}}%
  \expandafter\x
    \csname \romannumeral\the\count@ pt\expandafter\endcsname
    \csname @\romannumeral\the\count@ pt\endcsname
  \csname #3\endcsname}%
\fi
\fi\endgroup
\begin{picture}(4215,3000)(331,-2791)
\thinlines
\put(2746,-2446){\vector( 1, 0){1350}}
\put(871,-2446){\vector( 1, 0){1425}}
\put(646,-2521){\makebox(0,0)[lb]{\smash{\SetFigFont{10}{14.4}{rm}$\cal C$}}}
\put(2446,-2521){\makebox(0,0)[lb]{\smash{\SetFigFont{10}{14.4}{rm}$\cal C$}}}
\put(4246,-2521){\makebox(0,0)[lb]{\smash{\SetFigFont{10}{14.4}{rm}$\cal D$ .}}}
\put(1186,-2791){\makebox(0,0)[lb]{\smash{\SetFigFont{10}{14.4}{rm}$X\otimes -$}}}
\put(2731,-2791){\makebox(0,0)[lb]{\smash{\SetFigFont{10}{14.4}{rm}$Lan_{^n\otimes}(\otimes^n(F^n))$}}}
\multiput(1283,-1621)(2.76471,-8.29412){18}{\makebox(1.6667,11.6667){\SetFigFont{5}{6}{rm}.}}
\multiput(1330,-1762)(8.25000,2.06250){17}{\makebox(1.6667,11.6667){\SetFigFont{5}{6}{rm}.}}
\put(1313,-1703){\line( 1, 2){120}}
\put(1396,-1741){\line( 1, 2){120}}
\multiput(2175,-505)(2.76471,-8.29412){18}{\makebox(1.6667,11.6667){\SetFigFont{5}{6}{rm}.}}
\multiput(2222,-646)(8.25000,2.06250){17}{\makebox(1.6667,11.6667){\SetFigFont{5}{6}{rm}.}}
\put(2205,-587){\line( 1, 2){120}}
\put(2288,-625){\line( 1, 2){120}}
\multiput(3053,-1697)(2.76471,-8.29412){18}{\makebox(1.6667,11.6667){\SetFigFont{5}{6}{rm}.}}
\multiput(3100,-1838)(8.25000,2.06250){17}{\makebox(1.6667,11.6667){\SetFigFont{5}{6}{rm}.}}
\put(3083,-1779){\line( 1, 2){120}}
\put(3166,-1817){\line( 1, 2){120}}
\put(976,-136){\vector( 1, 0){1275}}
\put(2926,-136){\vector( 1, 0){1200}}
\put(691,-361){\vector( 0,-1){1890}}
\put(4381,-331){\vector( 0,-1){1935}}
\put(916,-316){\vector( 2,-1){1456}}
\put(2776,-1066){\vector( 2, 1){1440}}
\put(2506,-1381){\vector( 0,-1){840}}
\put(601,-211){\makebox(0,0)[lb]{\smash{\SetFigFont{10}{14.4}{rm}${\cal C}^{\boxtimes n}$}}}
\put(2401,-211){\makebox(0,0)[lb]{\smash{\SetFigFont{10}{14.4}{rm}${\cal D}^{\boxtimes n}$}}}
\put(4201,-211){\makebox(0,0)[lb]{\smash{\SetFigFont{10}{14.4}{rm}${\cal D}^{\boxtimes n}$}}}
\put(1426, 14){\makebox(0,0)[lb]{\smash{\SetFigFont{10}{14.4}{rm}$F^n$}}}
\put(2806, 29){\makebox(0,0)[lb]{\smash{\SetFigFont{10}{14.4}{rm}$F(X)\otimes - \boxtimes {\cal D}^{\boxtimes n-1}$}}}
\put(331,-766){\makebox(0,0)[lb]{\smash{\SetFigFont{10}{14.4}{rm}$^n\otimes$}}}
\put(4546,-796){\makebox(0,0)[lb]{\smash{\SetFigFont{10}{14.4}{rm}$\otimes^n$}}}
\put(2401,-1261){\makebox(0,0)[lb]{\smash{\SetFigFont{10}{14.4}{rm}${\cal C}^{\boxtimes n}$}}}
\put(1568,-1688){\makebox(0,0)[lb]{\smash{\SetFigFont{10}{14.4}{rm}\mbox{\boldmath $\alpha$}}}}
\put(3331,-1801){\makebox(0,0)[lb]{\smash{\SetFigFont{10}{14.4}{rm}$can$}}}
\put(2461,-556){\makebox(0,0)[lb]{\smash{\SetFigFont{10}{14.4}{rm}$\tilde{F}\boxtimes Id$}}}
\end{picture}

\end{center}

The coherence for the left action follows from the coherence of
$\tilde{F}$ and universality.  The construction of the right action
is completely analogous.  The additional coherence condition for
the $F,F$-bimodule structure follows trivially from the separation
of the actions.  (Note:  it fails for $n = 1$.)

We can now consider what condition on a natural tranformation
$\tilde{\phi}:^n\otimes(F^n)\Rightarrow M(\otimes^n)$ 
implies that the induced
natural transformation $\phi:Lan_{\otimes^n}(^n\otimes(F^n)
\Rightarrow M$ is a left (resp. right, bi-) module map.

We need all instances of 

\begin{center} 
\setlength{\unitlength}{3300sp}%
\begingroup\makeatletter\ifx\SetFigFont\undefined
\def\x#1#2#3#4#5#6#7\relax{\def\x{#1#2#3#4#5#6}}%
\expandafter\x\fmtname xxxxxx\relax \def\y{splain}%
\ifx\x\y   
\gdef\SetFigFont#1#2#3{%
  \ifnum #1<17\tiny\else \ifnum #1<20\small\else
  \ifnum #1<24\normalsize\else \ifnum #1<29\large\else
  \ifnum #1<34\Large\else \ifnum #1<41\LARGE\else
     \huge\fi\fi\fi\fi\fi\fi
  \csname #3\endcsname}%
\else
\gdef\SetFigFont#1#2#3{\begingroup
  \count@#1\relax \ifnum 25<\count@\count@25\fi
  \def\x{\endgroup\@setsize\SetFigFont{#2pt}}%
  \expandafter\x
    \csname \romannumeral\the\count@ pt\expandafter\endcsname
    \csname @\romannumeral\the\count@ pt\endcsname
  \csname #3\endcsname}%
\fi
\fi\endgroup
\begin{picture}(5025,1380)(226,-961)
\thinlines
\put(1426,164){\vector( 0,-1){825}}
\put(2776,314){\vector( 1, 0){1275}}
\put(2101,-886){\vector( 1, 0){2400}}
\put(5026, 89){\vector( 0,-1){750}}
\put(601,239){\makebox(0,0)[lb]{\smash{\SetFigFont{9}{14.4}{rm}$F(A) \otimes Lan_{^n\otimes}(\otimes^n(F^n))(B)$}}}
\put(4201,239){\makebox(0,0)[lb]{\smash{\SetFigFont{9}{14.4}{rm}$Lan_{^n\otimes}(\otimes^n(F^n))(A \otimes B)$}}}
\put(901,-961){\makebox(0,0)[lb]{\smash{\SetFigFont{9}{14.4}{rm}$F(A) \otimes M(B)$}}}
\put(226,-436){\makebox(0,0)[lb]{\smash{\SetFigFont{9}{14.4}{rm}$F(A) \otimes \phi_B$}}}
\put(5251,-286){\makebox(0,0)[lb]{\smash{\SetFigFont{9}{14.4}{rm}$\phi_{A\otimes B}$}}}
\put(4651,-961){\makebox(0,0)[lb]{\smash{\SetFigFont{9}{14.4}{rm}$M(A \otimes B)$ .}}}
\end{picture}

\end{center}

\noindent The top path round this diagram is induced by

\begin{center} 
\setlength{\unitlength}{3947sp}%
\begingroup\makeatletter\ifx\SetFigFont\undefined
\def\x#1#2#3#4#5#6#7\relax{\def\x{#1#2#3#4#5#6}}%
\expandafter\x\fmtname xxxxxx\relax \def\y{splain}%
\ifx\x\y   
\gdef\SetFigFont#1#2#3{%
  \ifnum #1<17\tiny\else \ifnum #1<20\small\else
  \ifnum #1<24\normalsize\else \ifnum #1<29\large\else
  \ifnum #1<34\Large\else \ifnum #1<41\LARGE\else
     \huge\fi\fi\fi\fi\fi\fi
  \csname #3\endcsname}%
\else
\gdef\SetFigFont#1#2#3{\begingroup
  \count@#1\relax \ifnum 25<\count@\count@25\fi
  \def\x{\endgroup\@setsize\SetFigFont{#2pt}}%
  \expandafter\x
    \csname \romannumeral\the\count@ pt\expandafter\endcsname
    \csname @\romannumeral\the\count@ pt\endcsname
  \csname #3\endcsname}%
\fi
\fi\endgroup
\begin{picture}(4410,2475)(271,-1876)
\thinlines
\multiput(1351,-811)(3.75000,-7.50000){21}{\makebox(1.6667,11.6667){\SetFigFont{5}{6}{rm}.}}
\multiput(1426,-961)(7.50000,3.75000){21}{\makebox(1.6667,11.6667){\SetFigFont{5}{6}{rm}.}}
\put(1391,-896){\line( 1, 4){ 75}}
\put(1496,-916){\line( 1, 4){ 75}}
\multiput(3396,-856)(3.75000,-7.50000){21}{\makebox(1.6667,11.6667){\SetFigFont{5}{6}{rm}.}}
\multiput(3471,-1006)(7.50000,3.75000){21}{\makebox(1.6667,11.6667){\SetFigFont{5}{6}{rm}.}}
\put(3436,-941){\line( 1, 4){ 75}}
\put(3541,-961){\line( 1, 4){ 75}}
\multiput(2336,-21)(3.75000,-7.50000){21}{\makebox(1.6667,11.6667){\SetFigFont{5}{6}{rm}.}}
\multiput(2411,-171)(7.50000,3.75000){21}{\makebox(1.6667,11.6667){\SetFigFont{5}{6}{rm}.}}
\put(2376,-106){\line( 1, 4){ 75}}
\put(2481,-126){\line( 1, 4){ 75}}
\put(976,314){\vector( 1, 0){1350}}
\put(2851,314){\vector( 1, 0){1275}}
\put(4426,164){\vector( 0,-1){1500}}
\put(751, 89){\vector( 0,-1){1425}}
\put(976,-1486){\vector( 1, 0){1500}}
\put(2776,-1486){\vector( 1, 0){1500}}
\put(2626,-661){\vector( 0,-1){600}}
\put(976,164){\vector( 3,-1){1350}}
\put(2851,-361){\vector( 3, 1){1350}}
\put(601,239){\makebox(0,0)[lb]{\smash{\SetFigFont{10}{14.4}{rm}${\cal C}^{\boxtimes n}$}}}
\put(4201,239){\makebox(0,0)[lb]{\smash{\SetFigFont{10}{14.4}{rm}${\cal D}^{\boxtimes n}$}}}
\put(2401,239){\makebox(0,0)[lb]{\smash{\SetFigFont{10}{14.4}{rm}${\cal D}^{\boxtimes n}$}}}
\put(751,-1561){\makebox(0,0)[lb]{\smash{\SetFigFont{10}{14.4}{rm}${\cal C}$}}}
\put(2551,-1561){\makebox(0,0)[lb]{\smash{\SetFigFont{10}{14.4}{rm}${\cal C}$}}}
\put(4351,-1561){\makebox(0,0)[lb]{\smash{\SetFigFont{10}{14.4}{rm}${\cal D}$ ,}}}
\put(2401,-511){\makebox(0,0)[lb]{\smash{\SetFigFont{10}{14.4}{rm}${\cal C}^{\boxtimes n}$}}}
\put(271,-586){\makebox(0,0)[lb]{\smash{\SetFigFont{10}{14.4}{rm}$^n\otimes$}}}
\put(1396,464){\makebox(0,0)[lb]{\smash{\SetFigFont{10}{14.4}{rm}$F^n$}}}
\put(2851,419){\makebox(0,0)[lb]{\smash{\SetFigFont{10}{14.4}{rm}$F(A)\otimes -\boxtimes {\cal D}^{\boxtimes n-1}$}}}
\put(4681,-721){\makebox(0,0)[lb]{\smash{\SetFigFont{10}{14.4}{rm}$\otimes^n$}}}
\put(1291,-1786){\makebox(0,0)[lb]{\smash{\SetFigFont{10}{14.4}{rm}$A\otimes -$}}}
\put(3316,-1876){\makebox(0,0)[lb]{\smash{\SetFigFont{10}{14.4}{rm}$M$}}}
\put(3736,-886){\makebox(0,0)[lb]{\smash{\SetFigFont{10}{14.4}{rm}$\tilde{\phi}$}}}
\put(1696,-871){\makebox(0,0)[lb]{\smash{\SetFigFont{10}{14.4}{rm}\mbox{\boldmath $\alpha$\unboldmath}}}}
\put(2626,-46){\makebox(0,0)[lb]{\smash{\SetFigFont{10}{14.4}{rm}$\tilde{F}\boxtimes Id$}}}
\end{picture}

\end{center}

\noindent while the bottom path is induced by

\begin{center} 
\setlength{\unitlength}{3947sp}%
\begingroup\makeatletter\ifx\SetFigFont\undefined
\def\x#1#2#3#4#5#6#7\relax{\def\x{#1#2#3#4#5#6}}%
\expandafter\x\fmtname xxxxxx\relax \def\y{splain}%
\ifx\x\y   
\gdef\SetFigFont#1#2#3{%
  \ifnum #1<17\tiny\else \ifnum #1<20\small\else
  \ifnum #1<24\normalsize\else \ifnum #1<29\large\else
  \ifnum #1<34\Large\else \ifnum #1<41\LARGE\else
     \huge\fi\fi\fi\fi\fi\fi
  \csname #3\endcsname}%
\else
\gdef\SetFigFont#1#2#3{\begingroup
  \count@#1\relax \ifnum 25<\count@\count@25\fi
  \def\x{\endgroup\@setsize\SetFigFont{#2pt}}%
  \expandafter\x
    \csname \romannumeral\the\count@ pt\expandafter\endcsname
    \csname @\romannumeral\the\count@ pt\endcsname
  \csname #3\endcsname}%
\fi
\fi\endgroup
\begin{picture}(4410,2475)(271,-1876)
\thinlines
\multiput(1396,-121)(3.75000,-7.50000){21}{\makebox(1.6667,11.6667){\SetFigFont{5}{6}{rm}.}}
\multiput(1471,-271)(7.50000,3.75000){21}{\makebox(1.6667,11.6667){\SetFigFont{5}{6}{rm}.}}
\put(1436,-206){\line( 1, 4){ 75}}
\put(1541,-226){\line( 1, 4){ 75}}
\multiput(3391,-121)(3.75000,-7.50000){21}{\makebox(1.6667,11.6667){\SetFigFont{5}{6}{rm}.}}
\multiput(3466,-271)(7.50000,3.75000){21}{\makebox(1.6667,11.6667){\SetFigFont{5}{6}{rm}.}}
\put(3431,-206){\line( 1, 4){ 75}}
\put(3536,-226){\line( 1, 4){ 75}}
\multiput(2386,-1096)(3.75000,-7.50000){21}{\makebox(1.6667,11.6667){\SetFigFont{5}{6}{rm}.}}
\multiput(2461,-1246)(7.50000,3.75000){21}{\makebox(1.6667,11.6667){\SetFigFont{5}{6}{rm}.}}
\put(2426,-1181){\line( 1, 4){ 75}}
\put(2531,-1201){\line( 1, 4){ 75}}
\put(976,314){\vector( 1, 0){1350}}
\put(2851,314){\vector( 1, 0){1275}}
\put(4426,164){\vector( 0,-1){1500}}
\put(751, 89){\vector( 0,-1){1425}}
\put(976,-1486){\vector( 1, 0){1500}}
\put(2776,-1486){\vector( 1, 0){1500}}
\put(1006,-1336){\vector( 3, 1){1438.200}}
\put(2791,-751){\vector( 3,-1){1484.100}}
\put(2581,119){\vector( 0,-1){585}}
\put(601,239){\makebox(0,0)[lb]{\smash{\SetFigFont{10}{14.4}{rm}${\cal C}^{\boxtimes n}$}}}
\put(4201,239){\makebox(0,0)[lb]{\smash{\SetFigFont{10}{14.4}{rm}${\cal D}^{\boxtimes n}$}}}
\put(2401,239){\makebox(0,0)[lb]{\smash{\SetFigFont{10}{14.4}{rm}${\cal D}^{\boxtimes n}$}}}
\put(751,-1561){\makebox(0,0)[lb]{\smash{\SetFigFont{10}{14.4}{rm}${\cal C}$}}}
\put(2551,-1561){\makebox(0,0)[lb]{\smash{\SetFigFont{10}{14.4}{rm}${\cal C}$}}}
\put(4351,-1561){\makebox(0,0)[lb]{\smash{\SetFigFont{10}{14.4}{rm}${\cal D} ,$}}}
\put(271,-586){\makebox(0,0)[lb]{\smash{\SetFigFont{10}{14.4}{rm}$^n\otimes$}}}
\put(1396,464){\makebox(0,0)[lb]{\smash{\SetFigFont{10}{14.4}{rm}$F^n$}}}
\put(2851,419){\makebox(0,0)[lb]{\smash{\SetFigFont{10}{14.4}{rm}$F(A)\otimes - \boxtimes {\cal D}^{\boxtimes n-1}$}}}
\put(4681,-721){\makebox(0,0)[lb]{\smash{\SetFigFont{10}{14.4}{rm}$\otimes^n$}}}
\put(1291,-1786){\makebox(0,0)[lb]{\smash{\SetFigFont{10}{14.4}{rm}$A\otimes -$}}}
\put(3316,-1876){\makebox(0,0)[lb]{\smash{\SetFigFont{10}{14.4}{rm}$M$}}}
\put(1711,-136){\makebox(0,0)[lb]{\smash{\SetFigFont{10}{14.4}{rm}$\tilde{\phi}$}}}
\put(3706,-121){\makebox(0,0)[lb]{\smash{\SetFigFont{10}{14.4}{rm}
\boldmath$\alpha$\unboldmath}}}
\put(2536,-751){\makebox(0,0)[lb]{\smash{\SetFigFont{10}{14.4}{rm}${\cal D}$}}}
\put(2716,-1141){\makebox(0,0)[lb]{\smash{\SetFigFont{10}{14.4}{rm}$\mu$}}}
\put(2746,-241){\makebox(0,0)[lb]{\smash{\SetFigFont{10}{14.4}{rm}$\otimes^n$}}}
\put(1381,-991){\makebox(0,0)[lb]{\smash{\SetFigFont{10}{14.4}{rm}$M$}}}
\put(3496,-931){\makebox(0,0)[lb]{\smash{\SetFigFont{10}{14.4}{rm}$F(A) \otimes -$}}}
\end{picture}

\end{center}

\noindent where in each case \boldmath$\alpha \;$\unboldmath represents
the appropriate unique natural transformation given by 
Mac Lane's coherence theorem \cite{Mac.coh,CWM}.

The equality between the two natural fillers is given on objects
by an equation between

\begin{eqnarray*}
\lefteqn{^n\otimes(F(A)\otimes F(X_1),F(X_2),\ldots ,F(X_n))} \\
  &
\stackrel{^n\otimes(\tilde{F},Id,\ldots ,Id)}{\longrightarrow} &
^n\otimes(F(A\otimes X_1),F(X_2),\ldots ,F(X_n)) \\
 & \stackrel{\tilde{\phi}_{A\otimes X_1,X_2,\ldots ,X_n}}{\longrightarrow}
	&
M(\otimes^n(A\otimes X_1),X_2,\ldots ,X_n)) \\
 & \stackrel{\mbox{\boldmath $\alpha$\unboldmath}}{\longrightarrow}
	&
M(\otimes^{n+1}(A,X_1,\ldots X_n))
\end{eqnarray*}

\noindent and

\begin{eqnarray*}
\lefteqn{ ^n\otimes(F(A)\otimes F(X_1),F(X_2),\ldots ,F(X_n))} \\ 
&
\stackrel{\mbox{\boldmath $\alpha$\unboldmath}}{\longrightarrow}&
F(A)\otimes[^n\otimes(F(X_1),\ldots ,F(X_n))] \\
&\stackrel{F(A)\otimes \tilde{\phi}_{X_1,\ldots ,X_n}}{\longrightarrow}
	&
F(A)\otimes M(\otimes^n(X_1,\ldots ,X_n))\\
&\stackrel{\mu_{A,\otimes^n(X_1,\ldots ,X_n)}}{\longrightarrow}
	&
M(\otimes^{n+1}(A,X_1,\ldots X_n))
\end{eqnarray*}

Suppressing the \boldmath $\alpha \;$\unboldmath and parenthesizations
by invoking Mac Lane's coherence theorem, this becomes

\begin{center} 
\setlength{\unitlength}{3300sp}%
\begingroup\makeatletter\ifx\SetFigFont\undefined
\def\x#1#2#3#4#5#6#7\relax{\def\x{#1#2#3#4#5#6}}%
\expandafter\x\fmtname xxxxxx\relax \def\y{splain}%
\ifx\x\y   
\gdef\SetFigFont#1#2#3{%
  \ifnum #1<17\tiny\else \ifnum #1<20\small\else
  \ifnum #1<24\normalsize\else \ifnum #1<29\large\else
  \ifnum #1<34\Large\else \ifnum #1<41\LARGE\else
     \huge\fi\fi\fi\fi\fi\fi
  \csname #3\endcsname}%
\else
\gdef\SetFigFont#1#2#3{\begingroup
  \count@#1\relax \ifnum 25<\count@\count@25\fi
  \def\x{\endgroup\@setsize\SetFigFont{#2pt}}%
  \expandafter\x
    \csname \romannumeral\the\count@ pt\expandafter\endcsname
    \csname @\romannumeral\the\count@ pt\endcsname
  \csname #3\endcsname}%
\fi
\fi\endgroup
\begin{picture}(6150,2535)(151,-1861)
\thinlines
\put(1801, 89){\vector( 0,-1){1350}}
\put(3076,-1486){\vector( 1, 0){1950}}
\put(3076,314){\vector( 1, 0){1650}}
\put(6001,164){\vector( 0,-1){1425}}
\put(601,239){\makebox(0,0)[lb]{\smash{\SetFigFont{9}{14.4}{rm}$F(A)\otimes F(X_1)\otimes \ldots\otimes F(X_n)$}}}
\put(4801,239){\makebox(0,0)[lb]{\smash{\SetFigFont{9}{14.4}{rm}$F(A\otimes X_1)\otimes F(X_2) \ldots\otimes F(X_n)$}}}
\put(751,-1561){\makebox(0,0)[lb]{\smash{\SetFigFont{9}{14.4}{rm}$F(A)\otimes M(X_1\otimes \ldots\otimes X_n)$}}}
\put(5101,-1561){\makebox(0,0)[lb]{\smash{\SetFigFont{9}{14.4}{rm}$M(A\otimes X_1\otimes \ldots\otimes X_n)$}}}
\put(151,-811){\makebox(0,0)[lb]{\smash{\SetFigFont{9}{14.4}{rm}$F(A)\otimes \tilde{\phi}_{X_1\ldots X_n}$}}}
\put(6301,-661){\makebox(0,0)[lb]{\smash{\SetFigFont{9}{14.4}{rm}$\tilde{\phi}_{A\otimes X_1\ldots X_n}$}}}
\put(3451,539){\makebox(0,0)[lb]{\smash{\SetFigFont{9}{14.4}{rm}$\tilde{F}\otimes Id$}}}
\put(3751,-1861){\makebox(0,0)[lb]{\smash{\SetFigFont{9}{14.4}{rm}$\mu$}}}
\end{picture}

\end{center}

This diagram, then, is the condition required for a map from
$F(X_1)\otimes \ldots \otimes F(X_n)$ to $M(X_1\otimes \ldots \otimes X_n)$ to
induce a left-module map from $Lan_{\otimes^n}(^n\otimes (F^n))$ to $M$.  The
condition for inducing right-module and bimodule maps are entirely similar.

Consideration of the classical case of associative algebras shows that one
wants to begin with some projectivity assumptions about the underlying
objects.  In particular we will assume throughout the following discussion
that $F(I)$ is projective, that for all $A \in Ob({\cal C})$ the functors
$F(A)\otimes -$ and $- \otimes F(A)$ preserve epis, and moreover that
$F(I)\otimes -$ and $- \otimes F(I)$ have epi-preserving right adjoints.  
We will then proceed
by attempting to construct the requisite lifts to 
show that the various
$Lan_{\otimes^n}( ^n\otimes (F^n))$ are projective. Along
the way we will discover what other hypotheses will be needed.  
At each point, we will look
only for hypotheses satisfied in the case of ${\cal E} = k${\bf -v.s.} 
once the size restriction mentioned above has been placed on the 
source category.

We begin with the case $n = 2$.


Now observe that a map from 
$Lan_{\otimes^2} ( ^2\otimes (F^2)) = Lan_{\otimes}(F\otimes F)$ to
a functor $M:{\cal C}\rightarrow {\cal E}$ corresponds canonically to a 
natural transformation from 

\[ F(-) \otimes F(-):{\cal C}\boxtimes {\cal C}\rightarrow
{\cal E} \] 

\noindent to 

\[M(- \otimes -):{\cal C}\boxtimes {\cal C}\rightarrow
{\cal E} \,. \]

Thus a lift of a map $f:Lan_{\otimes^2} ( ^2\otimes (F^2))\rightarrow N$ 
along an epi $q:M\rightarrow N$ corresponds to a lift of a map 
$\phi:F(-) \otimes F(-)\rightarrow N(-\otimes -)$ along 
$q_{-\otimes -}$, as in Figure \ref{lift.at.2}.  Thus, in
particular, one must have a lift of the map 
$\phi_{I,I}:F(I)\otimes F(I)\rightarrow N(I\otimes I)$ along
$q_{I\otimes I}:M(I\otimes I)\rightarrow N(I\otimes I)$.  

The existence of the particular lift with $A = B = I$ 
follows trivially from the hypotheses that
$F(I)$ be projective and that $- \otimes F(I)$ (or $F(I)\otimes -$) admit
an epi-preserving right adjoint.  Now, recall that we want a lift in the
category of $F,F$-bimodules over $\cal C$, so that $M$ and $N$ are bimodules
and $q$ is a bimodule map.  We will now try to use the actions to extend
the particular case to the general case.

Tensoring the lifting diagram of in the upper corner of Figure
\ref{lift.at.2}  with $F(A)$ on the left and $F(B)$ on the
right, applying both actions, and adjoining obvious isomorphisms gives
the lower diagram of Figure \ref{lift.at.2}.

\begin{figure}[htb] \centering
\setlength{\unitlength}{3300sp}%
\begingroup\makeatletter\ifx\SetFigFont\undefined%
\gdef\SetFigFont#1#2#3#4#5{%
  \reset@font\fontsize{#1}{#2pt}%
  \fontfamily{#3}\fontseries{#4}\fontshape{#5}%
  \selectfont}%
\fi\endgroup%
\begin{picture}(5475,6084)(1201,-5590)
\thinlines
\put(3301,164){\vector( 0,-1){900}}
\put(1951,-886){\vector( 1, 0){900}}
\put(2626,-886){\vector( 1, 0){225}}
\put(2701,-886){\vector( 1, 0){150}}
\put(2851,164){\vector(-3,-2){1125}}
\put(2851,-886){\vector( 1, 0){0}}
\put(2626,-886){\vector( 1, 0){300}}
\put(6376,-1561){\vector( 0,-1){975}}
\put(4126,-2761){\vector( 1, 0){1125}}
\put(5476,-1561){\vector(-2,-1){1800}}
\put(2776,-3061){\vector( 0,-1){750}}
\put(6376,-2986){\vector( 0,-1){750}}
\put(3676,-4036){\vector( 1, 0){1650}}
\put(2776,-4336){\vector( 0,-1){675}}
\put(6301,-4261){\vector( 0,-1){750}}
\put(3301,-5236){\vector( 1, 0){2325}}
\put(6376,-1111){\vector( 0, 1){450}}
\put(6376,-286){\vector( 0, 1){300}}
\put(1201,-961){\makebox(0,0)[lb]{\smash{\SetFigFont{9}{14.4}{rm}$M(I \otimes I)$}}}
\put(3001,-961){\makebox(0,0)[lb]{\smash{\SetFigFont{9}{14.4}{rm}$N(I \otimes I)$}}}
\put(2926,314){\makebox(0,0)[lb]{\smash{\SetFigFont{9}{14.4}{rm}$F(I) \otimes F(I)$}}}
\put(3526,-286){\makebox(0,0)[lb]{\smash{\SetFigFont{9}{14.4}{rm}$\phi_{I,I}$}}}
\put(1426,-136){\makebox(0,0)[lb]{\smash{\SetFigFont{9}{14.4}{rm}$\psi_{I,I}$}}}
\put(5401,-1411){\makebox(0,0)[lb]{\smash{\SetFigFont{9}{14.4}{rm}$F(A) \otimes F(I) \otimes F(I) \otimes F(B)$}}}
\put(5476,-2836){\makebox(0,0)[lb]{\smash{\SetFigFont{9}{14.4}{rm}$F(A) \otimes N(I \otimes I) \otimes F(B)$}}}
\put(2176,-2836){\makebox(0,0)[lb]{\smash{\SetFigFont{9}{14.4}{rm}$F(A) \otimes M(I \otimes I) \otimes F(B)$}}}
\put(5551,-586){\makebox(0,0)[lb]{\smash{\SetFigFont{9}{14.4}{rm}$F(A \otimes I) \otimes F(I \otimes B)$}}}
\put(5926, 89){\makebox(0,0)[lb]{\smash{\SetFigFont{9}{14.4}{rm}$F(A) \otimes F(B)$}}}
\put(5626,-4036){\makebox(0,0)[lb]{\smash{\SetFigFont{9}{14.4}{rm}$N(A \otimes I \otimes I \otimes B)$}}}
\put(2176,-4111){\makebox(0,0)[lb]{\smash{\SetFigFont{9}{14.4}{rm}$M(A \otimes I \otimes I \otimes B)$}}}
\put(2326,-5311){\makebox(0,0)[lb]{\smash{\SetFigFont{9}{14.4}{rm}$M(A \otimes B)$}}}
\put(5851,-5311){\makebox(0,0)[lb]{\smash{\SetFigFont{9}{14.4}{rm}$N(A \otimes B)$}}}
\put(6676,-2086){\makebox(0,0)[lb]{\smash{\SetFigFont{9}{14.4}{rm}$F(A) \otimes \phi_{I,I} \otimes F(B)$}}}
\put(2551,-2086){\makebox(0,0)[lb]{\smash{\SetFigFont{9}{14.4}{rm}$F(A) \otimes \psi_{I,I} \otimes F(B)$}}}
\put(6601,-961){\makebox(0,0)[lb]{\smash{\SetFigFont{9}{14.4}{rm}$\tilde{F}_{A,I} \otimes \tilde{F}_{I,B}$}}}
\put(6676,-211){\makebox(0,0)[lb]{\smash{\SetFigFont{9}{14.4}{rm}$\rho \otimes \lambda$}}}
\put(2101,-1186){\makebox(0,0)[lb]{\smash{\SetFigFont{9}{14.4}{rm}$q_{I \otimes I}$}}}
\put(3901,-3136){\makebox(0,0)[lb]{\smash{\SetFigFont{9}{14.4}{rm}$F(A) \otimes q_{I \otimes I} \otimes F(B)$}}}
\put(3826,-4336){\makebox(0,0)[lb]{\smash{\SetFigFont{9}{14.4}{rm}$q_{A \otimes I \otimes I \otimes B}$}}}
\put(3826,-5536){\makebox(0,0)[lb]{\smash{\SetFigFont{9}{14.4}{rm}$q_{A \otimes B}$}}}
\put(6676,-3436){\makebox(0,0)[lb]{\smash{\SetFigFont{9}{14.4}{rm}$\nu_l \otimes F(B) \nu_r$}}}
\put(1276,-3511){\makebox(0,0)[lb]{\smash{\SetFigFont{9}{14.4}{rm}$\mu_l \otimes F(B) \mu_r$}}}
\put(1576,-4711){\makebox(0,0)[lb]{\smash{\SetFigFont{9}{14.4}{rm}$M(\rho \otimes \lambda)$}}}
\put(6601,-4711){\makebox(0,0)[lb]{\smash{\SetFigFont{9}{14.4}{rm}$N(\rho \otimes \lambda)$}}}
\end{picture}
\caption{Lifting for $n = 2$ \label{lift.at.2}}
\end{figure}

This, however, is not quite what we want:  the map between

\[ F(A\otimes I)\otimes F(I\otimes B) \] 

\noindent and 

\[ F(A)\otimes F(I)\otimes F(I)\otimes F(B) \]

\noindent  runs in the wrong direction.  If $F$ is a strong
monoidal functor this is easy to remedy, since we can invert the map.
Taking this as a hypothesis would, however, exclude the classical
case of associative algebras.  We therefore assume the weaker hypothesis
that $\tilde{F}_{I,A}$ and $\tilde{F}_{A,I}$  are epi for all $A$ and
admits a natural splitting $\check{F}_{I,A}$ (resp. $\check{F}_{I,A}$).

Adjoining the splitting to the diagram almost completes the construction.
We still must verify that the resulting composite from $F(A)\otimes F(B)$
to $N(A\otimes B)$ is ${\phi}_{A,B}$.  Now, consider the diagram

\begin{center}
\setlength{\unitlength}{3300sp}%
\begingroup\makeatletter\ifx\SetFigFont\undefined%
\gdef\SetFigFont#1#2#3#4#5{%
  \reset@font\fontsize{#1}{#2pt}%
  \fontfamily{#3}\fontseries{#4}\fontshape{#5}%
  \selectfont}%
\fi\endgroup%
\begin{picture}(7287,4455)(151,-4036)
\thinlines
\put(2901, 89){\vector( 1, 0){1150}} 
\put(4051,-61){\vector(-1, 0){1150}} 
\put(7126, 14){\vector(-1, 0){1100}} 
\put(1576,-211){\vector( 0,-1){750}}
\put(1576,-1411){\vector( 0,-1){750}}
\put(1576,-2686){\vector( 0,-1){1050}}
\put(4501,-361){\vector(-4,-3){2400}}
\put(7426,-286){\vector(-3,-2){5175}}
\put(4201,-61){\makebox(0,0)[lb]{\smash{\SetFigFont{9}{14.4}{rm}$F(A \otimes I) \otimes F(B \otimes I)$}}}
\put(7201,-61){\makebox(0,0)[lb]{\smash{\SetFigFont{9}{14.4}{rm}$F(A) \otimes F(B)$}}}
\put(601,-61){\makebox(0,0)[lb]{\smash{\SetFigFont{9}{14.4}{rm}$F(A) \otimes F(I) \otimes F(I) \otimes F(B)$}}}
\put(1201,-4036){\makebox(0,0)[lb]{\smash{\SetFigFont{9}{14.4}{rm}$N(A \otimes B)$}}}
\put(901,-2461){\makebox(0,0)[lb]{\smash{\SetFigFont{9}{14.4}{rm}$N(A \otimes I \otimes I \otimes B)$}}}
\put(676,-1261){\makebox(0,0)[lb]{\smash{\SetFigFont{9}{14.4}{rm}$F(A) \otimes N(I \otimes I) \otimes F(B)$}}}
\put(4876,-2611){\makebox(0,0)[lb]{\smash{\SetFigFont{9}{14.4}{rm}$\phi_{A,B}$}}}
\put(3301,-1636){\makebox(0,0)[lb]{\smash{\SetFigFont{9}{14.4}{rm}$\phi_{A \otimes I,B \otimes I}$}}}
\put(226,-1861){\makebox(0,0)[lb]{\smash{\SetFigFont{9}{14.4}{rm}$\nu_l \otimes F(B) \nu_r$}}}
\put(1651,-661){\makebox(0,0)[lb]{\smash{\SetFigFont{9}{14.4}{rm}$F(A) \otimes \phi_{I,I} \otimes F(B)$}}}
\put(6001,239){\makebox(0,0)[lb]{\smash{\SetFigFont{9}{14.4}{rm}$F(\rho^{-1}) \otimes F(\lambda^{-1})$}}}
\put(151,-3361){\makebox(0,0)[lb]{\smash{\SetFigFont{9}{14.4}{rm}$N(\rho \otimes \lambda)$}}}
\end{picture}
\end{center}

\noindent  The bottom square commutes by naturality of $\phi$ and the 
invertibility of $F(\rho)$ and $F(\lambda)$, while the
square beginning at $F(A)\otimes F(I)\otimes F(I)\otimes F(B)$ and 
ending at $N(A\otimes I\otimes I \otimes B)$ commutes
since this is the $n = 2$ case of the condition above which $\phi$
satisfies since it induces an $F,F$-bimodule map.  From this it follows
that the square beginning at $F(A\otimes I)\otimes F(I\otimes B)$
and ending at $N(A\otimes I\otimes I \otimes B)$ commutes
since $\tilde{F}(\check{F}) = Id$.

For the general case, we will proceed by induction. 

Fix $k$.  Now, suppose for all epimorphic
$F,F$-bimodule maps
$q:M\rightarrow N$, and all natural transformations
$\psi:F(X_1)\otimes \ldots F(X_k) \rightarrow N(X_1\otimes \ldots
\otimes X_k)$ satisfying the condition given above by which the induced
map from the Kan extension to $N$ is an $F,F$-bimodule map, we have
a lift $\hat{\psi}:F(X_1)\otimes \ldots F(X_k) \rightarrow M(X_1\otimes
\ldots \otimes X_k)$.

Now, tensor the lifting diagram on the right by $F(X_{k+1})$, and apply
the actions on the bimodules $M$ and $N$.  We claim that the resulting 
diagram then gives the map from $F(X_1)\otimes \ldots \otimes F(X_{k+1})$
to $M(X_1\otimes \ldots \otimes X_{k+1})$ which induces the required lift.

Plainly it induces a lift in the ambient functor category, so we only 
need show that it satisfies the conditions to induce left and right module
maps.  The former follows from the hypothesis that the original lift
satisfied the same condition.  The latter follows from the coherence 
property of the action on $M$. 

Fitting together all of the foregoing discussion we have shown:

\begin{thm}
If $F:{\cal C}\rightarrow {\cal D}$ is a monoidal functor with 
a small source, targetted
in a cocomplete abelian category $\cal D$ such that
for all $X\in Ob({\cal C})$ $F(X) \otimes -$ and $- \otimes F(X)$ are
exact and cocontinuous, $F(I)$ is projective, $F(I) \otimes -$ and $- \otimes
F(I)$ have epi-preserving right adjoints, and $\tilde{F}_{I,A}$ and 
$\tilde{F}_{A,I}$ are split epis with splitting natural in $A$, then
for all $n \geq 2$, $Lan_{\otimes^n}( ^n\otimes (F^n))$ is a projective
$F,F$-bimodule.
\end{thm}

First, let us observe that the rather baroque seeming technical hypotheses
are satisfied in two very natural cases:  strong monoidal functors 
targetted at a category of vectorspaces and algebras over a field, regarded as 
monoidal functors.

One annoying feature of this theorem is the size restriction on the
source, which makes the theorem inapplicable to identity functors
on large categories.  This size restriction can be relaxed somewhat by 
considering the construction of left Kan extensions as colimits.
In particular, recall from \cite{CWM} that the left Kan extension $Lan_K(T)$,
for $K:{\cal M}\rightarrow {\cal C}$ and $T:{\cal M}\rightarrow {\cal A}$
is given on objects by 
$Lan_K(T)(c) = {\rm colim}(K\downarrow c)\stackrel{P}{\rightarrow} 
{\cal M}\stackrel{T}{\rightarrow} {\cal A}$, where $P$ is the obvious 
projection functor from the comma category to $\cal M$.  It thus suffices 
for all of the comma categories $\otimes^n \downarrow c$ arising in the
Kan extensions used to admit small final subcategories and for $\cal D$
to admit colimits over all diagrams of no greater than the supremum of
the sizes of these final subcategories.

It is easy to see that any semisimple category with a small set of simple
objects satisfies the small-final subcategory condition on all of the
relevant comma categories.

Thus we state the seemingly more technical

\begin{thm} \label{strong.proj}
Suppose $F:{\cal C}\rightarrow {\cal D}$ is a monoidal functor such that
for all $n\in {\Bbb N}$ and all $c \in Ob({\cal C})$, the 
comma category $\otimes^n \downarrow c$ admits a small final subcategory 
$J_{n,c}$,
and the target
is an abelian category $\cal D$, which admits all colimits of diagrams with
cardinality less than $\sup |Arr(J_{n,c})|$ or all small diagrams if no
supremum of the cardinalities exists, and such that
for all $X\in Ob({\cal C})$ $F(X) \otimes -$ and $- \otimes F(X)$ are
exact and cocontinuous, $F(I)$ is projective, $F(I) \otimes -$ and $- \otimes
F(I)$ have epi-preserving right adjoints, and $\tilde{F}_{I,A}$ and 
$\tilde{F}_{A,I}$ are split epis with splitting natural in $A$.  In this case
 $Lan_{\otimes^n}( ^n\otimes (F^n))$ is a projective
$F,F$ -bimodule for all $n \geq 2$.
\end{thm}

Despite the technical nature of the theorem, it now applies to identity 
functors on many abelian monoidal
categories of interest.

The point of this result is, of course, to show that in many cases the 
cohomology of a monoidal functor with coefficients in a bimodule, and
in particular the deformation cohomology of a monoidal functor or 
monoidal category, is given by right derived functors.

To show this, we must see that the complex $X^\bullet(F,M)$ is
obtained by applying $Nat[-,M]$ to a projective resolution whose objects
are the $Lan_{\otimes^\bullet}(^\bullet\otimes (F^\bullet))$.

First, we show that the coboundary maps are induced by maps between the
Kan extensions:

\begin{thm}
The coboundary of the complex $X^\bullet(F,M)$ is induced by a map
of $F,F$-bimodules

\[ \partial:Lan_{\otimes^{n+1}}(^{n+1}\otimes(F^{n+1}))\rightarrow 
Lan_{\otimes^{n}}(^n\otimes(F^n)). \]
\end{thm}
 
\noindent {\bf proof:} The key is to consider the universal property 
defining $Lan_{\otimes^{n}}(^n\otimes(F^n))$, and to find natural 
transformations filling the square

\begin{center}
\setlength{\unitlength}{3947sp}%
\begingroup\makeatletter\ifx\SetFigFont\undefined
\def\x#1#2#3#4#5#6#7\relax{\def\x{#1#2#3#4#5#6}}%
\expandafter\x\fmtname xxxxxx\relax \def\y{splain}%
\ifx\x\y   
\gdef\SetFigFont#1#2#3{%
  \ifnum #1<17\tiny\else \ifnum #1<20\small\else
  \ifnum #1<24\normalsize\else \ifnum #1<29\large\else
  \ifnum #1<34\Large\else \ifnum #1<41\LARGE\else
     \huge\fi\fi\fi\fi\fi\fi
  \csname #3\endcsname}%
\else
\gdef\SetFigFont#1#2#3{\begingroup
  \count@#1\relax \ifnum 25<\count@\count@25\fi
  \def\x{\endgroup\@setsize\SetFigFont{#2pt}}%
  \expandafter\x
    \csname \romannumeral\the\count@ pt\expandafter\endcsname
    \csname @\romannumeral\the\count@ pt\endcsname
  \csname #3\endcsname}%
\fi
\fi\endgroup
\begin{picture}(3375,2535)(676,-2461)
\thinlines
\put(1276,-511){\vector( 0,-1){1350}}
\put(1501,-2086){\vector( 1, 0){2100}}
\put(1951,-286){\vector( 1, 0){1500}}
\put(3751,-436){\vector( 0,-1){1425}}
\put(2251,-1036){\line( 0,-1){225}}
\put(2251,-1261){\line( 1, 0){225}}
\put(2251,-1111){\line( 1, 1){225}}
\put(2386,-1246){\line( 1, 1){225}}
\put(1201,-361){\makebox(0,0)[lb]{\smash
{\SetFigFont{10}{14.4}{rm}${\cal C}^{\boxtimes n+1}$}}}
\put(3601,-361){\makebox(0,0)[lb]{\smash{\SetFigFont{10}{14.4}{rm}${\cal D}^{\boxtimes n+1}$}}}
\put(1276,-2161){\makebox(0,0)[lb]{\smash{\SetFigFont{10}{14.4}{rm}$\cal C$}}}
\put(3751,-2161){\makebox(0,0)[lb]{\smash{\SetFigFont{10}{14.4}{rm}$\cal D$}}}
\put(676,-1336){\makebox(0,0)[lb]{\smash
{\SetFigFont{10}{14.4}{rm}$\otimes^{n+1}$}}}
\put(4051,-1336){\makebox(0,0)[lb]{\smash
{\SetFigFont{10}{14.4}{rm}$^{n+1}\otimes$}}}
\put(2401,-61){\makebox(0,0)[lb]{\smash{\SetFigFont{10}{14.4}{rm}$F^{n+1}$}}}
\put(1786,-2461){\makebox(0,0)[lb]{\smash{\SetFigFont{10}{14.4}{rm}$Lan_{\otimes^n}(^n \otimes (F^n))$}}}
\end{picture}

\end{center}

\noindent in such a way that their composition with any 
$\phi:Lan_{\otimes^{n}}(^n\otimes(F^n))\rightarrow M$ is induced by each 
term of $\delta(\tilde{\phi})$, where $\tilde{\phi}: ^n\otimes (F^n)\Rightarrow
M(\otimes^n)$ is the corresponding map.

All but the first and last terms have fillers of the form

\begin{center}
\setlength{\unitlength}{3947sp}%
\begingroup\makeatletter\ifx\SetFigFont\undefined
\def\x#1#2#3#4#5#6#7\relax{\def\x{#1#2#3#4#5#6}}%
\expandafter\x\fmtname xxxxxx\relax \def\y{splain}%
\ifx\x\y   
\gdef\SetFigFont#1#2#3{%
  \ifnum #1<17\tiny\else \ifnum #1<20\small\else
  \ifnum #1<24\normalsize\else \ifnum #1<29\large\else
  \ifnum #1<34\Large\else \ifnum #1<41\LARGE\else
     \huge\fi\fi\fi\fi\fi\fi
  \csname #3\endcsname}%
\else
\gdef\SetFigFont#1#2#3{\begingroup
  \count@#1\relax \ifnum 25<\count@\count@25\fi
  \def\x{\endgroup\@setsize\SetFigFont{#2pt}}%
  \expandafter\x
    \csname \romannumeral\the\count@ pt\expandafter\endcsname
    \csname @\romannumeral\the\count@ pt\endcsname
  \csname #3\endcsname}%
\fi
\fi\endgroup
\begin{picture}(4350,3210)(676,-3196)
\thinlines
\put(3006,-791){\line( 0,-1){120}}
\put(3011,-911){\line( 1, 0){ 95}}
\multiput(3016,-831)(7.07692,5.30769){14}{\makebox(1.6667,11.6667){\SetFigFont{5}{6}{rm}.}}
\multiput(3066,-901)(7.07692,5.30769){14}{\makebox(1.6667,11.6667){\SetFigFont{5}{6}{rm}.}}
\put(1481,-1586){\line( 0,-1){120}}
\put(1486,-1706){\line( 1, 0){ 95}}
\multiput(1491,-1626)(7.07692,5.30769){14}{\makebox(1.6667,11.6667){\SetFigFont{5}{6}{rm}.}}
\multiput(1541,-1696)(7.07692,5.30769){14}{\makebox(1.6667,11.6667){\SetFigFont{5}{6}{rm}.}}
\put(3006,-1961){\line( 0,-1){120}}
\put(3011,-2081){\line( 1, 0){ 95}}
\multiput(3016,-2001)(7.07692,5.30769){14}{\makebox(1.6667,11.6667){\SetFigFont{5}{6}{rm}.}}
\multiput(3066,-2071)(7.07692,5.30769){14}{\makebox(1.6667,11.6667){\SetFigFont{5}{6}{rm}.}}
\put(4291,-1351){\line( 0,-1){120}}
\put(4296,-1471){\line( 1, 0){ 95}}
\multiput(4301,-1391)(7.07692,5.30769){14}{\makebox(1.6667,11.6667){\SetFigFont{5}{6}{rm}.}}
\multiput(4351,-1461)(7.07692,5.30769){14}{\makebox(1.6667,11.6667){\SetFigFont{5}{6}{rm}.}}
\put(1351,-526){\vector( 0,-1){2085}}
\put(2191,-1561){\vector(-2,-3){720}}
\put(1621,-526){\vector( 2,-3){430}}
\put(2581,-1351){\vector( 1, 0){960}}
\put(1576,-2866){\vector( 1, 0){2925}}
\put(1906,-286){\vector( 1, 0){2490}}
\put(4486,-496){\vector(-2,-3){450}}
\put(4711,-541){\vector( 0,-1){1980}}
\put(3961,-1606){\vector( 2,-3){620}}
\put(1201,-361){\makebox(0,0)[lb]{\smash{\SetFigFont{10}{14.4}{rm}${\cal C}^{\boxtimes n+1}$}}}
\put(676,-1336){\makebox(0,0)[lb]{\smash{\SetFigFont{10}{14.4}{rm}$\otimes^{n+1}$}}}
\put(1336,-2956){\makebox(0,0)[lb]{\smash{\SetFigFont{10}{14.4}{rm}$\cal C$}}}
\put(4576,-361){\makebox(0,0)[lb]{\smash{\SetFigFont{10}{14.4}{rm}${\cal D}^{\boxtimes n+1}$}}}
\put(5026,-1576){\makebox(0,0)[lb]{\smash{\SetFigFont{10}{14.4}{rm}$^{n+1} \otimes$}}}
\put(2311,-3196){\makebox(0,0)[lb]{\smash{\SetFigFont{10}{14.4}{rm}$Lan_{\otimes^n}(^n\otimes (F^n))$}}}
\put(3076,-121){\makebox(0,0)[lb]{\smash{\SetFigFont{10}{14.4}{rm}$F^{n+1}$}}}
\put(2056,-1426){\makebox(0,0)[lb]{\smash{\SetFigFont{10}{14.4}{rm}${\cal C}^{\boxtimes n}$}}}
\put(4651,-2926){\makebox(0,0)[lb]{\smash{\SetFigFont{10}{14.4}{rm}$\cal D$}}}
\put(3691,-1441){\makebox(0,0)[lb]{\smash{\SetFigFont{10}{14.4}{rm}${\cal D}^{\boxtimes n}$}}}
\put(1606,-1846){\makebox(0,0)[lb]{\smash{\SetFigFont{10}{14.4}{rm}\boldmath$\alpha$\unboldmath}}}
\put(4261,-1666){\makebox(0,0)[lb]{\smash{\SetFigFont{10}{14.4}{rm}\boldmath$\alpha$\unboldmath}}}
\put(3136,-2236){\makebox(0,0)[lb]{\smash{\SetFigFont{10}{14.4}{rm}$can$}}}
\put(3181,-1006){\makebox(0,0)[lb]{\smash{\SetFigFont{10}{14.4}{rm}$. . .\tilde{F}. . .$}}}
\put(1831,-676){\makebox(0,0)[lb]{\smash{\SetFigFont{10}{14.4}{rm}$. . . \otimes . . .$}}}
\put(3646,-661){\makebox(0,0)[lb]{\smash{\SetFigFont{10}{14.4}{rm}$. . . \otimes . . .$}}}
\put(1861,-2431){\makebox(0,0)[lb]{\smash{\SetFigFont{10}{14.4}{rm}$\otimes^n$}}}
\put(3916,-2266){\makebox(0,0)[lb]{\smash{\SetFigFont{10}{14.4}{rm}$^n\otimes $}}}
\end{picture}

\end{center}

The first has the filler

\begin{center}
\setlength{\unitlength}{3947sp}%
\begingroup\makeatletter\ifx\SetFigFont\undefined
\def\x#1#2#3#4#5#6#7\relax{\def\x{#1#2#3#4#5#6}}%
\expandafter\x\fmtname xxxxxx\relax \def\y{splain}%
\ifx\x\y   
\gdef\SetFigFont#1#2#3{%
  \ifnum #1<17\tiny\else \ifnum #1<20\small\else
  \ifnum #1<24\normalsize\else \ifnum #1<29\large\else
  \ifnum #1<34\Large\else \ifnum #1<41\LARGE\else
     \huge\fi\fi\fi\fi\fi\fi
  \csname #3\endcsname}%
\else
\gdef\SetFigFont#1#2#3{\begingroup
  \count@#1\relax \ifnum 25<\count@\count@25\fi
  \def\x{\endgroup\@setsize\SetFigFont{#2pt}}%
  \expandafter\x
    \csname \romannumeral\the\count@ pt\expandafter\endcsname
    \csname @\romannumeral\the\count@ pt\endcsname
  \csname #3\endcsname}%
\fi
\fi\endgroup
\begin{picture}(4350,3210)(676,-3196)
\thinlines
\put(3006,-791){\line( 0,-1){120}}
\put(3011,-911){\line( 1, 0){ 95}}
\multiput(3016,-831)(7.07692,5.30769){14}{\makebox(1.6667,11.6667){\SetFigFont{5}{6}{rm}.}}
\multiput(3066,-901)(7.07692,5.30769){14}{\makebox(1.6667,11.6667){\SetFigFont{5}{6}{rm}.}}
\put(1481,-1586){\line( 0,-1){120}}
\put(1486,-1706){\line( 1, 0){ 95}}
\multiput(1491,-1626)(7.07692,5.30769){14}{\makebox(1.6667,11.6667){\SetFigFont{5}{6}{rm}.}}
\multiput(1541,-1696)(7.07692,5.30769){14}{\makebox(1.6667,11.6667){\SetFigFont{5}{6}{rm}.}}
\put(3006,-1961){\line( 0,-1){120}}
\put(3011,-2081){\line( 1, 0){ 95}}
\multiput(3016,-2001)(7.07692,5.30769){14}{\makebox(1.6667,11.6667){\SetFigFont{5}{6}{rm}.}}
\multiput(3066,-2071)(7.07692,5.30769){14}{\makebox(1.6667,11.6667){\SetFigFont{5}{6}{rm}.}}
\put(4291,-1351){\line( 0,-1){120}}
\put(4296,-1471){\line( 1, 0){ 95}}
\multiput(4301,-1391)(7.07692,5.30769){14}{\makebox(1.6667,11.6667){\SetFigFont{5}{6}{rm}.}}
\multiput(4351,-1461)(7.07692,5.30769){14}{\makebox(1.6667,11.6667){\SetFigFont{5}{6}{rm}.}}
\put(1351,-526){\vector( 0,-1){2085}}
\put(2191,-1561){\vector(-2,-3){720}}
\put(1621,-526){\vector( 2,-3){430}}
\put(2581,-1351){\vector( 1, 0){960}}
\put(1576,-2866){\vector( 1, 0){2925}}
\put(1906,-286){\vector( 1, 0){2490}}
\put(4486,-496){\vector(-2,-3){450}}
\put(4711,-541){\vector( 0,-1){1980}}
\put(3961,-1606){\vector( 2,-3){620}}
\put(1201,-361){\makebox(0,0)[lb]{\smash{\SetFigFont{10}{14.4}{rm}${\cal C}^{\boxtimes n+1}$}}}
\put(676,-1336){\makebox(0,0)[lb]{\smash{\SetFigFont{10}{14.4}{rm}$\otimes^{n+1}$}}}
\put(1336,-2956){\makebox(0,0)[lb]{\smash{\SetFigFont{10}{14.4}{rm}$\cal C$}}}
\put(4576,-361){\makebox(0,0)[lb]{\smash{\SetFigFont{10}{14.4}{rm}${\cal D}^{\boxtimes n+1}$}}}
\put(5026,-1576){\makebox(0,0)[lb]{\smash{\SetFigFont{10}{14.4}{rm}$^{n+1} \otimes$}}}
\put(2311,-3196){\makebox(0,0)[lb]{\smash{\SetFigFont{10}{14.4}{rm}$Lan_{\otimes^n}(^n\otimes (F^n))$}}}
\put(3076,-121){\makebox(0,0)[lb]{\smash{\SetFigFont{10}{14.4}{rm}$F^{n+1}$}}}
\put(4651,-2926){\makebox(0,0)[lb]{\smash{\SetFigFont{10}{14.4}{rm}$\cal D$}}}
\put(1606,-1846){\makebox(0,0)[lb]{\smash{\SetFigFont{10}{14.4}{rm}\boldmath$\alpha$\unboldmath}}}
\put(4261,-1666){\makebox(0,0)[lb]{\smash{\SetFigFont{10}{14.4}{rm}\boldmath$\alpha$\unboldmath}}}
\put(1861,-2431){\makebox(0,0)[lb]{\smash{\SetFigFont{10}{14.4}{rm}$\otimes $}}}
\put(3916,-2266){\makebox(0,0)[lb]{\smash{\SetFigFont{10}{14.4}{rm}$ \otimes$}}}
\put(2056,-1426){\makebox(0,0)[lb]{\smash{\SetFigFont{10}{14.4}{rm}${\cal C}^{\boxtimes 2}$}}}
\put(3691,-1441){\makebox(0,0)[lb]{\smash{\SetFigFont{10}{14.4}{rm}${\cal D}^{\boxtimes 2}$}}}
\put(2626,-1261){\makebox(0,0)[lb]{\smash{\SetFigFont{10}{14.4}{rm}$F \boxtimes Lan$}}}
\put(2476,-661){\makebox(0,0)[lb]{\smash{\SetFigFont{10}{14.4}{rm}$Id_F \boxtimes can$}}}
\put(3151,-2236){\makebox(0,0)[lb]{\smash{\SetFigFont{10}{14.4}{rm}$\mu_l$}}}
\end{picture}

\end{center}

\noindent while the last has a similar one with left and right reversed.
In both of the diagrams above, each \boldmath$\alpha$\unboldmath, denotes
an appropriate map given by Mac Lane's coherence theorem. 

It is immediate that the composition of these with $\phi$ has the
desired property.  Therefore we can let the boundary map be the map between
Kan extensions induced by their alternating sum. $\Box$

Finally, we must show that the sequence of these maps actually forms a
resolution in the category of $F,F$-bimodules over $\cal C$.  

To do this, we use

\begin{lemma}
In any abelian category $\cal A$, given a complex 

\[ A\stackrel{\beta}{\longrightarrow} B \stackrel{\gamma}{\longrightarrow} C\]
 
\noindent such that for all $M \in Ob({\cal A})$

\[ Hom(C,M) \stackrel{Hom(\gamma,M)}{\longrightarrow} Hom(B,M)
\stackrel{Hom(\beta,M)}{\longrightarrow} Hom(A,M) \]

\noindent is exact in {\bf Ab}, the original complex
is itself exact in $\cal A$.
\end{lemma}

\noindent {\bf proof:} Suppose the original sequence is not exact, 
that is, $Im(\beta)$ is a
proper subobject of $ker(\gamma)$.  Now, let $M = coker(\beta)$.  The
canonical map $B\rightarrow M$ is in the kernel of $Hom(\beta, M)$, but 
is not in the image of $Hom(\gamma,M)$. All such maps have kernels 
containing $ker(\gamma)$. $\Box$

Now, we apply this lemma to the complex 
$Lan_{\otimes^\bullet}(^\bullet \otimes(F^\bullet)), \partial^\bullet$
to obtain

\begin{thm}
Under the hypotheses of Theorem \ref{strong.proj}, the complex

\[ Lan_{\otimes^\bullet}(^\bullet \otimes(F^\bullet)), \partial^\bullet \]

is a projective resolution of $F$ as an $F,F$-bimodule.
\end{thm}

\noindent {\bf proof:} 
Now by the lemma, it suffices to show that the complex

\[ Hom_{F,F-{\bf bimod}}(Lan_{\otimes^\bullet}(^\bullet \otimes(F^\bullet)),M),
 Hom_{F,F-{\bf bimod}}(\partial^\bullet,M) \]

 For general $n$, observe that by the bimodule coherence 
conditions, the first and second and the last and penultimate terms in
the expression for the coboundary cancel in pairs.  Using this, it is
easy to see that any $n+1$-cocycle $\phi_{X_1,X_2,\ldots ,X_{n+1}}$ is
the coboundary of $\phi_{X_1,I,X_2,\ldots ,X_{n+1}}$.

The cases of $1$, $2$ and $3$-cocycles must be handled separately:
For $1$-cocycles, use the coherence condition to cancel the first two
terms of the coboundary.  The remaining term is $0$ if and only if $\phi$ is
$0$, since tensoring with $F(X)$ is exact.

For $2$-cocycles, first observe that any $2$-cochain is a cocycle (by the
pairwise cancellation noted above).  But application of the coherence 
conditions to cancel two terms and rewrite the remaining term shows that any
$\phi_{X,Y}$ is the coboundary of the 1-cochain $\phi_{X,I}$.

For $3$-cocycles, the pairwise cancellation of terms in the coboundary
formula reduces the cocycle condition to 
$\lceil \phi_{X,Y\otimes Z, W}\rceil = 0$
for all $X,Y,Z,W \in Ob({\cal C})$.  Specializing to $Y = I$ then shows
that $\phi = 0$.  But this taken with the previous observation that all
$2$-cochains are cocyles completes the proof of the exactness of the
sequence of $Hom(-,M)$'s, and by the previous lemma, of the theorem. $\Box$

As a consequence we now have

\begin{thm}

If $H^\bullet(F,M)$ denotes the cohomology of $(X^\bullet(F,M),\delta)$,
then

\[ H\bullet(F,M) = {\bf R}Nat[-,M](F) . \]

\end{thm}

Although we have no cause to pursue the matter, this last result allows us
to generalize the deformation cohomology of a monoidal functor to a 
special cohomology of an $F,F$-bimodule $N$ with coefficients in
another $F,F$-bimodule $M$, given by right derived functors of $Nat[-,M]$.
Notice that this is not simply {\bf Ext} for the abelian category of
$F,F$-bimodules, since it involves natural transformations which are not
bimodule maps.

\clearpage 
  \section{Conclusions}
\vspace*{1cm}

The preceding results do not, in and of themselves, provide much help in
calculating deformations of monoidal categories, monoidal functors or
braided monoidal categories.  They do, however, move the subject into the
realm of classical homological algebra. 

As such, it may be hoped that they will further the development of 
categorical deformation theory, both as a subject in its own right, and 
in its applications to the theory of Vassiliev invariants 
(cf. \cite{Vas,Y.def}) 
and
to the quest of non-trivial Hopf categories (cf. \cite{CF,CY.def}) which 
was its
original motivation.

\clearpage

\begin{flushleft}
\bibliography{Book}
\end{flushleft}

\end{document}